\input ./style/arxiv-general.cfg
\documentclass[aop,MSNbibl,seceqn,dvips]{arximspdf}
\makeatletter
   \@ifpackageloaded{graphicx}{}{\usepackage{graphicx}}
\makeatother

\doi{10.1214/15-AOP1030}
\volume{45}
\issue{1}
\pubyear{2017}
\firstpage{4}
\lastpage{55}
\docsubty{FLA}

\makeatletter
\newcommand{\rrvert}{\vert}
\newcommand{\llvert}{\vert}
\newcommand{\eqref}[1]{(\ref{#1})}
 \def\sB{{\mathcal B}}

 \def\sT{{\mathcal T}} \def\sU{{\mathcal U}}

 \def\bE{{\mathbb E}}

 \def\bN{{\mathbb N}} 
\def\bP{{\mathbb P}}  \def\bR{{\mathbb R}}
\def\bS{{\mathbb S}}  
  
 \def\bZ{{\mathbb Z}}

\def\al{\alpha}
\def\lam{\lambda} 
 \def\eps{\varepsilon}
\def\th{\theta} 

 \def\gam{\gamma}
 \def\om{\omega}

\def\to{\rightarrow}

\def\diam{\operatorname{diam}}

\def\half{{\frac{1}2}}

\newtheorem{thmm}{Theorem}[section]
\newproclaim{defn}[thmm]{Definition}
\newproclaim{assu}{Assumption}
\newtheorem{cor}[thmm]{Corollary}

\newtheorem{lem}[thmm]{Lemma}
\newtheorem{propn}[thmm]{Proposition}
\newproclaim{rem}[thmm]{Remark}
\newproclaim{exam}[thmm]{Examples}

\def\SLE{\mathrm{SLE}}
\def\dU{d_\sU}
\def\dS{d_\sU^S}
\def\dE{d_E}

\def\uT{\underline{\mathcal{T}}}
\def\uU{\underline{\mathcal{U}}}
\def\bP{\mathbf{P}}
\makeatother

\begin{document}
\begin{frontmatter}

\title{Subsequential scaling limits of simple random walk on the
two-dimensional uniform spanning tree}
\runtitle{Random walk on the uniform spanning tree}

\begin{aug}
\author[A]{\fnms{M.~T.}~\snm{Barlow}\thanksref{T1}\ead[label=e1]{barlow@math.ubc.ca}},
\author[B]{\fnms{D.~A.}~\snm{Croydon}\thanksref{T2}\ead[label=e2]{d.a.croydon@warwick.ac.uk}}
\and
\author[C]{\fnms{T.}~\snm{Kumagai}\corref{}\thanksref{T3}\ead[label=e3]{kumagai@kurims.kyoto-u.ac.jp}}
\runauthor{M.~T.~Barlow, D.~A. Croydon and T. Kumagai}
\thankstext{T1}{Supported in part by NSERC (Canada).}
\thankstext{T2}{Supported in part by EPSRC First Grant EP/K029657/1.}
\thankstext{T3}{Supported in part by the JSPS Grant-in-Aid for Scientific Research (A) 25247007.}
\affiliation{University of British Columbia, University of Warwick and Kyoto University}
\address[A]{M. T. Barlow\\
Department of Mathematics\\
University of British Columbia\\
Vancouver, B.C., V6T 1Z2\\
Canada\\
\printead{e1}}
\address[B]{D. A. Croydon\\
Department of Statistics\\
University of Warwick\\
Coventry, CV4 7AL\\
United Kingdom\\
\printead{e2}}
\address[C]{T. Kumagai\\
Research Institute for Mathematical Sciences\\
Kyoto University\\
Kyoto 606-8502\\
Japan\\
\printead{e3}}
\end{aug}

%
\received{\smonth{7} \syear{2014}}
%
\revised{\smonth{4} \syear{2015}}

%
\begin{abstract}
The first main result of this paper is that the law of the (rescaled)
two-dimensional uniform spanning tree is tight in a space whose
elements are measured,
rooted real trees continuously embedded into Euclidean space.
Various properties of the intrinsic metrics, measures and embeddings of the
subsequential limits in this space are obtained, with it being proved
in particular
that the Hausdorff dimension of any limit in its intrinsic metric is
almost surely equal to $8/5$.
In addition, the tightness result is applied to deduce that the
annealed law of the simple
random walk on the two-dimensional uniform spanning tree is tight under
a suitable rescaling.
For the limiting processes, which are diffusions on random real trees
embedded into
Euclidean space, detailed transition density estimates are derived.
\end{abstract}

%
\end{frontmatter}

\section{Introduction}\label{sec1}

The study of uniform spanning trees (USTs) has a long history;
in the 1840s Kirchhoff used them in his classic paper \cite{Kirchoff}
on electrical resistance. Much of the recent theory in the probability
literature is based on the discovery that paths in the UST have the
same law as loop erased random walks.
Using this connection, algorithms to construct the UST from random walks
have been given in \cite{Aldous0,Broder,Wilson}. See \cite{BLPS} for
a survey
of the properties of the UST, and a description of Wilson's algorithm,
which will be important for this article, and \cite{Law99} for a survey
of the properties of the loop erased random walk (LERW). We also remark
that USTs can be considered as a boundary case of the random cluster
model; see \cite{Haggstrom}.

In \cite{Schramm}, Schramm studied the scaling limit of the UST in $\bZ
^2$, and this led
him to introduce the SLE process. In \cite{LSW}, it was proved that the
LERW in $\bZ^2$ has $\SLE_2$ as its scaling limit, and this connection
was used in \cite{BM10,Mas09} to improve earlier results of Kenyon
\cite{Kenyon}
on the growth function of two-dimensional LERW. In \cite{BM11},
this good control on the length of LERW paths, combined with Wilson's
algorithm, was used
to obtain volume growth and resistance estimates for the
two-dimensional UST $\sU$.
Using the connection between random walks and electrical resistance,
and the methods
of \cite{BJKS,KM}, these bounds then led to heat kernel bounds for $\sU$.

In this paper, we study scaling limits of $\sU$, as well as the random
walk on it.
While very significant progress in this direction was made on the first
topic in \cite{ABNW,Schramm},
those papers are focused on the topological properties of the
scaling limit as a subset of $\mathbb{R}^2$. Here, we work in a
framework that allows us to
describe properties of the joint scaling limit of the corresponding
intrinsic metric, uniform
measure and simple random walk.

We begin by introducing our main notation.
Throughout this article, $\sU$ will represent the uniform spanning tree
on $\mathbb{Z}^2$, and $\bP$
the probability measure on the probability space on which this is built.
As proved in \cite{Pemantle}, $\sU$ is the local limit of the uniform
spanning tree on $[-n,n]^2\cap\mathbb{Z}^2$ (equipped with
nearest-neighbour bonds) as $n\rightarrow\infty$.
We note that $\sU$ is $\bP$-a.s. indeed a spanning tree of $\mathbb
{Z}^2$,  that is, it is a graph with vertex set $\mathbb{Z}^2$, and
any two of its vertices are connected by a unique path in
$\sU$. We will denote by $\dU$ the intrinsic (shortest path) metric on
the graph~$\sU$, and $\mu_\sU$ the uniform measure on $\sU$ (i.e., the
measure which places a unit mass at each vertex).

To describe the scaling limit of the metric measure space $(\sU, \dU,
\mu_\sU)$,
we work with a Gromov--Hausdorff-type topology of the kind that has
proved useful for studying real trees.
(See \cite{BBI} for an introduction to the classical theory, and \cite
{Evans} for its application to real trees.)
In particular, we will build on the notions of
Gromov--Hausdorff--Prohorov topology of \cite{ADH,Evans,Miermont},
and the topology for spatial trees of~\cite{DL} (cf. the spectral
Gromov--Hausdorff topology of \cite{CHK}).
We extend the metric space $(\sU,\dU)$ to a complete and locally
compact real tree by adding unit line segments along edges. The measure
$\mu_\sU$ is then viewed as a locally finite (atomic) Borel measure on
this space.
To retain information about $\sU$ in the Euclidean topology, we
consider $(\sU,\dU)$ as a spatial tree, that is, as an abstract real
tree embedded into $\mathbb{R}^2$ via a continuous map
$\phi_\sU:\sU\rightarrow\mathbb{R}^2$, which we take in our example to
be just the identity on
vertices, with linear interpolation along edges.
In addition, we will suppose the space $(\sU,\dU)$ is rooted at the
origin of $\mathbb{Z}^2$.
Thus, we define a random quintuplet $(\sU,\dU,\mu_\sU,\phi_\sU, 0)$,
and our
first result (Theorem~\ref{main1} below) is that the law of this object
is tight under rescaling in the
appropriate space of ``measured, rooted spatial trees.'' The principal
advantage of working in this topology is that it allows us to preserve
information about the intrinsic metric $\dU$ and measure $\mu_\sU$;
these parts of the picture were missing
from the earlier
scaling results of \cite{ABNW,Schramm}.

The final ingredient we need in order to state our first main result
comes from the growth function for LERW in $\bZ^2$. This is the function
$G_2(r) = \bE|L_r|$, where $|L_r|$ is the length of a LERW run from 0
until it first exits the ball of radius $r$.
In particular, from the results in \cite{LE,Mas09} we have (see \cite{barLRW},
Corollary~3.15) that there exist constants $c_1,c_2\in
(0,\infty)$ such that
%
\begin{equation}
\label{g2cont} c_1 r^{\kappa} \le G_2(r) \le
c_2 r^{\kappa},
\end{equation}
where the growth exponent $\kappa:=5/4$.
This exponent plays a key role in the comparison of the intrinsic and
Euclidean metrics
on the UST.
We remark that in \cite{BM11}, where the key result of \cite{LE} was
not available, the heat kernel estimates on $\sU$ take on a more
complicated form involving the function $G_2$ and functions derived
from it.

\begin{thmm}\label{main1} If $\bP_\delta$ is the law of the measured, rooted
spatial tree
$(\sU,\delta^{\kappa}\dU, \delta^{2}\mu_\sU,\delta\phi_\sU,0)$ under
$\mathbf{P}$, then the collection
$(\mathbf{P}_\delta)_{\delta\in(0,1)}$ is tight.
\end{thmm}

As already noted, this theorem extends the results of \cite{ABNW,Schramm} to include scaling of the intrinsic metric and uniform
measure. We further note that the tightness in \cite{ABNW,Schramm} was
essentially a finite-dimensional statement, since it described the
shape in Euclidean space of the tree spanning a finite number of
points, while the result above establishes tightness for the entire space.

\begin{rem}
To extend the above theorem to a full convergence result, and establish
that the scaling limit satisfies the
obvious scale invariance properties, it would be sufficient to
characterise the limit uniquely from
a suitable finite-dimensional convergence result.
We expect that such a characterisation will be possible once it is
known that
two-dimensional loop-erased random walk converges as a process. Proving
this is an open problem, but see \cite{AKM,LS,LZ} for recent progress
on proving the convergence of LERW to the $\SLE_2$ curve in its
``natural parameterisation.''
\end{rem}

The tightness in Theorem~\ref{main1} implies the existence of
subsequential scaling limits for the
collection $(\mathbf{P}_\delta)_{\delta\in(0,1)}$ of laws on measured,
rooted spatial trees as \mbox{$\delta\rightarrow0$}.
The following theorem gives a number of properties of these limits.
We note that (a)(ii) translates part of \cite{ABNW}, Theorem~1.2, into
our setting, and
the topological aspects of (c)(i) and (c)(ii) are a restatement of
parts of
\cite{Schramm}, Theorem~1.6. [In particular, the set $\phi_\mathcal
{T}(\mathcal{T}^o)$ that appears
in the statement of our result is identical to Schramm's notion of the
``trunk'' for the UST scaling limit; see Lemma~\ref{trunklem}.]
We do not expect the powers of logarithms and log-logarithms in
\eqref{globalmeas} and \eqref{e:rootmeas} to be optimal.
We write $\operatorname{ deg}_\mathcal{T}(x)$ for the degree of a point $x$ in a
real tree $\mathcal{T}$, that is, the number of connected components of
$\mathcal{T}\setminus\{x\}$,
$|A|$ to represent the cardinality of a subset $A\subseteq\mathcal
{T}$, and $\mathcal{L}$ to represent Lebesgue measure on $\mathbb{R}^2$.

\begin{thmm}\label{main12}
If $\tilde{\mathbf{P}}$ is a subsequential limit of $(\mathbf{P}_\delta
)_{\delta\in(0,1)}$,
then for $\tilde{\mathbf{P}}$-a.e. measured, rooted spatial tree
$(\mathcal{T},d_\mathcal{T},\mu_{\mathcal{T}},\phi_\mathcal{T},\rho
_\mathcal{T})$ it holds that:
\begin{enumerate}[(a)]
\item[(a)]
\begin{enumerate}[(iii)]
\item[(i)] the Hausdorff dimension of the complete and locally compact real
tree $(\mathcal{T},d_\mathcal{T})$ is given by
\begin{equation}
\label{e:dfdef} d_f:=\frac{2}{\kappa}=\frac{8}{5};
\end{equation}
\item[(ii)]$(\mathcal{T},d_\mathcal{T})$ has precisely one end at infinity
[i.e., there exists a unique isometric embedding of $\mathbb{R}_+$ into
$(\mathcal{T},d_\mathcal{T})$ that maps $0$ to $\rho_\mathcal{T}$];
\end{enumerate}
\item[(b)]
\begin{enumerate}[(iii)]
\item[(i)] the locally finite Borel measure $\mu_\mathcal{T}$ on
$(\mathcal{T},d_\mathcal{T})$ is nonatomic and supported on the
leaves of
$\mathcal{T}$, that is, $\mu_\mathcal{T}(\mathcal{T}^o)=0$, where
$\mathcal{T}^o:=\mathcal{T}\setminus\{x\in\mathcal{T}: \operatorname{
deg}_\mathcal{T}(x)=1\}$;
\item[(ii)] given $R>0$, there exists a random $r_0(\mathcal{T})>0$ and
deterministic $c_1,c_2\in(0,\infty)$ such that
%
\begin{equation}
\label{globalmeas} c_1r^{d_f}\bigl(\log r^{-1}
\bigr)^{-80} \leq\mu_\mathcal{T} \bigl(B_\mathcal {T}(x,r)
\bigr)\leq c_2 r^{d_f}\bigl(\log r^{-1}
\bigr)^{80},
\end{equation}
for every $x\in B_\mathcal{T}(\rho_\mathcal{T},R)$ and $r\in
(0,r_0(\mathcal{T}))$, where $B_\mathcal{T}(x,r)$ is the open ball
centred at $x$ with radius $r$ in $(\mathcal{T},d_\mathcal{T})$;
\item[(iii)] there exists a random $r_0(\mathcal{T})>0$ and deterministic
$c_1,c_2\in(0,\infty)$ such that
\begin{equation}
\label{e:rootmeas} c_1r^{d_f}\bigl(\log\log r^{-1}
\bigr)^{-9} \leq\mu_\mathcal{T} \bigl(B_\mathcal {T}(
\rho_\mathcal{T},r) \bigr)\leq c_2r^{d_f}\bigl(\log
\log r^{-1}\bigr)^{3},
\end{equation}
for every $r\in(0,r_0(\mathcal{T}))$;
\end{enumerate}
\item[(c)]
\begin{enumerate}[(iii)]
\item[(i)] the restriction of the continuous map $\phi_\mathcal{T}:\mathcal
{T}\rightarrow\mathbb{R}^2$ to $\mathcal{T}^o$ is a homeomorphism
between $\mathcal{T}^o$ (equipped with the topology induced by the
metric $d_\mathcal{T}$) and its image $\phi_\mathcal{T}(\mathcal{T}^o)$
(equipped with the Euclidean topology), the latter of which is dense in
$\mathbb{R}^2$;
\item[(ii)]$\max_{x\in\mathcal{T}}\operatorname{ deg}_\mathcal{T}(x)=3=\max_{x\in
\mathbb{R}^2}|\phi_\mathcal{T}^{-1}(x)|$;
\item[(iii)]$\mu_\mathcal{T}=\mathcal{L}\circ\phi_\mathcal{T}$.
\end{enumerate}
\end{enumerate}
\end{thmm}

\begin{figure}

\includegraphics{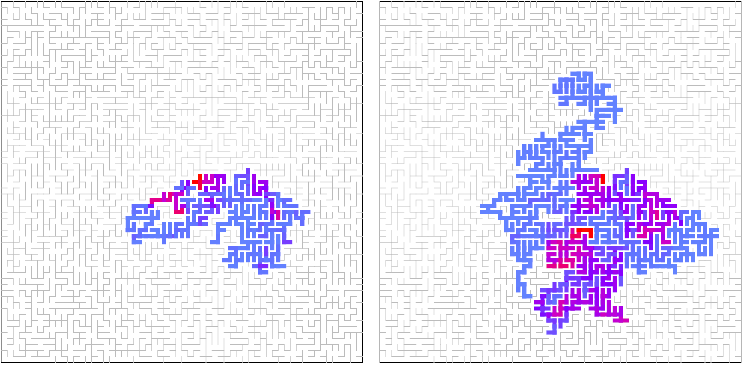}

\caption{The range of a realisation of the simple random walk on
uniform spanning tree on a $60\times60$ box (with wired boundary
conditions), shown after 5000 and 50,000 steps. From most to least
crossed edges, colours blend from red to blue.}\label{srwfig}
\end{figure}

The second topic of this paper is the scaling limit of the simple
random walk (SRW) on the two-dimensional UST. For a given realisation
of the graph $\sU$, the SRW on $\sU$ is the discrete time Markov
process $X^\sU=((X_n)_{n\geq0},({P}_x^\sU)_{x\in\mathbb{Z}^2})$ which
at each time step jumps from its current location to a uniformly chosen
neighbour in $\sU$ (considered as a graph); see Figure~\ref{srwfig}.
For $x\in\mathbb{Z}^2$, the law ${P}_x^\sU$ is called the \emph{quenched} law of the simple random walk on $\sU$ started at $x$. Since
$0$ is always an element of $\sU$, we can define the \emph{annealed} or
\emph{averaged} law $\mathbb{P}$ as the semi-direct product of the
environment law $\mathbf{P}$ and the quenched law ${P}_0^\sU$ by setting
%
\begin{equation}
\label{annlaw} \mathbb{P} (\cdot ):=\int P_0^\sU(\cdot)
\,d\mathbf{P}.
\end{equation}
It is this measure for which we will deduce scaling behaviour.

Techniques for deriving the scaling limits of random walks on the kinds
of trees generated
by critical branching processes have previously been developed in \cite
{Croy1,Croy2,Croy3};
see also \cite{K14}, Chapter~7, for a survey.
In the present work, we adapt these to prove a general result of the
following form; see Theorem~\ref{srwthm} below for details and some
additional technical conditions.
If we have a sequence of graph trees $(T_n)$, $n\geq1$, each equipped
with its intrinsic metric
$d_{T_n}$, a measure $\mu_{{T}_n}$, an embedding $\phi
_{{T}_n}:T_n\rightarrow\mathbb{R}^2$ and a distinguished root vertex
$\rho_n$, for which there exist null sequences
$(a_n)_{n\geq1}$, $(b_n)_{n\geq1}$, $(c_n)_{n\geq1}$ with
$b_n=o(a_n)$ such that
$({T}_n,a_nd_{{T}_n},b_n\mu_{{T}_n},c_n\phi_{{T}_n},\rho
_{{T}_n})\rightarrow(\mathcal{T},d_\mathcal{T},\mu_{\mathcal{T}},\phi
_\mathcal{T},\rho_\mathcal{T})$
in the space of measured, rooted spatial trees, then the corresponding
rescaled random walks
$(c_n\phi_{T_n}(X^{T_n}_{t/a_nb_n}))_{t\geq0}$ converge in distribution.
Further, the limiting process can be written as $(\phi_{\mathcal
{T}}(X^{\mathcal{T}}_{t}))_{t\geq0}$, where $X^\mathcal{T}=((X^\mathcal
{T}_t)_{t\geq0},
(P^\mathcal{T}_x)_{x\in\mathcal{T}})$ is the canonical Brownian motion
on $(\mathcal{T},d_\mathcal{T},\mu_\mathcal{T})$, as constructed in
\cite{AEW}, for example, cf. \cite{Kigami}.
(We give a brief introduction to Brownian motion on measured real trees
at the start of Section~\ref{srwsec}.)

Combining Theorems \ref{main1} and \ref{srwthm} we obtain
the following theorem, which establishes the existence of subsequential
scaling limits
for the annealed law of the simple random walk on $\sU$. Given the volume
estimates~\eqref{globalmeas}, the general results of \cite{Croy-1} yield
sub-diffusive transition density bounds for the limiting diffusion.
These demonstrate that, uniformly over bounded regions of space, the
transition density in
question has at most
logarithmic fluctuations from the leading order polynomial terms in
both the on-diagonal and exponential off-diagonal decay parts.
In Section~\ref{hklimsec}, we also deduce pointwise on-diagonal
estimates with
only log-logarithmic fluctuations (cf. the discrete result of \cite{BM11},
Theorem~4.5(a)),
as well as annealed on-diagonal polynomial bounds. We note the
similarity between these results
and the transition density estimates for the Brownian continuum random tree
given in~\cite{Croy0}.

\begin{thmm}\label{main2}
If $(\mathbf{P}_{\delta_i})_{i\geq1}$ is a convergent sequence with
limit $\tilde{\mathbf{P}}$,
then the following statements hold:

\textup{(a)} The annealed law of $(\phi_\mathcal{T}(X^\mathcal{T}_t))_{t\geq0}$,
where $X^\mathcal{T}$ is Brownian motion on
$(\mathcal{T},d_\mathcal{T},\mu_\mathcal{T})$ started from $\rho
_\mathcal{T}$, that is,
%
\begin{equation}
\label{annp} \tilde{\mathbb{P}} (\cdot ) :=\int P_{\rho_\mathcal{T}}^\mathcal{T}
\circ\phi_\mathcal{T}^{-1}(\cdot ) \,d\tilde{\mathbf{P}},
\end{equation}
is a well-defined probability measure on $C(\mathbb{R}_+,\mathbb
{R}^2)$.

\textup{(b)} If $\mathbb{P}_\delta$ is defined to be the law of
$(\delta X^\sU_{\delta^{-\kappa d_w}t})_{t\geq0}$ under $\mathbb{P}$,
where the walk dimension $d_w$ of $\sU$ is defined by
\[
d_w:=1+d_f=\tfrac{13}{5},
\]
then $(\mathbb{P}_{\delta_i})_{i\geq1}$ converges to $\tilde{\mathbb
{P}}$.

\textup{(c)} $\tilde{\mathbf{P}}$-a.s., the process $X^\mathcal{T}$ is recurrent
and admits a jointly continuous transition density $(p_t^\mathcal
{T}(x,y))_{x,y\in\mathcal{T},t>0}$. Moreover, it $\tilde{\mathbf
{P}}$-a.s. holds that, for any $R>0$, there exist random constants
$c_i(\mathcal{T})$ and $t_0(\mathcal{T}) \in(0,\infty)$
and deterministic constants $\theta_1,\theta_2,\theta_3,\theta_4\in
(0,\infty)$
(not depending on $R$) such that
\begin{eqnarray*}
p_t^\mathcal{T}(x,y)&\leq& c_1(\mathcal{T})
t^{-d_f/d_w}\ell \bigl(t^{-1}\bigr)^{\theta_1}\\
&&{}\times \exp \biggl
\{-c_2(\mathcal{T}) \biggl(\frac{d_\mathcal
{T}(x,y)^{d_w}}{t} \biggr)^{1/(d_w-1)}
\ell\bigl(d_\mathcal{T}(x,y)/t\bigr)^{-\theta_2} \biggr\},
\\
p_t^\mathcal{T}(x,y)&\geq& c_3(\mathcal{T})
t^{-d_f/d_w}\ell \bigl(t^{-1}\bigr)^{-\theta_3} \\
&&{}\times\exp \biggl
\{-c_4(\mathcal{T}) \biggl(\frac{d_\mathcal
{T}(x,y)^{d_w}}{t} \biggr)^{1/(d_w-1)}
\ell\bigl(d_\mathcal{T}(x,y)/t\bigr)^{\theta
_4} \biggr\},
\end{eqnarray*}
for all $x,y\in B_\mathcal{T}(\rho_\mathcal{T},R)$, $t\in(0,t_0(\mathcal
{T}))$, where $\ell(x):=1\vee\log x$.
\end{thmm}

\begin{rem}
If follows that for $\tilde{\mathbf{P}}$-a.e. realisation of
$(\mathcal{T},d_\mathcal{T},\mu_{\mathcal{T}},\phi_\mathcal{T},\rho
_\mathcal{T})$,
we have that
\[
-\lim_{t\rightarrow0}\frac{2\log p^\mathcal{T}_t(x,x)}{\log t}=\frac
{2d_f}{1+d_f}=
\frac{16}{13} \qquad\mbox{for every $x\in\mathcal{T}$. }
\]
Using the language of diffusions on fractals, this means that the
spectral dimension
of the limiting tree is $\tilde{\mathbf{P}}$-a.s. equal to $16/13$,
which is the same as for the discrete model (see \cite{BM11}).
\end{rem}

The remainder of this article is organised as follows.
In Section~\ref{estsec}, we prove some key estimates for $\sU$,
which enable us to compare distances in the Euclidean and intrinsic
metrics on this set. These allow us to extend some of the volume
estimates of \cite{BM11}. In Section~\ref{topsec}, we introduce our
topology for measured, rooted spatial trees,
and in Section~\ref{tightsec} we prove
tightness in this topology for the rescaled trees.
The properties of limiting trees are studied in Section~\ref{propsec}.
Following this, we turn our attention to the simple random walk on $\sU$,
establishing in Section~\ref{srwsec} a general convergence result for
simple random walks on measured, rooted spatial trees and applying this
to the two-dimensional UST.
In addition, we explain how this convergence result can be applied to
branching random walks and trees without embeddings.
In Section~\ref{hklimsec}, we then derive the transition density
estimates for the limiting diffusion.

We write $c$ or $c_i$ for constants in $(0,\infty)$; these will be
universal and nonrandom, but may change in
value from line to line. We use the notation $c_i(\mathcal{T})$ for
(random) constants which depend on the tree $\mathcal{T}$.

\section{UST estimates}\label{estsec}

In this section, we obtain estimates
for the two-\break dimensional UST $\sU$, which improve those in \cite{BM11}.
Our arguments will depend heavily on Wilson's algorithm, which gives the
construction of $\mathcal{U}$ in terms of LERW.
In particular, we can construct $\mathcal{U}$ by first running an
infinite loop-erased random walk
from $0$ to $\infty$ (for details of this see \cite{Mas09}), and then,
sequentially running through vertices $x\in\mathbb{Z}^2\setminus\{0\}$,
adding a loop-erased random walk path from $x\in\mathbb{Z}^2$ to the
part of the tree already created.
We remark that $\sU$ is a one-ended tree; see \cite{BLPS}.

We will consider three metrics on $\sU$, which we now introduce.
We define $\dE$ to be the Euclidean metric on $\bZ^2$,
and write $B_E(x,r)=\{y: \dE(x,y) \le r \}$ for balls in this metric.
For $x \in\bZ^2$, we let $\gam(x,y)$ be the unique path in $\sU$
between $x$ and $y$.
We define the intrinsic (shortest path) metric $\dU$ by setting $\dU
(x,y):=|\gamma(x,y)|$, that is,
the number of edges on the path $\gam(x,y)$, and write $B_\sU(x,r)$ for balls
in this metric.
Finally, it will also be helpful to use a modification of a metric
introduced by Schramm in \cite{Schramm},
given by
%
\begin{equation}
\label{sm} \dS(x,y):= \diam\bigl( \gam(x,y)\bigr),
\end{equation}
where the right-hand side refers to the diameter of $\gam(x,y)$ in the
metric $\dE$.

We begin by recalling the comparison between $B_\sU(0,r^{1/\kappa})$
and $B_E(0,r)$ and the estimates on the size of $|B_\sU(0,r)|$ from
\cite{BM11}.

\begin{thmm}[(See \cite{BM11}, Theorems 1.1, 1.2)]\label{vollemust}
\textup{(a)} There exist $c_1,c_2$ such
that for every $r\geq1$ and $\lambda\geq1$,
\begin{eqnarray*}
\bP \bigl(B_\sU\bigl(0, \lam^{-1} r^\kappa\bigr)
\not\subset B_E(0,r) \bigr)& \leq& c_1e^{-c_2\lambda^{2/3} },
\\
\bP \bigl(B_E(0, r) \not\subset B_\sU\bigl(0, \lam
r^{\kappa}\bigr) \bigr)& \leq &c_1\lambda^{-1/5}.
\end{eqnarray*}
\textup{(b)} There exist $c_1,c_2$ such that for every $r\geq1$ and $\lambda
\geq1$,
\begin{eqnarray*}
\bP \bigl(\bigl|B_\sU(0,r)\bigr| \geq\lam r^{d_f} \bigr)& \leq&
c_1e^{-c_2\lambda
^{1/3}},
\\
\bP \bigl(\bigl|B_\sU(0,r)\bigr| \leq\lam^{-1}r^{d_f}
\bigr)& \leq &c_1e^{-c_2\lambda^{1/9}}.
\end{eqnarray*}
\end{thmm}

Since the law of $\sU$ is translation invariant, the above result also
holds for $B_\sU(x,r)$ for any
$x \in\bZ^d$. However, we wish to have these bounds (for suitable
$r,n$) for every
$x \in B_E(0,n)$; obtaining such uniform estimates is one of the main
goals of this section. If we use a simple union bound, as, for example,
in \cite{BM11}, (4.47), we obtain an error estimate of the form $n^2
\exp( -\lam^{c})$, which is only small when $\lam\gg(\log n)^{1/c}$.
To improve this, for a suitable $\delta=\delta(\lam)>0$ we choose a
$\delta$-cover $D$ of $B_E(0,n)$ with $|D| \le c \delta^{-2}$. (Recall
that a subset $A$ of $\bZ^2$ is called a \emph{$\lam$-cover} if every
point of $\bZ^2$ is within distance $\lam$ of a point of $A$.) We then
obtain good behaviour of $B_\sU(x,r)$ for all $x \in D$, except on a
set of probability $|D| \exp( -\lam^{c})$. Using a ``filling in lemma''
(see Lemma~\ref{L:fillin} below), together with some additional bounds,
we are able to extend this good behaviour
to $B_\sU(y,r)$ for all $y \in B_E(0,n)$; see Proposition~\ref
{P:keypropals} below for a uniform version of part (b) in particular.
An example of the kind of additional result that we need is that if
$d_E(x,y)=3r$ then every path in $\sU$ between $B_E(x,r)$ and
$B_E(y,r)$ is of length at least $c r^{\kappa}$, except on a set of
trees of small
probability; note that \cite{BM10}, Theorem~1.2, shows that with high
probability the unique infinite self avoiding path in $\sU$ started
from $x$ takes at least $c r^{\kappa}$ steps to escape a Euclidean ball
of radius $r$, but again this does not readily extend to a uniform
bound and so further work is required.

We proceed by introducing some further notation and results from \cite
{BM11}. Let $\gam_x = \gam(x, \infty)$ be the unique infinite self
avoiding path in $\sU$ started at $x$;
by Wilson's algorithm $\gam_x$ has the law of the loop-erased random
walk from $x$ to $\infty$.
Write $\gam_x[i]$ for the $i$th point on $\gam_x$, and let
$\tau_{y,r}=\tau_{y,r}(\gam_x) = \min\{ i: \gam_x[i] \notin B_E(y, r)
\}$.
Whenever we use notation such as $\gam_x[\tau_{y,r}]$, the exit time
$\tau_{y,r}$ will always
be for the path $\gam_x$. We define the segment of the path $\gam_x$
between its $i$th and $j$th
points by $\gam_x[i,j]= (\gam_x[i], \gam_x[i+1], \ldots, \gam_x[j])$, and
define $\gam_x[i, \infty)$ in a similar fashion. For such paths, the
following was established in \cite{BM11}.

\begin{lem}[(See \cite{BM11}, Lemma~2.4)]\label{L:pathsize}
There exists $c_1$ such that for every
$r\geq1$ and $k\geq2$,
\[
\bP\bigl( \gam_x[\tau_{x,kr},\infty) \cap
B_E(x,r)\neq\varnothing\bigr) \le c_1k^{-1}.
\]
\end{lem}

We next give the filling in lemma that we will use several times. This
is a small extension of
\cite{BM11}, Proposition~3.2. Note that \cite{BM10}, Proposition~6.2,
shows that the function $G(r)$ considered in \cite{BM11} is comparable
with the function $G_2(r)$ appearing in (\ref{g2cont}).

%
\begin{lem}
\label{L:fillin} There exist constants $c_1,c_2\in(0,\infty)$ such
that for each
$\delta\le1$ the following holds. Let $r \geq1$, and $U_0$ be a
fixed tree in $\bZ^2$
with the property that $\dE(x, U_0) \le\delta r$ for each $x \in B_E(0,r)$.
Let $\sU$ be the random spanning tree in $\bZ^2$ obtained by running
Wilson's algorithm with root $U_0$ (i.e., starting from the tree $U_0$).
Then there exists an event $G$ such that $\bP(G^c) \le c_1 e^{-c_2
\delta^{-1/3}}$, and on $G$ we have that for all $x \in B_E(0,r/2)$,
\[
\dU(x, U_0) \le \bigl( \delta^{1/2} r
\bigr)^\kappa; \qquad\dS(x, U_0) \le \delta^{1/2} r;\qquad
\gam(x, U_0) \subset B_E(0,r).
\]
\end{lem}

\begin{pf}
Except for the bound involving $\dS$ this is proved in \cite{BM11}.
[Note that the hypothesis there that $U_0$ connects $0$ to
$B_E(0,2r)^c$ is
unnecessary.] The proof of the $\dS$ bound is similar.
\end{pf}

The following sequence of lemmas will improve the results in \cite{BM11}
on the comparison of the metrics $\dU$, $d_E$ and $d^S_\sU$.
Fix for now $r, k\geq1$ and $x \in\bZ^2$,
and choose points $z_j$ on $\gam_x$ so that $z_0=x$ and $z_j = \gam
_x[s_j]$, where $ s_j = \min\{ i: \dE( \gam_x[i], \{ z_0, \ldots, z_{j-1}
\}) \ge r/k \}$.
Let $N=N(r,k)=\max\{j: s_j \le\tau_{x,r}(\gam_x) \}$. Moreover, define
a collection of disjoint balls $\sB_{r,k}= \{ B_j=B_E(z_j,r/3k), j=1,
\ldots, N(r,k) \}$.
These depend on the path $\gam_x$, and when we need to
recall this we will write $\sB_{r,k}(\gam_x)$. Let $ a= 1+ k^{-1/8}$,
and set
\[
F_1(x,r,k) =\bigl\{\mbox{$\gam_x[
\tau_{x,ar}, \infty]$ hits fewer than $k^{1/2} $ of
$B_1, \ldots,B_{N(r,k)}$}\bigr\}.
\]

\begin{lem}\label{L:F1} There exist constants $c_1,c_2\in(0,\infty)$ such
that, if $r,k\geq1$
and $x\in\mathbb{Z}^2$, then
\[
\bP \bigl( F_1(x,r,k)^c \bigr) \le c_1 e
^{-c_2 k^{1/8} }.
\]
\end{lem}

\begin{pf*}{Proof (\textup{see} \cite{BM11}, \textup{Lemma}~3.7)}
Write $\tau_{s} = \tau_{x,s}(\gamma_x)$, and let $b= e^{k^{1/8}}\geq2$.
Then by Lemma~\ref{L:pathsize},
%
\begin{equation}
\label{e:expret} \bP\bigl( \gam_x[\tau_{br}, \infty] \cap
B_E(x,r) \neq\varnothing\bigr) \le c b^{-1} =
ce^{-k^{1/8}}.
\end{equation}
If $\gam_x[\tau_{ar}, \infty]$ hits more than $k^{1/2}$ balls from the
family $\sB_{r,k}(\gam_x)$, then either $\gam_x$ hits $B_E(0,r)$ after
time $\tau_{br}$,
or $\gam_x[\tau_{ar}, \tau_{br}]$ hits more than $k^{1/2}$ balls.
Given \eqref{e:expret}, it is therefore sufficient to prove that
%
\begin{equation}
\label{e:43b} \bP\bigl( \gam_x[\tau_{ar},
\tau_{br}] \mbox{ hits more than $k^{1/2}$ balls} \bigr) \le
c_1 e ^{-c_2 k^{1/8} }.
\end{equation}
Let $S$ be a simple random walk on $\mathbb{Z}^2$ started at $x$, $L'$
be the loop-erasure of
$S[0, \tau_{x,4br}(S)]$, and $L'' = L'[ \tau_{x,ar}(L'), \tau_{x,br}(L')]$.
Then by \cite{Mas09}, Corollary~4.5, in order to prove \eqref{e:43b},
it is sufficient to prove that
\[
\bP\bigl( L'' \mbox{ hits more than
$k^{1/2}$ balls in $\sB_{r,k}\bigl(L'\bigr)$}
\bigr) \le c_1 e ^{-c_2 k^{1/8} }.
\]
Define stopping times for $S$ by letting $T_0=\tau_{x,ar}(S)$ and for
$j \ge1$,
setting $R_j = \min\{ n \ge T_{j-1}: S_n \in B_E(x, r) \}$ and $T_j =
\min\{ n \ge R_{j}: S_n \notin B_E(x, ar) \}$.
Note that the balls in $\sB_{r,k}(L')$ can only be hit by $S$ in the intervals
$[R_j, T_j]$ for $j \ge1$. Let $M =\min\{j: R_j \ge\tau_{x,4 br}(S) \}$.
Then, by the result of \cite{Lawlerbook}, Exercise~1.6.8,
$ \bP( M = j+1 | M>j) \ge { c(\log(a r) -\log r)}/({\log( 4br) -
\log r})\ge c k^{-2/8}$.
Hence, $\bP( M \ge k^{3/8}) \le c_1 \exp( - c_2 k^{1/8} )$.
Now, for each $j\ge1$, let $L_j$ be the loop-erasure of $S[0, T_j]$,
$\al_j$ be the first exit by $L_j$ from $B_E(x,ar)$, and $\beta_j$ be
the number of steps in $L_j$.
If $L''$ hits more than $k^{1/2}$ balls in $\sB_{r,k}(L')$, then there
must exist some $j \le M$ such that $L_j[\al_j, \beta_j]$ hits more
than $k^{1/2}$
of the balls in the collection $\sB_{r,k}(L_j)$. Hence, if $M \le
k^{3/8}$ and $L''$ hits more than $k^{1/2}$ balls in $\sB_{r,k}(L')$,
then $S$ must hit more than $k^{1/8}$ balls in $\sB_{r,k}(L_j)$ in one
of the intervals $[R_j, T_j]$, without hitting the path $L_j[0, \al
_j]$. However, by Beurling's estimate (see \cite{Lawlerbook}, Lemma~2.5.3, e.g.), the probability of this event is less
than $c_1 \exp(-c_2 k^{1/8})$.
Combining these estimates completes the proof.
\end{pf*}

\begin{figure}

\includegraphics{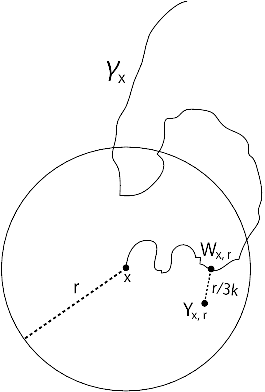}

\caption{A sample of $A_1(x,r,k)$ in Lemma \protect\ref{L:basicLP}.}\label{srkk}
\end{figure}

Our next lemma shows that if $D_0$ is $\delta r$-cover of
$B_{E}(x,2r)$, then with high probability we can find points $Y_{x,r}$
and $W_{x,r}$
which are close to the boundary of $B_E(x,r)$ and to each other, and
such that $Y_{x,r} \in D_0$
and $W_{x,r}\in\gamma_x\cap B_E(x,r)$. (See Figure~\ref{srkk}.) In the
proof, we refer to the event $F_2(x,r,k) = \{ 8k^{-1/4} r^\kappa\le
\tau_{x,r}(\gam_x) \le k^{1/4}r^\kappa \}$. From \cite{BM10}, Theorems
5.8, 6.1, we have
\begin{equation}
\label{e:F2bnd} \bP\bigl( F_2(x,r,k)^c \bigr) \le
c_1 \exp\bigl( - c_2 k^{1/6}\bigr).
\end{equation}

\begin{lem}\label{L:basicLP}
Let $r \ge1$, $k\geq2$, $x \in\bZ^2$, and $D_0 \subset\bZ^2$ satisfy
$B_E(x,2r) \subset\bigcup_{y \in D_0} B_E(y,r/18k)$. Then there exists an
event $A_1=A_1(x,r,k)$, defined in \eqref{e:A1def} below, which satisfies
\begin{equation}
\label{e:Acub} \bP\bigl( A_1^c \bigr) \le
e^{-k^{1/8}},
\end{equation}
and on $A_1(x,r,k)$ there exists $T \le\tau_{x,r}(\gam_x)$ such that,
writing $W_{x,r}=\gam_x(T)$:
\begin{longlist}[(a)]
\item[(a)] $k^{-1/4} r^\kappa\le T \le k^{1/4} r^\kappa$;

\item[(b)] $ a^{-2}r \le d_E(x,W_{x,r}) \le r$;

\item[(c)] there exists $Y_{x,r} \in D_0$ such that $d_E( Y_{x,r}, W_{x,r})
\le{r}/{3k}$, $\dS(Y_{x,r}, W_{x,r} ) \le{2r}/{3k}$ and also $\dU
(Y_{x,r},W_{x,r})\le c_1 (r/k)^\kappa$.
\end{longlist}
\end{lem}

\begin{pf} Fix $k\geq1$ and recall that $a=1+ k^{-1/8}$. Suppose
that the event
\begin{equation}
\label{e:F122} F_1(x,r/a,k) \cap F_2
\bigl(x,r/a^2,k\bigr) \cap F_2(x,r,k)
\end{equation}
occurs. Write $\tau_s = \tau_{x,s}(\gam_x)$, $T_1= \tau_{r/a^2}$, and
$T_2= \tau_{r/a}$. Let $J_0=J_0(\om)$ be the set of $j$
such that $z_j \in\gam_x[T_1,T_2]$ and $B_E(z_j, r/3ak ) \subset
B_E(x,r/a)\setminus B_E(x, r/a^2)$.
Then $|J_0| \ge c k^{7/8}$. Since $F_1(x,r/a,k)$ holds,
at most $k^{1/2}$ of the balls\break $(B_E(z_j, r/3ak),  j \in J_0)$ are hit
by $\gam_x[ \tau_{r}, \infty]$. So if $J=J(\om)$ is the set of $j \in
J_0$ such that
$ B_E(z_j, r/3ak) \cap\gam_x[ \tau_{r}, \infty] = \varnothing$, then
$|J| \ge k^{7/8}-k^{1/2}\geq ck^{3/4}$. For each $j\in J$, we can find
a point $y_j \in D_0$ with $\dE(y_j,z_j) \le r/18k$. Hence,
$B_E(y_j, r/18k) \cap \gam_x[T_1, T_2] \neq\varnothing$, while
$B_E(y_j, r/9 k) \cap \gam_x[ \tau_{r}, \infty] = \varnothing$. Note
that $B_E(y_j,r/9k)$ may
however intersect the path $\gam_x$ in the interval $[T_2,\tau_r]$.

For the remainder of the proof, it will be helpful to regard
$\gam_x$ as a fixed deterministic path which satisfies the conditions
in \eqref{e:F122}.
For each $j\in J$, let $X^j$ be a SRW started at $y_j$ and run until it
hits $\gam_x$, and
let $L^j$ be the loop-erasure of $X^j$.
Let
\[
H_j= \bigl\{ X^j \mbox{ hits $\gam_x$
before it exits } B_E(z_j, r/3ak), \bigl|L^j\bigr| \le
c_0 (r/3k)^\kappa \bigr\}.
\]
By \cite{BM11}, Theorem~2.2, we have [taking $D=\bZ^2\setminus\gam_x$
and $D'=D\cap B_E(z_j, r/3ak)$],
\[
\bP\bigl( \bigl|L^j \cap B_E(z_j, r/3ak) \bigr| >
\lam(r/k)^\kappa\bigr) \le c_1 \exp( -c_2 \lam).
\]
So, by Beurling's estimate (see \cite{Lawlerbook}, Lemma~2.5.3, e.g.), we can choose $c_0$ so that there exists $p>0$ such that $\bP
(H_j) \ge p$.

Recall now the implementation of Wilson's
algorithm using ``stacks'' (see \cite{Wilson}). For each $j$, assume we
have stack variables $\xi_{x,i}$ for $x \in B_E(z_j,r/3ak)$. We use
these to make a random walk path $X^j$ started at $y_j$ and run either
it hits $\gam_x$ or leaves $B_E(z_j,r/3ak)$. Thus, the event $H_j$ is
measurable with respect to $\sigma( \xi_{x,i}, i\ge1, x \in
B_E(z_j,r/3ak))$. We now consider the $y_j$ one at a time, and continue
until either we
obtain a success, or we have tried $k^{3/4}$ of the points $y_j$. Since
these events are
independent,
if $H$ is the event that we obtain a success, then
\begin{equation}
\label{eq:hcdef} \bP\bigl( H^c\bigr) \le (1-p)^{k^{3/4}} \le
c_1 \exp\bigl(-c_2 k^{3/4}\bigr).
\end{equation}
If $H$ occurs, with a success for $y_j$, set $Y=Y_{x,r}=y_j$, let
$W=W_{x,r}$ be the
point where $X^j$ hits $\gam_x$, and let $T$ be such that $\gam
_x(T)=W$. We take
\begin{equation}
\label{e:A1def} A_1=A_1(x,r,k) = H \cap
F_1(x,r/a,k) \cap F_2\bigl(x,r/a^2,k\bigr)
\cap F_2(x,r,k).
\end{equation}
By Lemma~\ref{L:F1}, \eqref{e:F2bnd} and \eqref{eq:hcdef},
we have the upper bound \eqref{e:Acub} on $\bP(A_1^c)$.

Finally, suppose that $A_1(x,r,k)$ occurs. By construction, we have $\dE
(Y,\break W) \le r/3k$,
and since the path $X^j$ lies inside $B_E(z_j,r/3k)$ we also have $\dS
(Y,W) \le2r/3k$.
The definition of the event $H_j$ gives that $\dU(Y,W) \le c(r/k)^\kappa$.
Since $X^j$ hits $\gam_x$ inside $B_E(z_j, r/3ak)$, we must have $W\in
B_E(x,r)\setminus B_E(x,a^{-2}r)$. Moreover, because $j \in J$, $T\le
\tau_r(\gam_x)$, so since $F_2(x,k,r)$ holds we
have $T \le k^{1/4} r^\kappa$. Since $B_E(z_j, r/3ak) \cap B_E(x,
r/a^2) =\varnothing$,
we must also have $T \ge\tau_{r/a^2}$, and so $T \ge 8k^{-1/4}
(r/a^2)^\kappa\ge k^{-1/4} r^\kappa$.
\end{pf}

The next lemma allows us to compare $d_\sU^S$ and $d_\sU$ on a large
family of paths in a ball.

\begin{lem}\label{L:LP1}
Let $r \ge1$, $k \ge8$, and $x \in\bZ^2$. Set $M_1=e^{k^{1/8}/27}$,
$M_2=e^{k^{1/8}/3}$, $R_i=rM_i$, and let $D_0 \subset B_E(x,2R_2)$
satisfy $|D_0| \le c k^2 M_2^2$, $|D_0 \cap B_E(x, 2R_1)| \le c k^2
M_1^2$, $B_E(x,2 R_2) \subset\bigcup_{y \in D_0} B_E(y, r/18k)$. Write
$D_1= D_0 \cap B_E(x,2 R_1)$.
Then there exist constants $b_1, b_2$ and an event $A_2=A_2(x,r,k)$ with
\begin{equation}
\label{e:pacub1} \bP\bigl(A_2^c\bigr) \le c \exp\bigl( -
k^{1/8}/4\bigr),
\end{equation}
such that on $A_2$ the following holds for every $y \in D_1$:
\begin{longlist}[(a)]
\item[(a)] $\gam_y[ \tau_{x, R_2},\infty] \cap B_E(x,4R_1) =\varnothing$;

\item[(b)] if $x_1, x_2 \in\gam_y[0, \tau_{x,R_2}]$ and $\dS(x_1, x_2) > b_2
r$, then
$\dU(x_1, x_2) \ge\half k^{-1/4} r^\kappa$;

\item[(c)] if $x_1, x_2 \in\gam_y[0, \tau_{x,R_2}]$ and $\dS(x_1, x_2) < b_1
r$, then
$\dU(x_1, x_2) \le2 k^{1/4} r^\kappa$.
\end{longlist}
\end{lem}

\begin{pf} For $y \in D_1$, let $ F_3(y,r,k)= \{ \gam_y[ \tau_{x,
R_2},\infty] \cap B_E(x,4R_1) =\varnothing\}$.
By Lemma~\ref{L:pathsize}, we have $\bP( F_3^c) \le{c M_1}/{M_2}$. Now set
\[
A_2= \biggl( \bigcap_{y \in D_0}
A_1(y,r,k) \biggr) \cap \biggl( \bigcap
_{y \in D_1} F_3(y,r,k) \biggr),
\]
where $A_1(y,r,k)$ is the event defined by \eqref{e:A1def}. From \eqref
{e:Acub}, we note that
\[
\bP\bigl( A_2^c\bigr) \le c k^2
M_2^2 e^{-k^{1/8}} + c M_1^2
k^2 M_1 M_2^{-1} \le c \exp\bigl(
- k^{1/8}/4\bigr).
\]

Now suppose that $A_2$ holds, and let $y \in D_1$. It is immediate
that (a) holds.
Write $W_0=Y_0=y$, and let $Y_1=Y_{Y_0,r}$ and $W_1=W_{Y_0,r}$ be the
points given by the
event $A_1(Y_0,r,k)$. Similarly write $Y_{j+1}$ and $W_{j+1}$ for the
points given by the event $A_1(Y_j,k,r)$
for $j \ge1$, and continue until we have for some $N=N_y$ that $W_N
\notin B_E(x,3 R_2/2)$.
Note that both $\dS$ and $\dU$ are monotone on the path $\gam_y$, in
the sense that if
$x_1, x_2 \in\gam_y$ and $x_3 \in\gam(x_1,x_2)$ then for $\rho= \dS$
or $\rho=\dU$
then $\rho(x_1, x_3) \le\rho(x_1, x_2)$.
This is immediate for $\dU$ and easily proved from the definition of
$\dS$.

The construction of the $(Y_j, W_j)$ gives that
\begin{eqnarray*}
\frac{r}{a^2} &\le&\dS(Y_j, W_{j+1}) \le r,\qquad
\dS(Y_j, W_j) \le\frac{2r}{k},
\\
k^{-1/4} r^\kappa&\le&\dU(Y_j, W_{j+1}) \le
k^{1/4} r^\kappa, \qquad\dU (Y_j, W_j) \le c
(r/k)^{\kappa}.
\end{eqnarray*}
Thus, we have
\begin{eqnarray*}
\dS(W_j, W_{j+1}) &\le&\dS(W_j,
Y_j) + \dS(Y_j,W_{j+1}) \le \frac
{2r}{k}
+ r = \half b_2 r,
\\
\dS(W_j, W_{j+1})&\ge& \dS(Y_j,W_{j+1})
- \dS(W_j, Y_j)\ge r/a^2- \frac
{2r}{k}
= b_1 r.
\end{eqnarray*}
Here, we have used the equations above to define $b_1$ and $b_2$.
Similarly, we have
\begin{eqnarray*}
\dU(W_j, W_{j+1})&\leq&\dU(Y_j,W_{j+1})
\leq k^{1/4} r^\kappa,
\\
\dU(W_j, W_{j+1})&\geq&\dU(Y_j,
W_{j+1})- \dU(Y_j, W_j) \ge k^{-1/4}
r^\kappa - c (r/k)^{\kappa}\ge\tfrac{1}{2} k^{-1/4}
r^\kappa.
\end{eqnarray*}

Let $x_1, x_2 \in\gam_y[0, \tau_{x,3R_2/2}]$. We can assume that $x_1
\in\gam(y,x_2)$.
Let $j= \min\{i: W_i \in\gam(x_1, \infty)\}$. If $x_2 \in\gam(x_1, W_{j+1})$,
then $\dS(x_1,x_2) \le\dS(W_{j-1},W_j) + \dS(W_j, W_{j+1}) \le b_2r$.
So if $\dS(x_1,x_2) > b_2 r$,
then both $W_j$ and $W_{j+1}$ are on the path $\gam(x_1,x_2)$, and so
$\dU(x_1,x_2) \ge\half k^{-1/4} r^\kappa$, proving (b). Similarly, if
both $W_j$ and $W_{j+1}$ are on the path $\gam(x_1,x_2)$, then
we have $\dS(x_1,x_2) \ge b_1 r$. So if $\dS(x_1,x_2) < b_1 r$, then
$W_{j+1}\in\gam(x_2,\infty)$, and
hence $\dU(x_1,x_2) \le2 k^{1/4} r^\kappa$. \end{pf}

We now extend this result to all paths $\gam_x$ in a ball.

\begin{lem}\label{L:pathball}
Let $r \ge1$, $k \ge8$, $x_0 \in\bZ^2$, $M_i, R_i$, and $b_1,b_2$ be
as in Lem\-ma~\ref{L:LP1}. Then there exist constants $b_3, b_4$
(depending on $k$) and an event $A_3=A_3(x_0,r,k)$ with
\begin{equation}
\label{e:pacub2} \bP\bigl(A_3^c\bigr) \le c_1
\exp\bigl( - c_2k^{1/8}\bigr),
\end{equation}
such that on $A_3$ the following holds for every $x \in B_E(x_0,R_1)$:
\begin{longlist}[(a)]
\item[(a)] $\gam_x[ \tau_{x_0, R_2},\infty] \cap B_E(x_0,4R_1) =\varnothing$;

\item[(b)] If $x_1, x_2 \in\gam_x[0, \tau_{x_0,R_2}]$ and $\dS(x_1, x_2) >
b_3 r$, then
$\dU(x_1, x_2) \ge\half k^{-1/4} r^\kappa$;

\item[(c)] If $x_1, x_2 \in\gam_x[0, \tau_{x_0,R_2}]$ and $\dS(x_1, x_2) <
b_1 r$, then
$\dU(x_1, x_2) \le b_4 k^{1/4} r^\kappa$;

\item[(d)] If $x_1, x_2 \in B_E(x_0, R_1)$ and $\dS(x_1, x_2) > 2 b_3 r$, then
$\dU(x_1, x_2) \ge\half k^{-1/4} r^\kappa$;

\item[(e)] If $x_1, x_2 \in B_E(x_0, R_1)$ and $\dS(x_1, x_2) < b_1 r$, then
$\dU(x_1, x_2) \le2b_4 k^{1/4} r^\kappa$.
\end{longlist}
\end{lem}

\begin{pf} We begin by choosing a set $D_0$ which satisfies the conditions
of Lemma~\ref{L:LP1}.
Let $A_2(x_0,r,k)$ be the event defined in that lemma, and let $U_0$ be
the random tree obtained by applying Wilson's algorithm with initial
points in $D_1=D_0 \cap B(x_0,2R_1)$. Let $z \in B_E(x_0,R_1)$. We now
apply the filling in Lemma~\ref{L:fillin} to $B_E(z,r)$ taking $\delta
= 1/18k$. Let $G(z)$ be the ``good'' event given by the lemma; we have
\begin{equation}
\label{e:pacub3} \bP\bigl( G(z)^c\bigr) \le c \exp\bigl( -c
k^{1/3}\bigr).\vadjust{\goodbreak}
\end{equation}
Now choose $z_i$, $i=1, \ldots, N$ so that $N \le c M_1^2$ and
$B_E(x_0, R_1) \subset\bigcup_i B_E(z_i, r/4)$, and let
$ A_3= A_2(x_0,r,k) \cap(\bigcap_{i=1}^N G(z_i))$. The bound \eqref
{e:pacub2} then follows from \eqref{e:pacub1} and \eqref{e:pacub3}.

Let $x \in B(x_0, R_1)$, and let $W_x$ be the point where $\gam_x$
first hits the tree $U_0$. Since $G(z_i)$ holds for some $z_i$ with
$d_E(x,z_i) \le r/4$, we have by Lemma~\ref{L:fillin} that
$\dS(x,W_x) \le c k^{-1/2} r$, $\dU(x,W_x)\le c ( k^{-1/2} r)^\kappa$.
Since $U_0 =\bigcup_{y\in D_1}\gam_y$ there must exist a $y\in D_1$ such
that $W_x \in\gam_y$. Let $W_j$ be the points given in the proof of
Lemma~\ref{L:LP1}. By property Lemma~\ref{L:LP1}(a), we have that $\gam
_y$ does not return to $B_E(x_0,4R_1)$ after leaving $B_E(x_0, R_2)$
and, therefore, there exists $j$ such that $W_x \in\gam(W_{j-1},
W_j)$. (We take $W_x = W_j$ if $W_x$ is one of the points $W_i$.) Note
also that property (a) of $\gam_x$ follows from the same property for
$\gam_y$.

Let $x_1, x_2$ be on the path $\gam_x[0, \tau_{x_0, R_2}]$; we can
assume that $x_1 \in\gam(x,x_2)$. If $x_1 \in\gam(W_x, \infty)$ then both
$x_1$ and $x_2$ are in $\gam_y$, and so properties (b) and (c) follow
from Lemma~\ref{L:LP1}. So suppose that $x_1\in\gam(x,W_x)$. If $x_2
\in\gam(x, W_{j+1})$, then
\begin{eqnarray*}
\dS(x_1,x_2) &\le&\dS(x,W_x) +
\dS(W_x, W_{j+1})
\\
&\le& \dS(x,W_x) + \dS(W_{j-1}, W_j) +
\dS(W_j, W_{j+1})
\\
&\le& c k^{-1/2} r +b_2 r \leq(c+b_2)r =
b_3 r.
\end{eqnarray*}
So if $\dS(x_1,x_2) >b_3r$, then $x_2\in\gam(W_{j+1},\infty)$, and hence
$\dU(x_1,x_2) \ge\dU(W_j,\break   W_{j+1}) \ge\half k^{-1/4} r^\kappa$.
Similarly, if $x_2\in\gam(W_{j+1},\infty)$, then
$\dS(x_1,x_2) \ge\dS(W_j,\break  W_{j+1}) \ge b_1r$. So if $\dS(x_1,x_2) < b_1r$,
then $x_2 \in\gam(x, W_{j+1})$, and so
\begin{eqnarray*}
\dU(x_1,x_2) &\le&\dU(x,W_x) +
\dU(W_{j-1}, W_j) + \dU(W_j, W_{j+1})
\\
&\le& c \bigl(k^{-1/2} r\bigr)^\kappa+ 2 k^{1/4}
r^\kappa\le b_4k^{1/4} r^\kappa.
\end{eqnarray*}
This proves properties (b) and (c) of $\gam_x$.

Finally, let $x_1, x_2 \in B_E(x_0, R_1)$, and let $W$ be the point
where $\gam_{x_1}$ and $\gam_{x_2}$ meet. If $\dS(x_1,x_2) > 2b_3r$ and
$W\in\gamma_{x_1}[0,\tau_{x_0, R_2}]\cap\gamma_{x_2}[0,\tau_{x_0,
R_2}]$, then we have $\max_i \dS(x_i,W) > b_3 r$, and so $\dU(x_1,x_2)
\ge\max_i \dU(x_i,W) \ge\half k^{-1/4} r^\kappa$. If, on the other
hand, $\dS(x_1,x_2) > 2b_3r$ and $W\notin\gamma_{x_i}[0,\tau_{x_0,
R_2}]$ for either $i=1$ or $i=2$, then set $W'=\gamma_{x_i}(\tau_{x_0,
R_2})$ for the relevant $i$. Note that $W'\in\gamma(x_1,x_2)\cap\gamma
_{x_i}[0,\tau_{x_0, R_2}]$ and $d_\sU^S(x_i,W')\geq b_3r$, and so $\dU
(x_1,x_2)\geq d_\sU(x_i,W')\geq \half k^{-1/4} r^\kappa$ in this case
as well. Similarly, if $\dS(x_1,x_2) < b_1 r$, then necessarily we have
$W\in\gamma_{x_1}[0,\tau_{x_0, R_2}]\cap\gamma_{x_2}[0,\tau_{x_0,
R_2}]$ and $\max_i \dS(x_i,W) < b_1 r$, which implies $\dU(x_1,x_2) \le
2b_4 k^{1/4} r^\kappa$. This proves properties (d) and (e).
\end{pf}

Note that there is a gap between the conditions (d) and (e) above. We
could fill this by a direct calculation, but instead we will handle
this in the next result by varying $r$.

\begin{propn} \label{thmm:keycomp}
Let $r \ge1$, $\lam\ge\lam_0$ (where $\lam_0$ is a large, finite
constant), $x_0 \in\bZ^2$, and $R = r e^{c_1 \lam^{1/2}}$.
There exists an event $A_4$ with $\bP(A_4^c) \le c \exp( -c_2 \lam^{1/2})$
such that on $A_4$, for all $x,y \in B_E(x_0,R)$,
%
\begin{eqnarray}
\label{e:dudsc} \lam^{-1} \dS(x,y)^\kappa&\le&\dU(x,y) \le\lam
\dS(x,y)^\kappa \qquad\mbox{if } r\le \dS(x,y) \le R,\nonumber
\\
\dU(x,y) &\le& \lam r^\kappa \qquad\mbox{if } \dS(x,y) \le r,
\\
\dU(x,y)& \ge& \lam^{-1} R^\kappa \qquad\mbox{if } \dS(x,y) \ge R.
\nonumber
\end{eqnarray}
\end{propn}

\begin{pf} Choose $k = c \lam^4$, let $m$ be such that $2^{m-1}< \exp(
k^{1/8}/27)\le2^m$, and
define $A_4 = \bigcap_{i=0}^m A_3(x_0, 2^i r,k)$. Then $\bP(A_4^c) \le\exp
( -c k^{1/8}) \le\break \exp(-c' \lam^{1/2})$.
Now let $x,y \in B_E(x_0, R)$, and suppose $r'=\dS(x,y)\leq R$. Then
choosing the largest $i\in\{0,1,\ldots,m\}$
so that $r' \ge2b_3 2^i r$, we have $\dU(x_1,x_2) \ge c \lam^{-1} (
2^i r)^\kappa\ge c \lam^{-1} (r')^\kappa$.
Similarly, we have $\dU(x_1,x_2) \le c \lam(r')^\kappa$. Replacing $c
\lam$ by $\lam$ this
gives~\eqref{e:dudsc}, and the other two inequalities follow. \end{pf}

One consequence of the above proposition is the following approximation
result, which shows that
if a set of points is an $r/18k^2$-cover in the Euclidean metric, then
it is also a cover with respect
to the metrics $d_\sU^S$ and $d_\sU$.

\begin{propn}
\label{P:AN2FI} Let $r\geq k\geq1$. Define
$R_1:=re^{k^{1/32}}$, $R_2:=re^{k^{1/16}}$,
and suppose $D_2\subseteq\mathbb{Z}^2$ satisfies
%
\begin{equation}
\label{inclusion} B_E(0,6R_2)\subseteq\bigcup
_{x\in D_2}B_E\bigl(x,r/18k^{2}\bigr).
\end{equation}
Then there exists an event $A_5=A_5(r,k)$ such that $\mathbf
{P}(A_5^c)\leq c_1e^{-c_2k^{1/16}}$ and on $A_5$ the following holds:
%
\begin{eqnarray}
\label{schramm2fill} \max_{x\in B_E(0,R_1)}d_\sU^S
(x,D_2)&\leq&\frac{2r}{k},
\\
\label{intrinsic2fill} \max_{x\in B_E(0,R_1)}d_\sU(x,D_2)
&\leq&\frac{4r^{\kappa}}{k^{1/4}}.
\end{eqnarray}
\end{propn}

\begin{pf} First, choose a subset $D_2'\subseteq D_2$ such that (\ref
{inclusion}) holds when $D_2$ is replaced by $D_2'$ and also
$|D_2'|\leq ck^4e^{2k^{1/16}}$. Set $A'(r,k):=\bigcap_{x\in
D_2'}A_1(x,r/k,k)$, where $A_1$ is defined in the statement of Lemma~\ref{L:basicLP}. From that result, we know that
%
\begin{equation}
\label{adash} \mathbf{P}\bigl(A'^c\bigr)\leq
ck^4e^{2k^{1/16}} \mathbf{P}\bigl(A_1(0,r/k,k)\bigr)
\leq ce^{-c k^{1/8}}.
\end{equation}
Moreover, if $A'$ holds, then for $x\in B_E(0,2R_1)\cap D_2'$ we can
define $(W_j,Y_j)_{j=0}^N$ similarly to the proof of Lemma~\ref{L:LP1}.
In particular, set $W_0=Y_0=x$, and let $W_j,Y_j$ be given by the event
$A_1(Y_{j-1},r/k,k)$, up to $j=N:=\inf\{m:d_E(x,W_m)>2R_2\}$. By
construction, it follows that
%
\begin{eqnarray}
\label{upperz} \max_{z\in\gamma_x(0,\tau_{x,R_2})}d_\sU^S(z,D_2)
&\leq&\max_{j=1,\ldots
,N}d_\sU^S(W_{j-1},W_j)
\nonumber
\\[-8pt]
\\[-8pt]
\nonumber
&\leq&\max_{j=1,\ldots,N}d_\sU^S(Y_{j-1},W_j)
\leq \frac{r}{k}.
\end{eqnarray}
Next, choose $D_2''\subseteq D_2'\cap B_E(0,2R_1)$ such that $
B_E(0,R_1)\subseteq\bigcup_{x\in D_2''}B_E(x,r/\break 18k^2)$ and $|D_2''|\leq c
k^4e^{2k^{1/32}}$. Set
$A''(r,k):=A'(r,k)\cap(\bigcap_{x\in D_2''}B(x,r,k))$, where $B(x,r,k):=\{
\gamma_x(\tau_{x,R_2},\infty)\cap B_E(0,2R_1)=\varnothing\}$. By applying
Lemma~\ref{L:pathsize} in conjunction with (\ref{adash}), we obtain
\[
\mathbf{P}\bigl(A''^c\bigr)\leq
ce^{-c k^{1/8}}+c k^4e^{2k^{1/32}} \mathbf {P}\bigl(
\gamma_0(\tau_{0,R_2},\infty)\cap B_E(0,4R_1)
\neq\varnothing\bigr)\leq c_1e^{-c_2k^{1/16}}.
\]
Define $U_0$ to be the subtree of $\sU$ spanned by $D_2''$ and suppose
$A''$ holds. If $x\in U_0\cap B_E(0,2R_1)$, then it must be the case
that $x\in\gamma_y(0,\tau_{y,R_2})$ for some $y\in D_2''$. Hence, by
(\ref{upperz}), it holds that $\max_{x\in U_0\cap B_E(0,2R_1)}d_\sU
^S(x,D_2) \leq r/k$. Now, by applying Lemma~\ref{L:fillin} with root
$U_0$, it is possible to deduce
\[
\mathbf{P} \biggl(\max_{x\in B_E(0,R_1)}d_{\sU}^S(x,U_0)>
\frac
{r}{k} \biggr)\leq Ce^{-c e^{k^{1/32}}}.
\]
So, if $A'''$ is defined to be the event that both $A''$ and 
$\max_{x\in B_E(0,R_1)}d_{\sU}^S(x,\break U_0)\leq r/k$ hold, then we have $\mathbf
{P}(A'''^c)\leq c_1e^{-c_2k^{1/16}}$ and also (\ref{schramm2fill})
holds on $A'''$.

To complete the proof, we will use Proposition~\ref{thmm:keycomp}
with $(x_0,r,\lambda)$ given by $(0,2r/k, k)$ to compare the relevant distances.
Since $R=2rk^{-1}e^{c_1k^{1/2}}\geq2R_1$ for large $k$, we find that
with probability
exceeding $1-ce^{-c_2k^{1/2}}$ it is the case that
\[
\mathop{\max_{x,y\in B_E(0,2R_1):}}_{d_{\sU}^S(x,y)
\le2r/k}d_{\sU}(x,y)\leq
k \biggl(\frac{2r}{k} \biggr)^{\kappa}\leq \frac{4r^{\kappa}}{k^{1/4}}.
\]
Note that if $A'''$ and the above inequality both hold, then so does
(\ref{intrinsic2fill}). Hence, in conjunction with the conclusion of
the previous paragraph, this completes the proof.
\end{pf}

We can now improve the volume estimates of \cite{BM11}. Recall from
\eqref{e:dfdef} that
$d_f=2/\kappa=8/5$, and define, for $\lambda, n\ge1$,
\begin{eqnarray*}
\tilde A(\lambda, n)&:=&\bigl\{\omega: \lambda^{-1}R^{d_f}
\le\bigl|B_\sU(x,R)\bigr|\le \lambda R^{d_f}\\
&& \mbox{for all } x\in
B_E(0, n), R\in\bigl[e^{-\lambda^{1/40}}n^{\kappa},n^{\kappa}
\bigr]\bigr\}.
\end{eqnarray*}
The following result extends a fundamental estimate of \cite{BM11};
the key improvement is that the upper bound does not depend on $n$
(once $n$ is suitably large). Although we do not need to do so here, we
note that the same approach can also be
used to obtain a similar improvement of the resistance estimates in
\cite{BM11}.

\begin{propn}\label{P:keypropals}
There exist constants $c_1,c_2\in(0,\infty)$ such that
\[
\bP \bigl(\tilde A(\lambda, n)^c \bigr) \le c_1\exp
\bigl(-c_2\lambda^{1/80}\bigr)\qquad \mbox{for all } n\geq
e^{\lambda^{1/16}}.
\]
\end{propn}

\begin{pf}
Let $k=\lam$, $r=ne^{-\lam^{1/32}}$, and let
$R_1=n$, $R_2= r e^{k^{1/6}}$ and $D_2$ be as in Proposition~\ref
{P:AN2FI}, with
$|D_2|\leq ck^4e^{2k^{1/16}}$.
Set $m_0:=\inf\{m:k^m\geq e^{k^{1/32}} \}$.
Let $A_5(r,k)$
be the event given in the statement of Proposition~\ref{P:AN2FI}, and
\[
E(r,k):=\bigcap_{x\in D_2}\bigcap
_{m=1}^{m_0+1}\bigl\{k^{-1}
\bigl(rk^m\bigr)^{\kappa} \leq\bigl|B_\sU\bigl(x,
\bigl(rk^m\bigr)^{\kappa}\bigr)\bigr|\leq k\bigl(rk^m
\bigr)^{\kappa}\bigr\}.
\]
A simple union bound allows us to deduce from
Theorem~\ref{vollemust}(b) that
\[
\mathbf{P} \bigl(E(r,k)^c \bigr)\leq Ck^4e^{2k^{1/16}}
k^{1/32} ce^{-k^{1/9}}\leq C e^{-ck^{1/9}}.
\]
Consequently, we have $\mathbf{P}( E(r,k)^c \cup A_5(r,k)^c) \le c \exp
(-c \lam^{1/16} )$.

Suppose that $E(r,k)\cap A_5(r,k)$ holds. Let $x\in B_E(0,n)$, and $s
\in[rk^3, n]$.
Choose $m \in\{3, \ldots, m_0+1\}$ such that $s\in[rk^m,rk^{m+1})$.
Since $A_5(r,k)$ holds, there exists $y\in D_2$ with $d_\mathcal
{U}(x,y)\leq4r^{\kappa}/k^{1/4}$.
Hence,
\begin{eqnarray*}
\bigl|B_\sU \bigl(x,s^{\kappa} \bigr)\bigr|&\leq&\bigl\llvert
B_\sU \bigl(y,\bigl(rk^{m+1}\bigr)^{\kappa}+4r^{\kappa}/k^{1/4}
\bigr)\bigr\rrvert \leq\bigl\llvert B_\sU \bigl(y,
\bigl(rk^{m+2}\bigr)^{\kappa} \bigr)\bigr\rrvert \leq k
\bigl(rk^{m+2}\bigr)^2\\
&\leq& k^5 s^2.
\end{eqnarray*}
Similarly, $|B_\sU(x,s^{\kappa})|\geq k^{-5} s^2$. Since $(rk^3)^\kappa
\le n^\kappa\exp(-\lam^{1/40})$
it follows that $E(r,k)\cap A_5(r,k) \subset\tilde{A}(\lambda^5,n)$,
which completes the proof of the proposition. \end{pf}

From this, we can prove the following distributional measure bounds,
which will be used in the proof of Theorem~\ref{main12}(b)(ii).

\begin{cor}\label{P:measprop} Given $R>0$, there exist constants
$c_1,\ldots, c_7 \in(0,\infty)$ (depending on $R$)
such that for every $r\in(0,c_7)$,
%
\begin{eqnarray}
\label{measupperbound} &&\limsup_{\delta\rightarrow0} \mathbf{P} \Bigl(
\delta^{2} \min_{x\in B_E(0,\delta^{-1}R)} \mu_\sU
\bigl(B_\sU\bigl(x,\delta^{-\kappa} r\bigr) \bigr) \leq
c_1 r^{d_f}\bigl(\log r^{-1}\bigr)^{-80}
\Bigr)
\nonumber
\\[-8pt]
\\[-8pt]
\nonumber
&&\qquad\leq c_2r^{c_3},
\\
\label{measlowerbound} &&\limsup_{\delta\rightarrow0} \mathbf{P} \Bigl(
\delta^{2} \max_{x\in B_E(0,\delta^{-1}R)} \mu_\sU
\bigl(B_\sU\bigl(x,\delta^{-\kappa} r\bigr) \bigr) \geq
c_4 r^{d_f}\bigl(\log r^{-1}\bigr)^{80}
\Bigr)
\nonumber
\\[-8pt]
\\[-8pt]
\nonumber
&&\qquad\leq c_5r^{c_6}.
\end{eqnarray}
\end{cor}

\begin{pf} We just prove \eqref{measlowerbound}; the proof of
\eqref{measupperbound} is similar. Fix $R\geq1$, and suppose $r\in
(0,1)$, $\delta\in(0,1)$.
Define $n:=\delta^{-1}R$ and $\lambda:= (\log(R^{\kappa}/r))^{80}$.
Since $\delta^{-\kappa}r\in[e^{-\lambda^{1/40}}n^{\kappa},n^{\kappa}]$,
we have that, on $\tilde{A}(\lambda,n)$,
\[
\min_{x\in B_E(0,\delta^{-1}R)} \mu_\sU \bigl(B_\sU
\bigl(x,\delta^{-\kappa} r\bigr) \bigr) \geq\lambda^{-1}\delta
^{-2}r^{d_f} \geq c_1\delta^{-2}r^{d_f}
\bigl(\log r^{-1}\bigr)^{-80}.
\]
Hence, by Proposition~\ref{P:keypropals}, the left-hand side of (\ref
{measlowerbound}) is bounded above by $Ce^{-c\lambda^{1/80}}$.
\end{pf}

Let $N_\mathcal{U}(r,s)$ the minimum number of $d_\mathcal{U}$-balls of
radius $s$
required to cover $B_\mathcal{U}(0,r)$.
Another consequence of Proposition~\ref{P:keypropals} is the following
bound on
$N_\mathcal{U}(r,r/\lam)$.

\begin{lem}\label{disclem}
There exist constants $c_1,c_2,c_3,\lambda_0\in(0,\infty)$ such that,
for $r\geq e^{\kappa(\log\lambda)^{41/16}}$ and $\lambda\geq\lambda_0$,
\[
\bP \bigl(N_\mathcal{U}(r,r/\lambda)\ge c_1(\log
\lam)^{107} \lam^{d_f} \bigr)\leq c_2
e^{-c_3 (\log\lam
)^{41/80} }.
\]
\end{lem}

\begin{pf}
Let $\th\ge1$ be such that $2\lam\le\th^{-1} \exp( \th^{1/40})$. By
Theorem~\ref{vollemust}(a), we have that
\[
\bP \bigl( B_\sU(0,r) \not\subset B_E\bigl( 0,
\th^{1/\kappa} r^{1/\kappa} \bigr) \bigr) \le e^{-c \th^{2/3} }.
\]
Now it is straightforward to check that one can cover $B_\sU(0,r)$ by
balls $B_\sU(z_i,r/\lam)$, $i=1,\ldots, M$, such that $B_\sU(z_i,r/2\lam
)$ are disjoint and $z_i \in B_\sU(0,r)$. Moreover, it is necessarily
the case that $M\geq N_\sU(r, r/\lam)$. Setting $n^\kappa= \th r$, if
$\tilde A(\th,n)$ holds and $ B_\sU(0,r) \subset B_E(0,n)$
then we have $|B_\sU(0,r)| \le(\th r)^{d_f}$ and
$|B_\sU(z_i,r/2\lam)| \ge c \th^{-1} (r/\lam)^{d_f}$ for each $i$.
Thus, we deduce from
Proposition~\ref{P:keypropals} that
\begin{equation}
\label{e:ballbound} \bP\bigl( N_\sU(r, r/\lam) \ge c \th^{1+d_f}
\lam^{d_f} \bigr) \le c \exp\bigl( - c \th^{1/80} \bigr).
\end{equation}
Taking $\th= (\log\lam)^{41}$ completes the proof. \end{pf}

\begin{rem} \label{R:covers}
Taking $\th=\lam$ in \eqref{e:ballbound} gives the bound, for $r\geq
e^{\kappa\lambda^{1/16}}$ and $\lambda$ large,
\[
\bP\bigl( N_{\sU}(r, r/\lam) \ge c \lam^{1+2d_f} \bigr) \le c
\exp\bigl( - c \lam ^{1/80} \bigr).
\]
\end{rem}

\section{Topology for UST scaling limit}\label{topsec}

In this section, we introduce the topology on measured, rooted spatial
trees for which we
prove tightness for the law of the rescaled UST.
This topology is finer than that considered in \cite{ABNW,Schramm},
since it incorporates the full convergence of real trees embedded into
Euclidean space,
rather than merely the shape of subsets spanning a finite number of vertices.
This point will be important when it comes to the proof of Theorem~\ref{main2}.

We define $\mathbb{T}$ to be the collection of quintuplets of the form
\[
\uT=(\mathcal{T},d_\mathcal{T},\mu_{\mathcal{T}},\phi_\mathcal{T},
\rho _\mathcal{T}),
\]
where: $(\mathcal{T},d_\mathcal{T})$ is a complete and locally compact
real tree (see \cite{rrt},
Definition~1.1, e.g.); $\mu_\mathcal{T}$ is a locally finite
Borel measure on
$(\mathcal{T},d_\mathcal{T})$; $\phi_\mathcal{T}$ is a continuous map
from $(\mathcal{T},d_\mathcal{T})$ into a separable metric space
$(M,d_M)$; and $\rho_\mathcal{T}$ is a distinguished vertex in $\mathcal{T}$.
[Usually the image space $(M,d_M)$ we consider is $\mathbb{R}^2$
equipped with the Euclidean distance,
though we will also consider other image spaces at certain places in
our arguments.]
We call such a quintuplet a \emph{measured, rooted, spatial tree}.
Let $\mathbb{T}_c$ be the subset of $\mathbb{T}$ for which $(\mathcal
{T},d_\mathcal{T})$ is compact.
We will say that two elements of $\mathbb{T}$, $\uT$ and $\uT'$ say,
are equivalent if there exists an isometry $\pi:(\mathcal{T},d_\mathcal
{T})\rightarrow(\mathcal{T}',d_\mathcal{T}')$ for which $\mu_\mathcal
{T}\circ\pi^{-1}=\mu_\mathcal{T}'$, $\phi_\mathcal{T}=\phi_\mathcal
{T}'\circ\pi$ and also $\pi(\rho_\mathcal{T})=\rho_\mathcal{T}'$.

In order to introduce a topology on $\mathbb{T}$, we will start by
defining a topology on~$\mathbb{T}_c$. In particular, for two elements
of $\mathbb{T}_c$, we set $\Delta_c (\uT,\uT' )$ to be equal to
%
\begin{eqnarray}
\label{deltacdef}&& \mathop{\inf_{Z,\psi,\psi',\mathcal{C}:}}_{(\rho_\mathcal{T},\rho
_\mathcal{T}')\in\mathcal{C}} \Bigl
\{d_P^Z \bigl(\mu_{\mathcal{T}}\circ
\psi^{-1},\mu_{\mathcal
{T}}'\circ\psi'^{-1}
\bigr)\nonumber
\\[-8pt]
\\[-8pt]
\nonumber
&&\hspace*{23pt}\qquad{}+ \sup_{(x,x')\in\mathcal{C}} \bigl(d_Z \bigl(\psi(x),
\psi'\bigl(x'\bigr) \bigr)+d_M \bigl(
\phi_\mathcal{T}(x),\phi_\mathcal{T}'
\bigl(x'\bigr) \bigr) \bigr) \Bigr\},
\end{eqnarray}
where the infimum is taken over all metric spaces $Z=(Z,d_Z)$,
isometric embeddings $\psi:(\mathcal{T},d_\mathcal{T})\rightarrow Z$,
$\psi':(\mathcal{T}',d_\mathcal{T}')\rightarrow Z$, and correspondences
$\mathcal{C}$ between $\mathcal{T}$ and $\mathcal{T}'$, and we define
$d_P^Z$ to be the Prohorov distance between finite Borel measures on
$Z$. Note that, by a correspondence $\mathcal{C}$ between $\mathcal{T}$
and $\mathcal{T}'$, we mean a subset of $\mathcal{T}\times\mathcal
{T}'$ such that for every $x\in\mathcal{T}$ there exists at least one
$x'\in\mathcal{T}'$ such that $(x,x')\in\mathcal{C}$ and conversely
for every $x'\in\mathcal{T}'$ there exists at least one $x\in\mathcal
{T}$ such that $(x,x')\in\mathcal{C}$.

\begin{propn}\label{sepprop} The function $\Delta_c$ defines a metric on the
equivalence classes of $\mathbb{T}_c$. Moreover, the resulting metric
space is separable.
\end{propn}

\begin{pf} The proof of this result is almost identical to that of \cite{CHK},
Lemma~2.1, taking, in the notation of that paper, $I=\{1\}$ and
$q_1(x,y):=\phi_\mathcal{T}(x)$. The main change is that when
considering a correspondence
between $\mathcal{T}$ and $\mathcal{T}'$, one has to require that the
pair of roots
$(\rho_\mathcal{T},\rho_{\mathcal{T}}')$ is included,
and, when selecting the points $x_i$, $x_i'$ as in \cite{CHK},
one should take $x_1=\rho_\mathcal{T}$ and $x_1'=\rho_\mathcal{T}'$. A
second change is that in the proof of
separability, rather than approximating by metric spaces with a finite
number of vertices,
one should approximate by real trees formed of a finite number of line segments;
however, making these changes is routine and we omit the details.
\end{pf}

\begin{rem}
Even if $(M,d_M)$ is assumed to be complete, the space of equivalence
classes of $\mathbb{T}_c$
is not complete with respect to the metric $\Delta_c$ in general.
Indeed, suppose $(M,d_M)=(\mathbb{R}^2,d_E^{(2)})$ and consider
$([0,1],d_E^{(1)},\mathcal{L},f,0)\in\mathbb{T}_c$, where $d_E^{(d)}$
is the $d$-dimensional Euclidean distance, $\mathcal{L}$ is Lebesgue
measure on $[0,1]$, and $f: [0,1] \to\bR^2$ is any continuous
nonconstant function. If we replace $d_E^{(1)}$ by $\varepsilon
d_E^{(1)}$, then the sequence of elements in $\mathbb{T}_c$ that we
obtain is Cauchy as $\varepsilon\rightarrow0$, but does not have a
limit in $\mathbb{T}_c$. One way to ensure completeness would be to
restrict to a subset of $\mathbb{T}_c$ for which the functions $\phi
_\mathcal{T}$ satisfy an equi-continuity condition.
\end{rem}

To extend $\Delta_c$ to a metric on the equivalence classes of $\mathbb
{T}$, we consider bounded restrictions of elements of $\mathbb{T}$
(cf. \cite{ADH}).
Thus, for $\uT\in\mathbb{T}$,
let $\uT^{(r)}=(\mathcal{T}^{(r)},d_\mathcal{T}^{(r)},\mu_{\mathcal
{T}}^{(r)},\phi_\mathcal{T}^{(r)},\rho_\mathcal{T}^{(r)})$ be obtained
by taking: $\mathcal{T}^{(r)}$ to be the closed ball in
$(\mathcal{T},d_\mathcal{T})$ of radius $r$ centred at $\rho_\mathcal{T}$;
$d_\mathcal{T}^{(r)}$ $\mu_\mathcal{T}^{(r)}$ and $\phi_\mathcal
{T}^{(r)}$ to be the
restriction of $d_\mathcal{T}$, $\mu_\mathcal{T}$ and $\phi_\mathcal{T}$,
respectively, to $\mathcal{T}^{(r)}$, and $\rho_\mathcal{T}^{(r)}$ to be
equal to $\rho_\mathcal{T}$. As in \cite{ADH}, the fact that $(\mathcal
{T},d_\mathcal{T})$ is a real tree, and therefore
a length space, means we can apply the Hopf--Rinow theorem (which
implies that all closed,
bounded subsets of a complete and locally compact length space are
compact) to establish
that $\uT^{(r)}$ is an element of $\mathbb{T}_c$.
Furthermore, as in \cite{ADH}, Lemma~2.8, we can check the regularity
of this restriction with
respect to the metric $\Delta_c$.

\begin{lem}\label{cadlag}For any two elements of $\mathbb{T}$, $\uT$ and $\uT'$,
the function $r\mapsto\Delta_c(\uT^{(r)},\uT'^{(r)})$ is cadlag.
\end{lem}

\begin{pf}
By considering the natural embedding of $\mathcal{T}^{(r)}$ into
$\mathcal{T}^{(r+\varepsilon)}$,
along with the correspondence consisting of pairs $(x,x')$ such that
$x$ is the closest point in
$\mathcal{T}^{(r)}$ to $x'\in\mathcal{T}^{(r+\varepsilon)}$, we have,
as in \cite{ADH}, Lemma~5.2, that
\[
{\Delta_c \bigl(\uT^{(r)},\uT^{(r+\varepsilon)} \bigr)}\leq
\mu_\mathcal {T} \bigl(\mathcal{T}^{(r+\varepsilon)}\setminus
\mathcal{T}^{(r)} \bigr)+\varepsilon+\mathop{\sup_{x,x'\in\mathcal{T}^{(r+\varepsilon)}:}}_{
d_\mathcal{T}(x,x')\leq\varepsilon}\,d_M
\bigl(\phi_\mathcal{T}(x),\phi _\mathcal{T}\bigl(x'
\bigr) \bigr);
\]
given this, the proof is a straightforward adaption of the proof of
\cite{ADH}, Lemma~2.8.
\end{pf}

This result allows us to well define a function $\Delta$ on $\mathbb
{T}^2$ by setting
%
\begin{equation}
\label{delta} {\Delta \bigl(\uT,\uT' \bigr)}:=\int
_0^\infty e^{-r} \bigl(1\wedge
\Delta_c \bigl(\uT^{(r)},\uT'^{(r)}
\bigr) \bigr) \,dr.
\end{equation}

\begin{propn}\label{deltametric} The function $\Delta$ defines a metric on
the equivalence classes of $\mathbb{T}$. Moreover, the resulting metric
space is separable.
\end{propn}

\begin{pf} Again, the proof is similar to the corresponding result in \cite
{ADH}. Positivity, finiteness and symmetry of $\Delta$ are clear.
Moreover, the triangle inequality is easy to check from the definition
and the fact that the triangle inequality holds for $\Delta_c$. So, to
establish that $\Delta$ is a metric, it remains to prove positive
definiteness. To this end, suppose that $\uT$ and $\uT'$ are such that
the expression at (\ref{delta}) is equal to zero. From Lemma~\ref
{cadlag}, it follows that
$\Delta_c(\uT^{(r)},\uT'^{(r)})=0$ for every $r>0$. Consequently, for
each $r$, there exists an isometry
$\pi_r:(\mathcal{T}^{(r)},d_\mathcal{T}^{(r)})\rightarrow(\mathcal
{T'}^{(r)},d_\mathcal{T}'^{(r)})$ such that
$\mu_\mathcal{T}^{(r)}\circ\pi_r^{-1}=\mu_{\mathcal{T}}'^{(r)}$, $\phi
^{(r)}_\mathcal{T}=\phi_\mathcal{T}'^{(r)}\circ\pi_r$ and also $\pi
_r(\rho^{(r)}_\mathcal{T})=\rho_\mathcal{T}'^{(r)}$. For $n,k\geq1$,
let $(x^{n,k}_i)_{i=1}^{N(n,k)}$ be a finite $k^{-1}$-cover of $\mathcal
{T}^{(n)}$ containing the root $\rho_\mathcal{T}$ (such a collection
exists as a result of the compactness of $\mathcal{T}^{(n)}$). Since
$\pi_r$ is an isometry, we have that $(\pi_m(x^{n,k}_i))_{m\geq n}$ is
a bounded sequence for each $n,k\geq1$ and $1\leq i\leq N(n,k)$, and
so has a convergent subsequence. By a diagonal procedure, one can thus
find a subsequence $(m_j)_{j\geq1}$ such that $\pi(x^{n,k}_i)=\lim_{j\rightarrow\infty}\pi_{m_j}(x^{n,k}_i)$ exists for every $n,k\geq1$
and $1\leq i\leq N(n,k)$. From this construction, we obtain that $\pi$
is distance-preserving on $\{x^{n,k}_i: n,k\geq1, 1\leq i\leq
N(n,k)\}$ and, since the latter set is dense in $\mathcal{T}$, we can
extend it to a distance-preserving map on $\mathcal{T}$. Clearly, by
reversing the roles of $\mathcal{T}$ and $\mathcal{T}'$, it is also
possible to find a distance-preserving map from $\mathcal{T}'$ to
$\mathcal{T}$. Hence, $\pi$ must be an isometry. Moreover, it is clear
that this map is root-preserving, that is, $\pi(\rho_\mathcal{T})=\rho
_\mathcal{T}'$. To check that it is measure-preserving, that is, $\mu
_\mathcal{T}\circ\pi^{-1}=\mu_{\mathcal{T}}'$, one can follow an
identical argument to that applied in the proof of \cite{ADH}, Proposition~5.3, based on considering approximations to the measures $\mu
_\mathcal{T}^{(n)}$ and $\mu_{\mathcal{T}}'^{(n)}$ supported on
$(x^{n,k}_i)_{i=1}^{N(n,k)}$ and $(\pi(x^{n,k}_i))_{i=1}^{N(n,k)}$,
respectively. Finally, we note that the continuity of $\phi_{\mathcal
{T}}'$ implies
\[
\phi_\mathcal{T}'\bigl(\pi\bigl(x^{n,k}_i
\bigr)\bigr)=\lim_{j\rightarrow\infty} \phi_\mathcal{T}'^{(m_j)}
\circ\pi_{m_j}\bigl(x^{n,k}_i\bigr)= \lim
_{j\rightarrow\infty} \phi_\mathcal{T}^{(m_j)}
\bigl(x^{n,k}_i\bigr)=\phi_\mathcal{T}
\bigl(x^{n,k}_i\bigr).
\]
Since $\phi_\mathcal{T}$ is also continuous, it follows that $\phi
_\mathcal{T}=\phi_\mathcal{T}'\circ\pi$. Hence, we have shown that $\uT
$ and $\uT'$ are equivalent, and so $\Delta$ is indeed a metric on the
equivalence classes of $\mathbb{T}$.

For separability, we first note that $\Delta(\uT,
\uT^{(r)})\leq e^{-r}$, and so $\mathbb{T}_c$ is dense in $(\mathbb
{T},\Delta)$. Since $(\mathbb{T}_c,\Delta_c)$ is separable, it will
thus be sufficient to check that convergence in $(\mathbb{T}_c,\Delta
_c)$ implies convergence in $(\mathbb{T},\Delta)$ (cf. \cite{ADH},
Proposition~2.10). So let us start by supposing that we have a sequence
$\uT_n$ that converges to $\uT$ in $(\mathbb{T}_c,\Delta_c)$. In
particular, we can find a sequence of metric spaces $Z_n$, isometric
embeddings $\psi_n:\mathcal{T}\rightarrow Z_n$, $\psi'_n:\mathcal
{T}_n\rightarrow Z_n$ and correspondences $\mathcal{C}_n$ between
$\mathcal{T}$ and $\mathcal{T}_n$ containing $(\rho_\mathcal{T},\rho
_{\mathcal{T}_n})$ such that
%
\begin{eqnarray}
\label{yy}&& d_P^{Z_n}\bigl(\mu_{\mathcal{T}}\circ
\psi_n^{-1},\mu_{\mathcal{T}_n}\circ\psi
_n'^{-1}\bigr) \nonumber\\
&&\quad{}+\sup_{(x,x')\in\mathcal{C}_n}
\bigl(d_{Z_n}\bigl(\psi_n(x),\psi_n'
\bigl(x'\bigr)\bigr)+ d_M \bigl(\phi_{\mathcal{T}}(x),
\phi_{\mathcal{T}_n}\bigl(x'\bigr) \bigr) \bigr)\\
&&\qquad<
\varepsilon_n,\nonumber
\end{eqnarray}
where $\varepsilon_n\rightarrow0$. Now, define $\psi_n^{(r)}$ to be
the restriction of $\psi_n$ to $\mathcal{T}^{(r)}$, ${\psi_n'}^{(r)}$
to be the restriction of $\psi'_n$ to $\mathcal{T}_n^{(r)}$, and
$\mathcal{C}_n^{(r)}$ to be the collection of pairs $(x,x')$ such that:
either $x\in\mathcal{T}^{(r)}$ and $x'$ is the closest point in
$\mathcal{T}^{(r)}_n$ to an element $x''\in\mathcal{T}_n$ such that
$(x,x'')\in\mathcal{C}_n$; or $x'\in\mathcal{T}_n^{(r)}$ and $x$ is the
closest point in $\mathcal{T}^{(r)}$ to an element $x''\in\mathcal{T}$
such that $(x'',x')\in\mathcal{C}_n$. Note that $\psi_n^{(r)}$ and
${\psi'_n}^{(r)}$ are isometric embeddings of $\mathcal{T}^{(r)}$ and
$\mathcal{T}_n^{(r)}$, respectively, into $Z_n$, and that $\mathcal
{C}_n^{(r)}$ is a correspondence between $\mathcal{T}^{(r)}$ and
$\mathcal{T}_n^{(r)}$ such that $(\rho_{\mathcal{T}}^{(r)},\rho
_{\mathcal{T}_n}^{(r)})\in\mathcal{C}_n^{(r)}$. If we suppose that $x\in
\mathcal{T}^{(r)}$ and $x'$ is the closest point in $\mathcal
{T}^{(r)}_n$ to an element $x''\in\mathcal{T}_n$ such that $(x,x'')\in
\mathcal{C}_n$, then
\begin{eqnarray*}
d_{\mathcal{T}_n}\bigl(\rho_{\mathcal{T}_n},x''
\bigr)&\leq& d_{Z_n}\bigl(\psi_n'(
\rho_{\mathcal{T}_n}),\psi_n(\rho_{\mathcal{T}})\bigr)
+d_{Z_n}\bigl(\psi_n(\rho_{\mathcal{T}}),
\psi_n(x)\bigr)\\
&&{}+ d_{Z_n}\bigl(\psi_n(x),
\psi'_n\bigl(x''\bigr)
\bigr),
\end{eqnarray*}
which is bounded above by $r+2\varepsilon_n$. It follows that
$d_{\mathcal{T}_n}(x',x'')<2\varepsilon_n$ and, therefore, also
$d_{Z_n}(\psi_n(x),\psi'_n(x'))<3\varepsilon_n$. A similar argument
applies to the case when $x'\in\mathcal{T}_n^{(r)}$ and $x$ is the
closest point in $\mathcal{T}^{(r)}$ to an element $x''\in\mathcal{T}$
such that $(x'',x')\in\mathcal{C}_n$. Consequently, we obtain that
%
\begin{equation}
\label{yyy} \sup_{(x,x')\in\mathcal{C}_n^{(r)}}\,d_{Z_n}\bigl(
\psi^{(r)}_n(x),{\psi _n'}^{(r)}
\bigl(x'\bigr)\bigr)<3\varepsilon_n.
\end{equation}
From this, one can proceed as in the proof of \cite{ADH}, Proposition~2.10, to deduce that
\[
d_P^{Z_n} \bigl(\mu_{\mathcal{T}}^{(r)}\circ
\bigl(\psi^{(r)}_n\bigr)^{-1},\mu
_{\mathcal{T}_n}^{(r)}\circ \bigl({\psi_n'}^{(r)}
\bigr)^{-1} \bigr)<\varepsilon_n+\mu_{\mathcal{T}} \bigl(
\mathcal{T}^{(r+4\varepsilon_n)}\setminus \mathcal{T}^{(r-4\varepsilon_n)} \bigr).
\]
Moreover, it is also elementary to deduce from (\ref{yy}) and (\ref
{yyy}) that
\[
\sup_{(x,x')\in\mathcal{C}_n^{(r)}}\bigl\llvert \phi_{\mathcal{T}}^{(r)}(x)-
\phi_{\mathcal{T}_n}^{(r)}\bigl(x'\bigr)\bigr\rrvert \leq
\varepsilon_n+ \mathop{\sup_{(x,x')\in\mathcal{T}^{(r+4\varepsilon_n)}:}}_{d_\mathcal
{T}(x,x')<4\varepsilon_n}d_M
\bigl(\phi_{\mathcal{T}}(x), \phi_{\mathcal{T}}\bigl(x'\bigr)
\bigr).
\]
Hence, we have established that
%
\begin{eqnarray}
\label{extra} \Delta_c \bigl(\uT_n^{(r)},
\uT^{(r)} \bigr)& \leq&5 \varepsilon_n+\mu _{\mathcal{T}}
\bigl(\mathcal{T}^{(r+4\varepsilon_n)}
\setminus \mathcal{T}^{(r-4\varepsilon_n)} \bigr)
\nonumber
\\[-8pt]
\\[-8pt]
\nonumber
&&{}+
\mathop{\sup_{(x,x')\in
\mathcal{T}^{(r+4\varepsilon_n)}:}}_{d_\mathcal{T}(x,x')<4\varepsilon
_n}d_M \bigl(
\phi_{\mathcal{T}}(x), \phi_{\mathcal{T}}\bigl(x'\bigr) \bigr).
\end{eqnarray}
Since $\mu_\mathcal{T}$ is a finite measure, this expression must
converge to zero for all but at most a countable number of values of
$r$. Thus, dominated convergence implies that $\Delta(\uT_n, \uT
)\rightarrow0$, as desired.
\end{pf}

Next, under the additional assumption that $(M,d_M)$ is proper
(i.e., every closed ball in $M$ is compact), we provide a sufficient
condition for a subset $\mathcal{A}$ of $\mathbb{T}$ to be relatively
compact with respect to the topology induced by $\Delta$. This extends
the corresponding result of \cite{ADH}, Theorem~2.11, to include the
spatial embedding.

\begin{lem}\label{relcompact} Suppose $(M,d_M)$ is proper. Let $\mathcal{A}$
be a subset of $\mathbb{T}$ such that, for every $r>0$:
\begin{longlist}[(iii)]
\item[(i)] for every $\varepsilon>0$, there exists a finite integer
$N(r,\varepsilon)$ such that for any element $\uT$ of $\mathcal{A}$
there is an
$\varepsilon$-cover of $\mathcal{T}^{(r)}$ of cardinality less than
$N(r,\varepsilon)$;

\item[(ii)] it holds that
\[
\sup_{\uT\in\mathcal{A}}\mu_\mathcal{T} \bigl(\mathcal{T}^{(r)}
\bigr)<\infty;
\]
\item[(iii)] $\{\phi_{\mathcal{T}}(\rho_\mathcal{T}): \uT\in\mathcal{A}\}$ is
a bounded subset of $M$, and for every $\varepsilon>0$, there exists a
$\delta=\delta(r,\varepsilon)>0$ such that
\[
\sup_{\uT\in\mathcal{A}} \mathop{\sup_{x,y\in\mathcal{T}^{(r)}:}}_{d_\mathcal{T}(x,y)\leq\delta
}d_M
\bigl(\phi_{\mathcal{T}}(x),\phi_\mathcal{T}(y) \bigr)<\varepsilon.
\]
Then $\mathcal{A}$ is relatively compact.
\end{longlist}
\end{lem}

\begin{pf} We follow closely the proof \cite{ADH}, Theorem~2.11. Suppose
that $\uT_n$ is a sequence in a set $\mathcal{A}\subseteq\mathbb{T}$
that is assumed to satisfy the properties listed in the statement of
the lemma. We can then define $U$ to be a countable index set such that
$\{x^n_u:u\in U\}$ is dense in $\mathcal{T}_n$ for each $n$ (we further
assume that $0\in U$ and $x_0^n=\rho_{\mathcal{T}_n}$), and also
introduce an abstract space $\mathcal{T}':=\{x_u:u\in U\}$ such that,
for some subsequence $(n_i)_{i\geq1}$,
%
\begin{equation}
\label{bb} d_{\mathcal{T}_{n_i}}\bigl(x_u^{n_i},x_v^{n_i}
\bigr)\rightarrow d_\mathcal{T}(x_u,x_v)
\end{equation}
for each pair of indices $u,v\in U$, where the right-hand side may be
taken as a definition of the function $d_\mathcal{T}:\mathcal{T}'\times
\mathcal{T}'\rightarrow\mathbb{R}_+$. In fact, $d_\mathcal{T}$ is a
quasi-metric on $\mathcal{T}'$, and so, with a slight abuse of
notation, we obtain a metric space $(\mathcal{T}',d_{\mathcal{T}})$ by
identifying points that are a $d_\mathcal{T}$-distance of zero apart.
Moreover, the argument of~\cite{ADH} gives us that the completion
$(\mathcal{T},d_{\mathcal{T}})$ of this metric space is locally
compact, and identifies $\rho_\mathcal{T}:=x_0$ as the root for the
space. It also describes how to construct a corresponding locally
finite Borel measure on $\mathcal{T}$, which we will call $\mu_\mathcal
{T}$. Now, from property (iii) and (\ref{bb}), it is easy to see that
$\phi_{\mathcal{T}_{n_i}}(x^{n_i}_u)$ is bounded for each~$u$, and so a
diagonal procedure yields that, by taking a further subsequence if
necessary, $\phi_{\mathcal{T}_{n_i}}(x^{n_i}_u)\rightarrow\phi_\mathcal
{T}(x_u)$, for each $u\in U$, where, similarly to the definition of
$d_\mathcal{T}$, the right-hand side provides a definition of $\phi
_\mathcal{T}(x_u)$ [that this function is well-defined on $\mathcal
{T}'$ is readily checked from (iii) and (\ref{bb})]. Moreover, it is
not difficult to check that
\[
\mathop{\sup_{x,y\in\mathcal{T'}^{(r)}:}}_{d_\mathcal{T}(x,y)\leq\delta
(r,\varepsilon)}d_M \bigl(
\phi_{\mathcal{T}}(x),\phi_\mathcal {T}(y) \bigr)<\varepsilon,
\]
and so the function can be extended continuously to the whole of
$\mathcal{T}$. In particular, we have so far constructed $\uT$, and to
check this is an element of $\mathbb{T}$, it remains to show that
$(\mathcal{T},d_\mathcal{T})$ is a real tree. However, in \cite{ADH}, Lemma~2.7, it is shown that $(\mathcal{T},d_\mathcal{T})$ is a length
space, and so it is connected. Moreover, the four-point condition for
the metric for $(\mathcal{T},d_\mathcal{T})$ follows from the
four-point condition that must hold for $(\mathcal{T}_n,d_{\mathcal
{T}_n})$ (see \cite{EPW}, (2.1)). It follows that $(\mathcal
{T},d_\mathcal{T})$ must be a real tree, as desired.

It remains to show that $\uT_{n_i}\rightarrow\uT$
in $(\mathbb{T},\Delta)$. For this it is sufficient to show that
$\uT_{n_i}^{(r)}\rightarrow\uT^{(r)}$ in $(\mathbb{T}_c,\Delta_c)$, at
least whenever $\mu_\mathcal{T}(\partial B_\mathcal{T}(\rho_\mathcal
{T},r))=0$. Again, this may be accomplished by following the argument
of \cite{ADH}, which involves introducing finite subsets
$U_{k,l}\subset U$ such that $\{x_u^{n_i}:u\in U_{k,l}\}$ and $\{
x_u:u\in U_{k,l}\}$ suitably well-approximate $\mathcal{T}^{(r)}_{n_i}$
and $\mathcal{T}^{(r)}$, respectively. Moreover, a consideration of the
correspondence between these finite sets given by $(x_u^{n_i},x_u)$,
$u\in U_{k,l}$, allows it to be deduced in our case that
\begin{eqnarray*}
{\lim_{i\rightarrow\infty}\Delta_c \bigl(\uT_{n_i}^{(r)},
\uT^{(r)} \bigr)}&\leq& 2\sup_{i\geq1} \mathop{\sup
_{x,y\in\mathcal{T}^{(r+\delta)}_{n_i}:}}_{d_{\mathcal
{T}_{n_i}}(x,y)\leq\delta}d_M \bigl(
\phi_{\mathcal{T}_{n_i}}(x),\phi _{\mathcal{T}_{n_i}}(y) \bigr)\\
&&{}+2\mathop{\sup
_{x,y\in\mathcal
{T}^{(r+\delta)}:}}_{d_\mathcal{T}(x,y)\leq\delta}d_M \bigl(\phi
_{\mathcal{T}}(x),\phi_\mathcal{T}(y) \bigr),
\end{eqnarray*}
for any $\delta>0$ [cf. the extra term involving the continuity of $\phi
_\mathcal{T}$ in (\ref{extra})]. Since the right-hand side can be made
arbitrarily small by suitable choice of $\delta$, this completes the proof.
\end{pf}

\begin{rem}The restriction to real trees for $(\mathbb{T}_c,\Delta_c)$ has
actually been unnecessary in this section so far, and so the same
topology could be extended to the setting where the metric space part
of an element---$(\mathcal{T},d_\mathcal{T})$---is simply assumed to
be a compact metric space. Similarly, for the topology, $(\mathbb
{T},\Delta)$, it would have been enough to assume that the metric space
part of an element is a locally compact length space (cf. \cite{ADH}).
In both cases, the restriction to the case where the metric space is a
real tree would then simply be the restriction to a closed subset of
the relevant topology (cf. \cite{Evans}, Lemma~4.22).
\end{rem}

To conclude this section, we present two consequences of convergence in
$(\mathbb{T}_c,\Delta_c)$, again assuming that $(M,d_M)$ is proper.
First, we prove convergence of the push-forward measures. In what
follows, $B_X(x,r)$ is the open ball in the metric space $X=(X,d_X)$
with radius $r$ centred at $x$.

\begin{lem}\label{measconv} Suppose $(M,d_M)$ is proper. If $\uT
_n\rightarrow\uT$ in $(\mathbb{T}_c,\Delta_c)$, then
%
\begin{equation}
\label{mconv} \mu_{\mathcal{T}_n}\circ\phi_{\mathcal{T}_n}^{-1}
\rightarrow\mu _{\mathcal{T}}\circ\phi_{\mathcal{T}}^{-1}
\end{equation}
weakly as Borel measures on $(M,d_M)$.
\end{lem}

\begin{pf}
Note first that if
$\uT_n\rightarrow\uT$ in $(\mathbb{T}_c,\Delta_c)$ then for each $n$
we can find
a measurable function $f_n:\mathcal{T}_n\rightarrow\mathcal{T}$ such that
$\mu_{\mathcal{T}_n}\circ f_n^{-1}\rightarrow\mu_{\mathcal{T}}$ weakly
as measures on $\mathcal{T}$, and also
%
\begin{equation}
\label{modcont} \sup_{x\in\mathcal{T}_n}d_M \bigl(
\phi_\mathcal{T}\bigl(f_n(x)\bigr), \phi_{\mathcal{T}_n}(x)
\bigr)\rightarrow0.
\end{equation}
Indeed, let $Z_n,\psi_n,\psi_n',\mathcal{C}_n$ be defined as in the
proof of Proposition~\ref{deltametric}, that is, so that (\ref{yy})
holds. Let $(x^n_i)_{i=1}^{N(n)}$ be a $\varepsilon_n$-cover
of $\mathcal{T}$. Set $A^n_1:=B_{Z_n}(\psi_n(x_1^n),2\varepsilon_n)$
and $A^n_i:=B_{Z_n}(\psi_n(x_i^n),2\varepsilon_n)\setminus A^n_{i-1}$
for $i=2,\ldots, N(n)$.
Then the sets $A^n_i$, $i=1,\ldots,N(n)$, are disjoint and their union
contains all those points in $Z_n$ within a distance $\varepsilon_n$ of
$\psi_n(\mathcal{T})$. In particular, they cover $\psi_n'(\mathcal
{T}_n)$, so one can define a (measurable) map $f_n:\mathcal
{T}_n\rightarrow\mathcal{T}$ by setting $f_n(x):=x_i^n$ if $\psi
_n'(x)\in A_i^n$. For this map, we have
\begin{eqnarray}
\label{e:distaa} d_P^\mathcal{T}\bigl(\mu_{\mathcal{T}_n}\circ
f_n^{-1},\mu_\mathcal{T}\bigr)&\leq&
d_P^{Z_n}\bigl(\mu_{\mathcal{T}_n}\circ
f_n^{-1}\circ\psi_n^{-1},
\mu_{\mathcal{T}_n}\circ\psi_n'^{-1}\bigr)
\nonumber
\\[-8pt]
\\[-8pt]
\nonumber
&&{}+
d_P^{Z_n}\bigl(\mu_{\mathcal{T}_n}\circ
\psi_n'^{-1},\mu_{\mathcal{T}}\circ
\psi_n^{-1}\bigr),
\end{eqnarray}
where $d_P^\mathcal{T}$ is the Prohorov distance between measures on
$\mathcal{T}$.
By (\ref{yy}), the second term in \eqref{e:distaa} is bounded above by
$\varepsilon_n$.
Moreover, by definition we have that $d_{Z_n}(\psi_n(f_n(x)),\psi
_n'(x))$ is strictly less than $2\varepsilon_n$ for all $x\in\mathcal
{T}_n$, and so the first term is bounded above by $2\varepsilon_n$.
This confirms that $\mu_{\mathcal{T}_n}\circ f_n^{-1}\rightarrow\mu
_{\mathcal{T}}$. Next, observe that if $f_n(x)=x_i^n$ and $(x',x)\in
\mathcal{C}_n$, then
\begin{eqnarray*}
d_M \bigl(\phi_\mathcal{T}\bigl(f_n(x)\bigr),
\phi_{\mathcal{T}_n}(x) \bigr)&\leq& \varepsilon_n+d_M
\bigl(\phi_\mathcal{T}\bigl(x_i^n\bigr),
\phi_{\mathcal
{T}}\bigl(x'\bigr) \bigr)\\
&\leq&\varepsilon_n+
\mathop{\sup_{(x,x')\in\mathcal
{T}:}}_{d_\mathcal{T}(x,x')<3\varepsilon_n}d_M \bigl(
\phi_{\mathcal{T}}(x), \phi_{\mathcal{T}}\bigl(x'\bigr) \bigr).
\end{eqnarray*}
By the continuity of $\phi_{\mathcal{T}}$, this upper bound converges
to zero as $n\rightarrow\infty$, and we have thereby established (\ref
{modcont}). As a consequence, if $g:M\rightarrow\mathbb{R}$ is
continuous and compactly supported, then
\begin{eqnarray*}
&&{\bigl\llvert \mu_{\mathcal{T}_n}\circ\phi_{\mathcal{T}_n}^{-1}(g) -
\mu_{\mathcal{T}}\circ\phi_{\mathcal{T}}^{-1}(g)\bigr\rrvert } \\
&&\qquad\leq
\int_{\mathcal{T}_n}\bigl\llvert g\bigl(\phi_{\mathcal{T}_n}(x)\bigr)-g
\bigl(\phi _\mathcal{T}\bigl(f_n(x)\bigr)\bigr)\bigr\rrvert
\mu_{\mathcal{T}_n}(dx)
\\
&&\qquad\quad{}+\biggl\llvert \int_{\mathcal{T}}g\bigl(\phi_\mathcal{T}(x)
\bigr)\mu_{\mathcal
{T}_n}\circ f_n^{-1}(dx)-\int
_\mathcal{T}g\bigl(\phi_{\mathcal{T}}(x)\bigr)\mu
_{\mathcal{T}}(dx)\biggr\rrvert
\\
&&\qquad\rightarrow 0,
\end{eqnarray*}
where the convergence of the first term in the upper bound to zero
follows from (\ref{modcont}) [and the fact that $\mu_{\mathcal
{T}_n}(\mathcal{T}_n)\rightarrow\mu_\mathcal{T}(\mathcal{T})<\infty$,
as follows from $\mu_{\mathcal{T}_n}\circ f_n^{-1}\rightarrow\mu
_{\mathcal{T}}$], and the convergence of the second term to zero also
follows from $\mu_{\mathcal{T}_n}\circ f_n^{-1}\rightarrow\mu_{\mathcal
{T}}$. This establishes that $\mu_{\mathcal{T}_n}\circ\phi_{\mathcal
{T}_n}^{-1}$ converges vaguely to $\mu_{\mathcal{T}}\circ\phi_{\mathcal
{T}}^{-1}$. Finally, since the masses of the measures in the sequence
converge to the mass of the limit, which is finite, it also
demonstrates weak convergence.
\end{pf}

\begin{rem}\label{measrem}
It is not difficult to extend the above proof to deduce that the
conclusion of (\ref{mconv})
holds in the sense of vague convergence of measures whenever $\uT
_n\rightarrow\uT$ in
$(\mathbb{T},\Delta)$, and in addition we have the following
condition which prevents an explosion of mass in a bounded region of
the proper space $(M,d_M)$: for each $r\in(0,\infty)$, there exists an
$R<\infty$ such that
%
\begin{equation}
\label{rR} \phi_{\mathcal{T}_n}^{-1} \bigl(B_M(
\rho_M,r) \bigr)\subseteq B_{\mathcal
{T}_n}(\rho_{\mathcal{T}_n},R)\qquad
\mbox{for all } n,
\end{equation}
where $\rho_M$ is a distinguished point in $M$. We will apply a
probabilistic version of such an argument to prove Lemma~\ref{measlem}.
\end{rem}

Our second result is that
convergence in $\mathbb{T}_c$ with respect to $\Delta_c$ implies
convergence in
a generalisation of the topology for path ensembles considered by
Schramm in \cite{Schramm}.
In that paper, the space $M$ considered was the one-point
compactification of $\mathbb{R}^2$,
$\mathbb{S}^2$ say.
This result will be used when we wish to transfer the results of \cite
{Schramm} to our setting.
Given a metric space $X$, write $\mathcal{H}(X)$ for the Hausdorff
space of compact
subsets of $X$. We write $\gamma_{\mathcal{T}}(x,y)$ for the unique path
between $x$ and $y$ in $\mathcal{T}$ (including its endpoints).

\begin{lem}\label{tlem} If we define
\[
\mathfrak{T} (\uT ):= \bigl\{ \bigl(\phi_\mathcal{T}(x),\phi
_\mathcal{T}(y),\phi_\mathcal{T}\bigl(\gamma_\mathcal{T}(x,y)
\bigr) \bigr): x,y\in\mathcal{T} \bigr\},
\]
then the convergence $\uT_n\rightarrow\uT$ in $(\mathbb{T}_c,\Delta
_c)$ implies that $\mathfrak{T}(\uT_n)\rightarrow\mathfrak{T}(\uT)$ in
$\mathcal{H}(M\times M\times\mathcal{H}(M))$.
\end{lem}

\begin{pf} Suppose that $\uT_n\rightarrow\uT$ holds in $(\mathbb
{T}_c,\Delta_c)$, and that $Z_n,\psi_n,\psi_n',\mathcal{C}_n$ are
defined as in the proof of Proposition~\ref{deltametric}, so that (\ref
{yy}) holds. We claim that if $(x,x_n),(y,y_n)\in\mathcal{C}_n$, then
%
\begin{eqnarray}
\label{dh}&& d_H^{M} \bigl(\phi_{\mathcal{T}}\bigl(
\gamma_\mathcal{T}(x,y)\bigr),\phi_{\mathcal
{T}_n}\bigl(
\gamma_{\mathcal{T}_n}(x_n,y_n)\bigr) \bigr)
\nonumber
\\[-8pt]
\\[-8pt]
\nonumber
&&\qquad\leq
\eta_n:=\varepsilon_n+ \mathop{\sup_{(z,z')\in\mathcal{T}:}}_{d_\mathcal{T}(z,z')<5\varepsilon
_n}d_M
\bigl(\phi_{\mathcal{T}}(z), \phi_{\mathcal{T}}\bigl(z'\bigr)
\bigr),
\end{eqnarray}
where $d_H^{M}$ is the Hausdorff distance between subsets of $M$. To
prove this, we start by considering $z\in\gamma_\mathcal{T}(x,y)$, and
defining $z_n$ to be any element of $\mathcal{T}_n$ such that
$(z,z_n)\in\mathcal{C}_n$. By applying (\ref{yy}) and the fact that the
metric $d_\mathcal{T}$ is additive along paths, we obtain
\begin{eqnarray*}
d_{\mathcal{T}_n}(x_n,z_n)+d_{\mathcal{T}_n}(z_n,y_n)&<&d_{\mathcal
{T}}(x,z)+d_{\mathcal{T}}(z,y)
+ 4\varepsilon_n=d_\mathcal{T}(x,y)+4\varepsilon_n
\\
&<&
d_{\mathcal
{T}_n}(x_n,y_n)+6\varepsilon_n.
\end{eqnarray*}
It follows that $z_n$ is within a distance of $3\varepsilon_n$ (with
respect to $d_{\mathcal{T}_n}$) of $\gamma_{\mathcal{T}_n}(x_n,y_n)$.
Now, if we let $z_n'$ be the closest point in $\gamma_{\mathcal
{T}_n}(x_n,y_n)$ to $z_n$, and $z_n''$ be such that $(z_n'',z_n')\in
\mathcal{C}_n$, then it is the case that $d_\mathcal{T}(z,z_n'')<
d_{\mathcal{T}_n}(z_n,z_n')+2\varepsilon_n<5\varepsilon_n$, and so
$d_M(\phi_{\mathcal{T}}(z),\phi_{\mathcal{T}_n}(z_n')) <\varepsilon
_n+d_M(\phi_{\mathcal{T}}(z),\phi_{\mathcal{T}}(z_n''))\leq\eta_n$.
Thus $\phi_\mathcal{T}(z)$ is within a $d_M$-distance $\eta_n$ of $\phi
_{\mathcal{T}_n}(\gamma_{\mathcal{T}_n}(x_n,y_n))$. A similar argument
yields that any point of $\phi_{\mathcal{T}_n}(\gamma_{\mathcal
{T}_n}(x_n,y_n))$ is within a $d_M$-distance $\eta_n$ of $\phi_{\mathcal
{T}}(\gamma_\mathcal{T}(x,y))$. This establishes~(\ref{dh}), from which
the result follows.
\end{pf}

\begin{rem}\label{otherrem}
As with Lemma~\ref{measconv}, this result is readily extended to
the noncompact case when $(M,d_M)$ is proper. Indeed, under this
assumption, if $\uT_n\rightarrow\uT$ in $(\mathbb{T},\Delta)$ and (\ref
{rR}) holds, then $\mathfrak{T}(\uT_n)\rightarrow\mathfrak{T}(\uT)$ in
$\mathcal{H}(\dot{M}\times\dot{M}\times\mathcal{H}(\dot{M}))$,
where $\dot{M}$ is defined to be the one-point compactification of $M$.
A probabilistic version of
this argument will be used to prove Lemma~\ref{push}.
\end{rem}

\begin{rem}While we do not need the result, we note that a similar
argument can be used to
relate convergence in our topology to convergence in the topology of
\cite{ABNW}.
This topology is similar to that of Schramm, but it incorporates
convergence of the shape of subtrees
spanning an arbitrary finite number of vertices, rather than just two.
\end{rem}

\section{Tightness of UST law under rescaling}\label{tightsec}

The aim of this section is to prove Theorem~\ref{main1}, that is, to
establish that the
law of the UST, considered as a measured, rooted, spatial tree, is
tight under rescaling.
The key estimates for this purpose were already established in Section~\ref{estsec}.
As discussed in the \hyperref[sec1]{Introduction}, here we extend $\mathcal{U}$ to a
(locally compact)
real tree by adding line segments of unit length along its edges, and
define $\phi_\mathcal{U}:\mathcal{U}\rightarrow\mathbb{R}^2$ to be the
identity map on
vertices with linear interpolation along edges.
Throughout this section, we suppose that the image space $(M,d_M)$
introduced in
Section~\ref{topsec} is $\mathbb{R}^2$ equipped with the Euclidean distance.

\begin{lem}\label{L:51}
For every $r>1$, $\varepsilon\in(0,\varepsilon_0)$, it holds that
%
\begin{eqnarray}
\label{uni} &&\lim_{N\rightarrow\infty}\liminf_{\delta\rightarrow0}
\mathbf{P} \bigl(\mbox{there exists a $\delta^{-\kappa}\varepsilon$-cover for
$B_\mathcal{U}\bigl(0,\delta^{-\kappa}r\bigr)$}\nonumber\\
&&\hspace*{58pt}\quad\mbox{of cardinality $\leq
N$} \bigr)\\
&&\qquad=1.\nonumber
\end{eqnarray}
\end{lem}

\begin{pf} Recalling the notation $N_\sU$ introduced above Lemma~\ref{disclem}.
we have that the probability in (\ref{uni}) is at least $\mathbf
{P}(N_\sU(\delta^{-\kappa}r,\delta^{-\kappa}\varepsilon)\leq N)$.
Let $\th=\th(N)$ be such that $c \th^{1+d_f} (r/\eps)^{d_f} = N$; then by
\eqref{e:ballbound} we have that
$\limsup_{\delta\rightarrow0}\mathbf{P}(N_\sU(\delta^{-\kappa}r,\delta
^{-\kappa}\varepsilon)\ge N) \le c \exp( -\th^{1/80})$,
and since $\lim_{N \to\infty} \th(N) = \infty$, this proves the result.
\end{pf}

\begin{lem}
For every $r<\infty$, it holds that
\[
\lim_{\lambda\rightarrow\infty}\liminf_{\delta\rightarrow0} \mathbf{P} \bigl(
\delta^2\mu_\mathcal{U} \bigl(B_\mathcal{U}\bigl(0,
\delta ^{-\kappa}r\bigr) \bigr)\leq\lambda \bigr)=1.
\]
\end{lem}

\begin{pf} This is a simple consequence of
Theorem~\ref{vollemust}(b).
\end{pf}

\begin{lem}\label{thmm:lem4-3}
For every $\varepsilon>0$, $r<\infty$, it holds that
\[
\lim_{\eta\rightarrow0}\liminf_{\delta\rightarrow0}\mathbf{P} \Bigl(
\mathop{\max_{x,y\in B_\mathcal{U}(0,\delta^{-\kappa}r):}}_{d_\mathcal
{U}(x,y)\leq\delta^{-\kappa}\eta}\bigl\llvert
\phi_\mathcal{U}(x)-\phi _{\mathcal{U}}(y)\bigr\rrvert \le
\delta^{-1}\varepsilon \Bigr)=1.
\]
\end{lem}

\begin{pf} Since $|\phi_\mathcal{U}(x)-\phi_\mathcal{U}(y)|\leq\dS(x,y)$
it is sufficient to prove that
%
\begin{equation}
\label{eq:fnofjf} \lim_{\eta\rightarrow0}\liminf_{\delta\rightarrow0}
\mathbf{P} \Bigl( \mathop{\max_{x,y\in B_\mathcal{U}(0,c_1\delta^{-\kappa}r):}}_{
d_\mathcal{U}(x,y)\leq c_2\delta^{-\kappa}\eta}\dS(x,y)>
\delta ^{-1}\varepsilon \Bigr)=0.
\end{equation}
Let $r'=\delta^{-1}\varepsilon$, and set $A_*(\lambda):=\{B_\sU
(0,c_1\delta^{-\kappa}r)\subset B_E(0,r'e^{c_1\lambda^{1/2}})\}$, where
$c_1$ is the constant of Proposition~\ref{thmm:keycomp}. By Theorem~\ref
{vollemust}(a),
\[
\mathbf{P} \bigl( \bigl\{B_\sU\bigl(0,c_1
\delta^{-\kappa}r\bigr)\subset B_E\bigl(0,(\lambda
r)^{1/\kappa}\delta^{-1}\bigr) \bigr\}^c \bigr) \le
c_2e^{-c_3\lambda^{2/3}}\qquad \forall\delta^\kappa\leq\lambda r,
\lambda\ge c_1,
\]
and in addition $(\lambda r)^{1/\kappa}\delta^{-1}\le r'e^{c_1\lambda
^{1/2}}$ for $\lambda$ large.
Thus, $\mathbf{P}(A_*(\lambda)^c)\le\break  c_2e^{-c_3\lambda^{2/3}}$ for all
$\delta^\kappa\leq\lambda r$, $\lambda\ge\lambda_1$, where $\lambda
_1$ is some large, finite constant.
Next, let $A_4$ be as in Proposition~\ref{thmm:keycomp} [taking
$(x_0,r,\lambda)$ in that result to be $(0,r',\lambda)$ in our current
parameterisation], so that $\bP(A_4^c) \le c_4 \exp( -c_5 \lam^{1/2})$
for all $\delta\leq\varepsilon$, $\lambda\geq\lambda_0$. Clearly, it
is enough to consider the event of \eqref{eq:fnofjf} on $A_*(\lambda
)\cap A_4$.
On $A_*(\lambda)\cap A_4$, if $x,y\in B_\mathcal{U}(0,c_1\delta^{-\kappa
}r)$ satisfy
$\dS(x,y)> \delta^{-1}\varepsilon=r'$, then by Proposition~\ref{thmm:keycomp},
$d_\mathcal{U}(x,y)\ge\lambda^{-1}\dS(x,y)^{\kappa}\ge\lambda
^{-1}\varepsilon^{\kappa}\delta^{-\kappa}$. Thus, by taking $\eta<
\lambda^{-1}\varepsilon^{\kappa}$, $d_\mathcal{U}(x,y)> \eta\delta
^{-\kappa}$, so \eqref{eq:fnofjf} is proved.
\end{pf}

\begin{pf*}{Proof of Theorem~\ref{main1}}
This is clear given the pre-compactness result of Lemma~\ref{relcompact},
and Lemmas \ref{L:51}--\ref{thmm:lem4-3}.
\end{pf*}

\section{Properties of limit measures}\label{propsec}

In this section, we establish properties of the limit measure of the
UST and will
prove Theorem~\ref{main12}. Throughout, we fix a sequence $\delta
_n\rightarrow0$ such that $(\mathbf{P}_{\delta_n})_{n\geq1}$
converges weakly [as measures on $(\mathbb{T},\Delta)$], and write $\uU
_{\delta_n}=(\mathcal{U},\delta_n^\kappa d_\mathcal{U},\delta_n^2\mu
_\mathcal{U},\delta_n\phi_\mathcal{U},0)$. Letting $\tilde{\mathbf{P}}$
be the relevant limiting law, we denote by $\uT=(\mathcal{T},d_\mathcal
{T},\mu_\mathcal{T},\phi_\mathcal{T},\rho_\mathcal{T})$ a random
variable with law $\tilde{\mathbf{P}}$. Again, we take the image space
$(M,d_M)$ of Section~\ref{topsec} to be $\mathbb{R}^2$ equipped with
the Euclidean distance. In many of the arguments, the following
coupling result will be useful.

\begin{lem}\label{couplem} There exist realisations of $(\uU_{\delta
_n})_{n\geq1}$ and $\uT$ built on the same probability space, with
probability measure $\mathbf{P}^*$ say, such that: for some subsequence
$(n_i)_{i\geq1}$ and divergent sequence $(r_j)_{j\geq1}$ it holds
that, $\mathbf{P}^*$-a.s.,
%
\begin{equation}
\label{dijdef} D_{i,j}:=\Delta_c \bigl(
\uU_{\delta_{n_i}}^{(r_j)},\uT^{(r_j)} \bigr)\rightarrow0
\end{equation}
as $i\rightarrow\infty$, for every $j\geq1$.
\end{lem}

\begin{pf} Recall that by the definition of $\tilde{\mathbf{P}}$ we have
that $\uU_{\delta_n}\rightarrow\uT$ in distribution (where the laws of
random variables on the left-hand side are considered under $\mathbf
{P}$, and those on the right under $\tilde{\mathbf{P}}$). Thus, since
the space $(\mathbb{T},\Delta)$ is separable (see Proposition~\ref
{deltametric}), we can suppose that we have versions of the random
variables built on a common probability space, with probability measure
$\mathbf{P}^*$, such that the convergence holds $\mathbf{P}^*$-a.s.
From the definition of $\Delta$ and Fubini's theorem, it follows that
$\int_{0}^\infty e^{-r}
(1\wedge\mathbf{E}^*(\Delta_c(\uU_{\delta_n}^{(r)},\break \uT
^{(r)})))\,dr\rightarrow0$.
Some standard analysis now yields that there exists a subsequence
$(n_i)_{i\geq1}$
such that for Lebesgue almost-every $r$, $\mathbf{E}^*(\Delta_c(\uU
_{\delta_{n_i}}^{(r)},\uT^{(r)}))\rightarrow0$.
In turn, letting $(r_j)_{j\geq1}$ be a divergent sequence such that
the above holds for every $r_j$,
a straightforward diagonalisation argument yields the result.
\end{pf}

We now show that the push-forward of $\mu_\mathcal{T}$ by $\phi_\mathcal
{T}$ is $\tilde{\mathbf{P}}$-a.s. equal to Lebesgue measure on $\mathbb{R}^2$.

\begin{lem}\label{measlem} $\tilde{\mathbf{P}}$-a.s., it holds that $\mu
_\mathcal{T}\circ\phi_\mathcal{T}^{-1}=\mathcal{L}$.
\end{lem}

\begin{pf} We first note that since $\delta^{2}\mu_\mathcal{U}\circ\phi
_\mathcal{U}^{-1}(\delta^{-1} \cdot)\rightarrow\mathcal{L}$ for any
realisation of the UST, it will suffice to show that
%
\begin{equation}
\label{a9} \delta_n^{2}\mu_\mathcal{U}\circ
\phi_\mathcal{U}^{-1}\bigl(\delta_n^{-1}
\cdot\bigr)\rightarrow\mu_\mathcal{T}\circ\phi_\mathcal{T}^{-1}
\end{equation}
in distribution with respect to the topology of vague convergence of
probability measures on $\mathbb{R}^2$. For this, it will be enough to
establish that, for any continuous, positive, compactly supported
function $f$,
%
\begin{equation}
\label{a10} \delta_n^{2}\int_{\mathbb{R}^2}f(
\delta_n x)\mu_\mathcal{U}\circ\phi _\mathcal{U}^{-1}(dx)
\rightarrow\int_{\mathbb{R}^2} f(x)\mu_\mathcal{T}\circ
\phi_\mathcal{T}^{-1}(dx)
\end{equation}
in distribution (see \cite{Kall}, Theorem~16.16, e.g.).

Applying the coupling of Lemma~\ref{couplem} in conjunction with Lemma~\ref{measconv}, we obtain
%
\begin{equation}
\label{a11} \delta_{n_i}^{2}\mu_\mathcal{U} \bigl(
\phi_\mathcal{U}^{-1} \bigl(\delta _{n_i}^{-1}
\cdot \bigr)\cap {B}_\mathcal{U} \bigl(0,\delta_{n_i}^{-\kappa}r_j
\bigr) \bigr)\rightarrow\mu_\mathcal{T} \bigl(\phi_\mathcal{T}^{-1}(
\cdot)\cap \mathcal{T}^{(r_j)} \bigr)
\end{equation}
weakly as measures on $\mathbb{R}^2$ as $i\rightarrow\infty$, for every
$r_j$, $\mathbf{P}^*$-a.s. In particular, this confirms that, for every
$r_j$, the above convergence holds in distribution (under the
convention that the laws of random variables on the left-hand side are
considered under $\mathbf{P}$, and those on the right under $\tilde
{\mathbf{P}}$). By monotonicity, we also clearly have $\tilde{\mathbf
{P}}$-a.s. that, for any positive measurable $f$,
%
\begin{equation}
\label{a12} \int_{\mathbb{R}^2}f(x)\mu_\mathcal{T} \bigl(
\phi_\mathcal{T}^{-1}(\cdot )\cap\mathcal{T}^{(r)}
\bigr) (dx)\rightarrow\int_{\mathbb{R}^2}f(x)\mu _\mathcal{T}
\circ\phi_\mathcal{T}^{-1}(dx),
\end{equation}
as $r\rightarrow\infty$.

As a consequence of (\ref{a11}) and (\ref{a12}), to establish the
convergence at (\ref{a10}) along the subsequence $(n_i)_{i\geq1}$, it
is sufficient to show that $\mu_\mathcal{T}\circ\phi_\mathcal{T}^{-1}$
is locally finite and also that, for any continuous, positive,
compactly supported function $f$,
%
\begin{equation}
\label{a13} \lim_{j\rightarrow\infty}\limsup_{i\rightarrow\infty}
\delta _{n_i}^{2}\biggl\llvert \mathbf{E} \biggl(\int
_{{{{B}_\mathcal{U} (0,\delta
_{n_i}^{-\kappa}r_j )^c}}}f\bigl(\delta_{n_i}\phi_\mathcal{U}( x)
\bigr)\mu _\mathcal{U}(dx) \biggr)\biggr\rrvert =0
\end{equation}
(cf. \cite{Bill}, Theorem~3.2). To show that the latter is true, first
choose $r$ such that the support of $f$ is contained within $\overline
{B}_E(0,r)$ (where we write $\overline{A}$ to represent the closure of
a set $A$), and define $A(i,j)$ to be the event that
%
\begin{equation}
\label{a14} \phi_{\mathcal{U}}^{-1} \bigl(\overline{B}_E
\bigl(0,\delta_{n_i}^{-1}r\bigr) \bigr)\subseteq{B}_\mathcal{U}
\bigl(0,\delta_{n_i}^{-\kappa}r_j \bigr)
\end{equation}
[i.e., similar to the inclusion at (\ref{rR})]. It is then the case
that the expression within the limits on the left-hand side of (\ref
{a13}) is equal to
\[
\delta_{n_i}^{2}\biggl\llvert \mathbf{E} \biggl(\int
_{{B}_\mathcal{U}
(0,\delta_{n_i}^{-\kappa}r_j )^c}f\bigl(\delta_{n_i}\phi_\mathcal{U}( x)
\bigr)\mu_\mathcal{U}(dx)\mathbf{1}_{A(i,j)^c} \biggr)\biggr\rrvert ,
\]
which is bounded above by
$\sup_{x\in\mathbb{R}^2}f(x)\delta_{n_i}^{2}\mu_\mathcal{U}(\overline
{B}_E(0,\delta_{n_i}^{-1}r))\mathbf{P}(A(i,j)^c)\leq  c\mathbf
{P}(A(i,j)^c)$ for some finite constant $c$. Consequently, since
%
\begin{equation}
\label{aij} \lim_{j\rightarrow\infty}\limsup_{i\rightarrow\infty}
\mathbf{P}\bigl(A(i,j)^c\bigr)=0
\end{equation}
by Theorem~\ref{vollemust}(a), we have proved (\ref{a13}), as desired.

To check that $\mu_\mathcal{T}\circ\phi_\mathcal{T}^{-1}$ is locally
finite, we will show that, for every $r>0$,
%
\begin{equation}
\label{toprove} \lim_{R\rightarrow\infty}\tilde{\mathbf{P}} \bigl(
\phi_{\mathcal
{T}}^{-1} \bigl(\overline{B}_E(0,r) \bigr)
\nsubseteq{\mathcal {T}}^{(R)} \bigr)=0.
\end{equation}
Suppose that (\ref{toprove}) is not true for some $r'>0$, with the
limit instead being equal to $\varepsilon>0$. It is then the case that
for every $R$, there exists an $R'$ such that
%
\begin{equation}
\label{event} \tilde{\mathbf{P}} \bigl(\phi_{\mathcal{T}}(x)\in\overline
{B}_E\bigl(0,r'\bigr)\mbox{ for some }x\in{
\mathcal{T}}^{(R')}\setminus{\mathcal {T}}^{(R)} \bigr)\geq
\varepsilon/2.
\end{equation}
Next, let us suppose that the sequences $(n_i)_{i\geq1}$ and
$(r_j)_{j\geq1}$ are given by Lemma~\ref{couplem}, and $D_{i,j}$, as defined by (\ref{dijdef}), is bounded
strictly above by $\delta$. We can then find a correspondence $\mathcal
{C}_{i,j}\subseteq{B}_\mathcal{U}(0,\delta_{n_i}^{-\kappa}r_j)\times
\mathcal{T}^{(r_j)}$ such that $|\delta_{n_i}^\kappa d_\mathcal
{U}(0,x)-d_\mathcal{T}(\rho_\mathcal{T},x')|<2\delta$ and also $|\delta
_{n_i}\phi_\mathcal{U}(x)-\phi_{\mathcal{T}}(x')|\leq\delta$ for every
$(x,x')\in\mathcal{C}_{i,j}$. It is then easy to check that if $r_j>R'$
and the event within the probability at (\ref{event}) holds, then so
does the event that
$\phi_{\mathcal{U}}(x)\in\overline{B}_E(0,\delta_{n_i}^{-1}(r'+\delta
))$ for some $x\in{B}_\mathcal{U}(0,\delta_{n_i}^{-\kappa
}r_j)\setminus{B}_\mathcal{U}(0,\delta_{n_i}^{-\kappa}(R-2\delta))$,
which is a subset of $\{\phi_{\mathcal{U}}^{-1}(\overline{B}_E(0,\delta
_{n_i}^{-1}(r'+\delta)))\nsubseteq{B}_\mathcal{U}(0,\delta
_{n_i}^{-\kappa}(R-2\delta))\}$. Since we know that $\mathbf
{P}^*(D_{i,j}>\delta)\rightarrow0$ as $i\rightarrow\infty$, it follows that
\[
\liminf_{i\rightarrow\infty}\mathbf{P} \bigl(\phi_{\mathcal
{U}}^{-1}
\bigl(\overline{B}_E\bigl(0,\delta_{n_i}^{-1}
\bigl(r'+\delta\bigr)\bigr)\bigr)\nsubseteq {B}_\mathcal{U}
\bigl(0,\delta_{n_i}^{-\kappa}(R-2\delta)\bigr) \bigr)\geq
\varepsilon/2.
\]
However, replacing $r'+\delta$ by $r$ and $R-2\delta$ by $r_j$ for
suitably large $j$, we see that this contradicts the statement at (\ref
{aij}). Consequently, (\ref{toprove}) must actually be true.

What we have proved already is enough to yield the lemma. We do note,
though, that for any subsequence $(n_i)_{i\geq1}$, we could have
applied the same argument to find a sub-subsequence $(n_{i_j})_{j\geq
1}$ along which the convergence at (\ref{a10}) holds. Since the limit
is identical for any such sub-subsequence, it must be the case that the
full sequence also converges to this limit, which thereby establishes
(\ref{a9}).
\end{pf}

Next, similar to (\ref{sm}), define a ``Schramm-metric'' on $\mathcal
{T}$ by setting, for \mbox{$x,y\in\mathcal{T}$},
\[
d_\mathcal{T}^S(x,y):=\operatorname{ diam} \bigl(\phi_\mathcal{T}
\bigl(\gamma_\mathcal {T}(x,y)\bigr) \bigr),
\]
where the diameter is in the Euclidean metric. It follows immediately
from the continuity of $\phi_\mathcal{T}$
that $d_\mathcal{T}^S$ takes values in $[0,\infty)$, and it is easy to
verify from the definition that $d^S_\sT$ is symmetric and satisfies
the triangle inequality. In the next two lemmas, we show that $d^S_\sT$
is a metric on $\sT$, and that it
gives the same topology as~$d_\sT$.

\begin{lem}\label{dsbound} For every $r,\eta>0$, we have
%
\begin{equation}
\label{ew2} \lim_{\varepsilon\rightarrow0}\tilde{\mathbf{P}} \Bigl(\mathop{\inf
_{x,y\in B_\mathcal{T}
(\rho_\mathcal{T},r):}}_{d_\mathcal{T}(x,y)\geq\eta}d_\mathcal {T}^S(x,y)<
\varepsilon \Bigr)=0.
\end{equation}
\end{lem}

\begin{pf} We start by proving the discrete analogue of the result: for
every $r,\eta>0$,
%
\begin{equation}
\label{disc} \lim_{\varepsilon\rightarrow0}\limsup_{\delta\rightarrow0}{
\mathbf{P}} \Bigl(\mathop{\inf_{x,y\in B_\mathcal{U}(0,\delta^{-\kappa}r):}}_{
d_\mathcal{U}(x,y)\geq\delta^{-\kappa}\eta}d_\mathcal{U}^S(x,y)<
\delta ^{-1}\varepsilon \Bigr)=0.
\end{equation}
We argue similar to the proof of Lemma~\ref{thmm:lem4-3}.
Again, it is enough to consider the event in $A_*(\lambda)\cap A_4$. On
$A_*(\lambda)\cap A_4$, if $x,y\in B_\mathcal{U}(0,c_1\delta^{-\kappa
}r)$ satisfy $\dS(x,y)\le \delta^{-1}\varepsilon$, then by Proposition~\ref{thmm:keycomp}, $d_\mathcal{U}(x,y)\le\lambda\varepsilon^{\kappa
}\delta^{-\kappa}$. Thus, by taking $\varepsilon$ small enough
so that $\lambda\varepsilon^{\kappa}<\eta$,
we have $d_\mathcal{U}(x,y)< \eta\delta^{-\kappa}$. This proves \eqref
{disc} with $r$ replaced by $c_1 r$, and a simple reparameterisation
yields the result.

To transfer to the continuous setting, let us suppose that the
sequences $(n_i)_{i\geq1}$
and $(r_j)_{j\geq1}$ are defined by Lemma~\ref{couplem} and that the
$\Delta_c$ distance between $\uU_{\delta_{n_i}}^{(r_j)}$ and $\uT
^{(r_j)}$, again denoted by $D_{i,j}$, is bounded strictly above by
$\delta$.
Similarly to (\ref{dh}), we can then find a correspondence $\mathcal
{C}_{i,j}\subseteq{B}_\mathcal{U}(0,\delta_{n_i}^{-\kappa}r_j)\times
\mathcal{T}^{(r_j)}$ such that
for every $(x,x'),(y,y')\in\mathcal{C}_{i,j}$, we have
\[
d_H^{\mathbb{R}^2}\bigl(\delta_{n_i}
\phi_\mathcal{U}\bigl(\gamma(x,y)\bigr), \phi_{\mathcal{T}}\bigl(
\gamma_{\mathcal{T}}\bigl(x',y'\bigr)\bigr)\bigr)\leq
\delta+ \sup_{(z,z')\in\mathcal{T}^{(r_j)}: d_\mathcal{T}(z,z')<5\delta}\bigl|\phi _{\mathcal{T}}(z)- \phi_{\mathcal{T}}
\bigl(z'\bigr)\bigr|,
\]
where $d_H^{\mathbb{R}^2}$ is the Hausdorff distance on $\mathbb{R}^2$,
and so
\[
\bigl\llvert \delta_{n_i}d_\mathcal{U}^S(x,y)-d_\mathcal{T}^S
\bigl(x',y'\bigr)\bigr\rrvert \leq2\delta+ 2\mathop{
\sup_{(z,z')\in\mathcal{T}^{(r_j)}:}}_{d_\mathcal
{T}(z,z')<5\delta}\bigl\llvert \phi_{\mathcal{T}}(z)-
\phi_{\mathcal{T}}\bigl(z'\bigr)\bigr\rrvert .
\]
We can further assume that $|\delta_{n_i}^{\kappa}d_\mathcal
{U}(x,y)-d_\mathcal{T}(x',y')|\leq2\delta$ for every $(x,x'), \break (y,y')\in
\mathcal{C}_{i,j}$. Next, fix $r,\eta>0$ and select $j$ so that
$r_j>r$. It then holds that the probability on the left-hand side of
(\ref{ew2}) is bounded above by
\begin{eqnarray*}
&&{\mathbf{P}} \Bigl(\mathop{\inf_{x,y\in{B}_\mathcal{U}(0,\delta
_{n_i}^{-\kappa}r_j):}}_{d_\mathcal{U}(x,y)\geq\delta_{n_i}^{-\kappa
}(\eta-2\delta)}d_\mathcal{U}^S(x,y)<2
\delta _{n_i}^{-1}\varepsilon \Bigr)\\
&&\qquad{}+\tilde{\mathbf{P}}
\Bigl(2\delta+ 2\mathop{\sup_{(z,z')\in\mathcal{T}^{(r_j)}:}}_{d_\mathcal
{T}(z,z')<5\delta} \bigl
\llvert \phi_{\mathcal{T}}(z)- \phi_{\mathcal{T}}\bigl(z'\bigr)
\bigr\rrvert >\varepsilon \Bigr) +\mathbf{P}^* (D_{i,j}>\delta ),
\end{eqnarray*}
where $\mathbf{P}^*$ is the coupling measure defined in the statement
of Lemma~\ref{couplem}. Now, by our choice of subsequence $(n_i)_{i\geq
1}$, the final expression converges to zero as $i\rightarrow\infty$,
for any value of $\delta>0$. Since $\phi_\mathcal{T}$ is continuous,
the second term converges to zero as $\delta\rightarrow0$, for any
value of $\varepsilon>0$. Hence, we can conclude
\[
\tilde{\mathbf{P}} \Bigl(\mathop{\inf_{x,y\in B_\mathcal{T}(\rho
_\mathcal{T},r):}}_{d_\mathcal{T}(x,y)\geq\eta}d_\mathcal
{T}^S(x,y)<\varepsilon \Bigr)\leq\limsup_{i\rightarrow\infty} {
\mathbf{P}} \Bigl(\mathop{\inf_{x,y\in{B}_\mathcal{U}(0,\delta
_{n_i}^{-\kappa}r_j):}}_{d_\mathcal{U}(x,y)\geq\delta_{n_i}^{-\kappa}\eta
/2}d_\mathcal{U}^S(x,y)<2
\delta_{n_i}^{-1}\varepsilon \Bigr).
\]
Since the upper bound converges to zero as $\varepsilon\rightarrow0$
by (\ref{disc}), this completes the proof.
\end{pf}

\begin{lem}\label{homeo} $\tilde{\mathbf{P}}$-a.s., $d_\mathcal{T}^S$ is a
metric on
$\mathcal{T}$, and the identity map from $(\mathcal{T},d_\mathcal{T})$
to $(\mathcal{T},d_\mathcal{T}^S)$ is a homeomorphism.
\end{lem}

\begin{pf} To establish that $d_\mathcal{T}^S$ is a metric, it remains to
check that it is positive definite. For this, we note $\tilde{\mathbf
{P}}(d_\mathcal{T}^S(x,y)=0\mbox{ for some }x,y\in\mathcal{T}\mbox{
with}\break d_\mathcal{T}(x,y)>0)$ is equal to
\[
\lim_{r\rightarrow\infty}\lim_{\eta\rightarrow0}\lim
_{\varepsilon
\rightarrow0}\tilde{\mathbf{P}} \bigl(d_\mathcal{T}^S(x,y)<
\varepsilon \mbox{ for some }x,y\in B_\mathcal{T}(\rho_\mathcal{T},r)
\mbox{ with }d_\mathcal{T}(x,y)\geq\eta \bigr),
\]
which in turn is equal to zero by Lemma~\ref{dsbound}. Next, we check
that the identity map from $(\mathcal{T},d_\mathcal{T})$ to $(\mathcal
{T},d_\mathcal{T}^S)$ is a homeomorphism. Clearly, it is a bijection.
Moreover, its continuity follows from the continuity of $\phi_\mathcal
{T}$. For the continuity of the inverse, we start by noting that a
simple Borel--Cantelli argument yields that, $\tilde{\mathbf{P}}$-a.s.,
for every $\eta>0$, there exists a $\varepsilon_\eta>0$ such that $\inf_{{x,y\in B_\mathcal{T}(\rho_\mathcal{T},r): d_\mathcal{T}(x,y)\geq
\eta}}d_\mathcal{T}^S(x,y)>\varepsilon_\eta$.
In particular, this implies that if $x,y\in B_\mathcal{T}(\rho_\mathcal
{T},r)$ and $d_\mathcal{T}^S(x,y)\leq\varepsilon_\eta$, then $d_\mathcal
{T}(x,y)<\eta$. Hence, the identity map from $(\mathcal{T},d_\mathcal
{T}^S)$ to $(\mathcal{T},d_\mathcal{T})$ is continuous, as desired.
\end{pf}

In order to transfer results from \cite{Schramm}, we now show that the
push-forward of $\tilde{\mathbf{P}}$ by the map $\mathfrak{T}$
introduced in Lemma~\ref{tlem} gives precisely a subsequential limiting
measure as considered in the latter paper. In particular, in \cite
{Schramm}, Schramm studied properties of the subsequential limits as
$\delta\rightarrow0$ of the laws of $\mathfrak{T}(\uU_{\delta})$,
viewed as probability measures on the space $\mathcal{H}(\mathbb
{S}^2\times\mathbb{S}^2\times\mathcal{H}(\mathbb{S}^2))$ (with
$\mathbb{S}^2$
the one-point compactification of $\mathbb{R}^2$). Whilst the space
$\mathcal{H}(\mathbb{S}^2\times\mathbb{S}^2\times\mathcal{H}(\mathbb
{S}^2))$ is compact, and so it is immediate that the laws of $(\mathfrak
{T}(\uU_{\delta}))_{\delta>0}$ are tight and admit such subsequential
limits, the next result shows that along the subsequence $(\delta
_n)_{n\geq1}$ we actually have convergence, with the limit being the
law of $\mathfrak{T}(\uT)$ under $\tilde{\mathbf{P}}$.

\begin{lem}\label{push} The laws of $(\mathfrak{T}(\uU_{\delta_n}))_{n\geq
1}$ under $\mathbf{P}$ converge to the law of $\mathfrak{T}(\uT)$ under
$\tilde{\mathbf{P}}$, weakly as probability measures on $\mathcal
{H}(\mathbb{S}^2\times\mathbb{S}^2\times\mathcal{H}(\mathbb{S}^2))$.
\end{lem}

\begin{pf} We again consider the coupling of Lemma~\ref{couplem}. Together
with Lemma~\ref{tlem}, this gives that there exists a divergent
sequence $(r_j)_{j\geq1}$ such that, for every $r_j$, $\mathbf
{P}^*$-a.s., $\mathfrak{T}(\uU_{\delta_{n_i}}^{(r_j)})\rightarrow
\mathfrak{T}(\uT^{(r_j)})$ in $\mathcal{H}(\mathbb{R}^2\times\mathbb
{R}^2\times\mathcal{H}(\mathbb{R}^2))$, and thus also in $\mathcal
{H}(\mathbb{S}^2\times\mathbb{S}^2\times\mathcal{H}(\mathbb{S}^2))$.

Let $d_{\mathbb{S}^2}$ be the usual metric on $\bS^2$.
Set $\delta_r:=\sup_{x,y\in\mathbb{S}^2\setminus B_E(0,r)}d_{\mathbb
{S}^2}(x,y)$, and note that $\delta_r\rightarrow0$ as $r\rightarrow
\infty$. Let also
\[
d_{\mathbb{S}^2\times\mathbb{S}^2\times\mathcal{H}(\mathbb
{S}^2)}\bigl((x,y,A),\bigl(x',y',A'
\bigr)\bigr):=d_{\mathbb{S}^2}\bigl(x,x'\bigr)+d_{\mathbb
{S}^2}
\bigl(y,y'\bigr)+d_H^{\mathbb{S}^2}
\bigl(A,A'\bigr),
\]
where $d_H^{\mathbb{S}^2}$ is the Hausdorff distance on $\mathcal
{H}(\mathbb{S}^2)$. Now, suppose that $i$ and $j$ are indices such that
the event $A(i,j)$ holds, where $A(i,j)$ is defined as in the proof of
Lemma~\ref{measlem} [see the definition at (\ref{a14}) in particular].
Denoting by $d_H^{\mathbb{S}^2\times\mathbb{S}^2\times\mathcal
{H}(\mathbb{S}^2)}$ the Hausdorff distance on $\mathcal{H}({\mathbb
{S}^2\times\mathbb{S}^2\times\mathcal{H}(\mathbb{S}^2)})$, we claim
that on $A(i,j)$,
\[
d_H^{\mathbb{S}^2\times\mathbb{S}^2\times\mathcal{H}(\mathbb
{S}^2)} \bigl(\mathfrak{T}\bigl(\uU_{\delta_{n_i}}^{(r_j)}
\bigr),\mathfrak{T}(\uU _{\delta_{n_i}}) \bigr)<3\delta_r,
\]
and similarly, if $\phi_{\mathcal{T}}^{-1} (\overline
{B}_E(0,r) )\subseteq\mathcal{T}^{(R)}$, then
\[
d_H^{\mathbb{S}^2\times\mathbb{S}^2\times\mathcal{H}(\mathbb
{S}^2)} \bigl(\mathfrak{T}\bigl(\uT^{(R)}
\bigr),\mathfrak{T}(\uT) \bigr)<3\delta_r.
\]
Since the two statements can be proved in the same way, let us consider
only the latter.
We need to show that if $x\in\mathcal{T}\setminus\mathcal{T}^{(R)}$ and
$y\in\mathcal{T}$ then
$(\phi_\mathcal{T}(x),\phi_\mathcal{T}(y), \phi_\mathcal{T}(\gamma
_{\mathcal{T}}(x,y)))$
is within a distance of $\delta_r$ of $\mathfrak{T}(\uT^{(R)})$ with
respect to the metric $d_{\mathbb{S}^2\times\mathbb{S}^2\times\mathcal
{H}(\mathbb{S}^2)}$. First, define $x_0$, $y_0$ to be the closest point
of $\mathcal{T}^{(R)}$ to $x$, $y$, respectively, so that the triple
$(\phi_\mathcal{T}(x_0),\phi_\mathcal{T}(y_0),\phi_\mathcal{T}(\gamma
_\mathcal{T}(x_0,y_0)))$ is an element of $\mathfrak{T}(\uT^{(R)})$. By
definition, we have that $\gamma_\mathcal{T}(x_0,x)\setminus\{x_0\}$
is a subset of $\mathcal{T}\setminus\mathcal{T}^{(R)}$, and so its
image under $\phi_\mathcal{T}$ must fall outside of $\overline
{B}_E(0,r)$. A similar observation holds in the case that $y\notin
\mathcal{T}^{(R)}$. It follows that\break 
$d_{\mathbb{S}^2\times\mathbb{S}^2\times\mathcal{H}(\mathbb
{S}^2)}((\phi_\mathcal{T}(x),\phi_\mathcal{T}(y),\phi_\mathcal{T}(\gamma
_{\mathcal{T}}(x,y))),
(\phi_\mathcal{T}(x_0),\phi_\mathcal{T}(y_0),\phi_\mathcal{T}(\gamma
_\mathcal{T}(x_0, y_0))))<3\delta_r$,
as desired.

Given (\ref{aij}), (\ref{toprove}) and the conclusions of the previous
two paragraphs, it is not difficult
to show that $\mathfrak{T}(\uU_{\delta_{n_i}})$ converges to $\mathfrak
{T}(\uT)$ in distribution.
The full convergence result can be obtained from this by applying a
subsequence argument
as in the proof of Lemma~\ref{measlem}.
\end{pf}

As a consequence of the previous lemma, we immediately inherit a number
of results from \cite{Schramm}.

\begin{lem}[(see \cite{Schramm}, Theorem~1.6, Corollary~10.4)]\label{schrammlem}
For $\tilde{\mathbf{P}}$-a.e. realisation of $\uT$,
the following properties are satisfied:
\begin{longlist}[(a)]
\item[(a)] For every $(a,b,\omega)\in\mathfrak{T}(\uT)$, if $a\neq b$, then
$\omega$ is a simple path, that is, homeomorphic to $[0,1]$. If $a=b$,
then $\omega$ is a single point or homeomorphic to a circle.

\item[(b)] Considered as a subset of $\mathbb{S}^2$,
%
\begin{equation}
\label{trunkdef} \mathrm{trunk}:=\bigcup_{(a,b,\omega)\in\mathfrak{T}(\uT)}\omega
{}\Big\backslash{}\{a,b\}
\end{equation}
is a dense topological tree.

\item[(c)] For each $x\in\mathrm{trunk}$, there are at most three connected
components of \mbox{$\mathrm{trunk}\setminus\{x\}$}.
\item[(d)] The Hausdorff dimension of $\mathrm{trunk}$ is in $(1,2)$.
\end{longlist}
\end{lem}

Note that, by construction, the set trunk defined at (\ref
{trunkdef}) is actually a subset of~$\mathbb{R}^2$, and is also a dense
topological tree when considered a subset of this space. In the
following lemma, we show further that trunk is topologically equivalent
to the set $\mathcal{T}^o$ introduced in the statement of Theorem~\ref{main12}.
For the proof of this result, we observe that $\mathcal{T}^o$ can
equivalently be defined by
%
\begin{equation}
\label{t0def} \mathcal{T}^o=\bigcup_{x,y\in\mathcal{T}}
\gamma_{\mathcal
{T}}(x,y){}\Big\backslash{}\{x,y\}.
\end{equation}
Define, for $x,y\in\mathrm{ trunk}$,
$d_\mathrm{ trunk}^S(x,y):=\operatorname{ diam}  (\gamma_\mathrm{ trunk}(x,y)
)$, where $\gamma_\mathrm{ trunk}(x,y)$ is the unique path between $x$ and
$y$ in $\mathrm{ trunk}$, and the diameter is taken with regards to the
Euclidean metric. We remark that although the metric $d_\mathrm{ trunk}^S$
behaves quite differently to the Euclidean one, the topologies these
two metrics induce on $\mathrm{ trunk}$ are the same; see the proof of
Theorem~\ref{main12} for details (cf. \cite{Schramm}, Remark~10.15).

\begin{lem}\label{trunklem} $\tilde{\mathbf{P}}$-a.s., $\phi_\mathcal{T}$ is
an isometry from $(\mathcal{T}^o,d_\mathcal{T}^S)$ to $(\mathrm{
trunk},d_\mathrm{ trunk}^S)$.
\end{lem}

\begin{pf} We start the proof by establishing that
%
\begin{equation}
\label{claim} \phi_\mathcal{T} \bigl(\gamma_{\mathcal{T}}(x,y)\setminus
\{x,y\} \bigr)=\phi_\mathcal{T} \bigl(\gamma_{\mathcal{T}}(x,y) \bigr)
\setminus\bigl\{\phi _\mathcal{T}(x),\phi_\mathcal{T}(y)\bigr\}
\end{equation}
for every $x,y\in\mathcal{T}$. The inclusion $\supseteq$ is easy, and
so we work toward showing $\subseteq$. Let $z\in\phi_\mathcal{T}
(\gamma_{\mathcal{T}}(x,y)\setminus\{x,y\} )$ for some $x\neq y$,
and suppose that it is also the case that $z=\phi_\mathcal{T}(x)$. By
assumption, we know that $z=\phi_\mathcal{T}(x')$ for some $x'\in\gamma
_{\mathcal{T}}(x,y)\setminus\{x,y\}$. Now, by Lemma~\ref{push},
because $\phi_\mathcal{T}(x')=\phi_\mathcal{T}(x)$ we can apply
Lemma~\ref{schrammlem}(a) to deduce that $\phi_\mathcal{T}(\gamma
_{\mathcal{T}}(x,x'))$ is either a single point or homeomorphic to a
circle. Actually, since $d_\mathcal{T}(x,x')>0$, Lemma~\ref{homeo}
tells us that $d_\mathcal{T}^S(x,x')>0$, and so it must be the latter
option that holds true. We continue to consider two cases. First, if
$\phi_\mathcal{T}(y)\neq\phi_\mathcal{T}(x)$, then
Lemma~\ref{schrammlem}(a) tells us that $\phi_\mathcal{T}(\gamma
_{\mathcal{T}}(x,y))$ must be a simple path. However, the circle $\phi
_\mathcal{T}(\gamma_{\mathcal{T}}(x,x'))$ is a subset of $\phi_\mathcal
{T}(\gamma_{\mathcal{T}}(x,y))$, and so we arrive at a contradiction.
Second, if $\phi_\mathcal{T}(y)=\phi_\mathcal{T}(x)$, then one can
again apply Lemma~\ref{homeo} to choose $x''\in\gamma_{\mathcal
{T}}(x,y)$ such that $\phi_\mathcal{T}(x'')\neq\phi_\mathcal{T}(x)$.
Clearly, we have that either $z\in\phi_\mathcal{T}(\gamma_{\mathcal
{T}}(x,x'')\setminus\{x,x''\})$ or $z\in\phi_\mathcal{T}(\gamma
_{\mathcal{T}}(x'',y)\setminus\{x'',y\})$. Since $x''\notin\{x,y\}$
and $\phi_\mathcal{T}(x'')\notin\break \{\phi_\mathcal{T}(x),\phi_\mathcal
{T}(y)\}$, the situation reduces to the first case, and yields another
contradiction. Hence, we cannot have $z=\phi_\mathcal{T}(x)$.
Similarly, $z\neq\phi_\mathcal{T}(y)$, so the claim at (\ref{claim})
is proved.

Now, from (\ref{trunkdef}), (\ref{t0def}) and (\ref{claim}), it is
clear that
\begin{eqnarray*}
\phi_\mathcal{T}\bigl(\mathcal{T}^o\bigr)&=&\bigcup
_{x,y\in\mathcal{T}}\phi_\mathcal {T} \bigl(\gamma_{\mathcal{T}}(x,y)
\setminus\{x,y\} \bigr)=\bigcup_{x,y\in\mathcal{T}}
\phi_\mathcal{T} \bigl(\gamma_{\mathcal{T}}(x,y) \bigr)\setminus\bigl\{
\phi _\mathcal{T}(x),\phi_\mathcal{T}(y)\bigr\}\\
&=& \mathrm{trunk},
\end{eqnarray*}
and so the map $\phi_\mathcal{T}:\mathcal{T}^o\rightarrow\mathrm{ trunk}$
is a surjection. To complete the proof, we will again use the fact
that, for every $x,y\in\mathcal{T}$ with $\phi_\mathcal{T}(x)\neq\phi
_\mathcal{T}(y)$, $\phi_\mathcal{T}(\gamma_{\mathcal{T}}(x,y))$ is a
simple path, and note that the proof of this result in \cite{Schramm}
includes showing that the endpoints of this path are $\phi_\mathcal
{T}(x)$ and $\phi_\mathcal{T}(y)$. In particular, if $x,y\in\mathcal
{T}^o$ are such that $\phi_\mathcal{T}(x)\neq\phi_\mathcal{T}(y)$, then
we know that the simple path $\phi_\mathcal{T}(\gamma_{\mathcal
{T}}(x,y))$ from $\phi_\mathcal{T}(x)$ to $\phi_\mathcal{T}(y)$ is
contained in $\mathrm{ trunk}$. Recalling that $\mathrm{ trunk}$ is a
topological tree, which implies there is be a unique path $\gamma_\mathrm{
trunk}(\phi_\mathcal{T}(x),\phi_\mathcal{T}(y))$ between $\phi_\mathcal
{T}(x)$ and $\phi_\mathcal{T}(y)$ within this set, it must be the case
that $\phi_\mathcal{T}(\gamma_{\mathcal{T}}(x,y))=\gamma_\mathrm{
trunk}(\phi_\mathcal{T}(x),\phi_\mathcal{T}(y))$. On the other hand, if
$x,y\in\mathcal{T}^o$ are such that $\phi_\mathcal{T}(x)=\phi_\mathcal
{T}(y)$, then, by
Lemma~\ref{schrammlem}(a), it must hold that $\phi_\mathcal{T}(\gamma
_{\mathcal{T}}(x,y))=\{\phi_\mathcal{T}(x)\}$, where we note that we
can exclude the possibility that $\phi_\mathcal{T}(\gamma_{\mathcal
{T}}(x,y))$ is homeomorphic to a circle, since $\mathrm{ trunk}$ is a
topological tree and cannot contain such a subset. Hence we obtain that
$\phi_\mathcal{T}(\gamma_{\mathcal{T}}(x,y))=\gamma_\mathrm{ trunk}(\phi
_\mathcal{T}(x),\phi_\mathcal{T}(y))$ in this case, also. Consequently,
$d_\mathcal{T}^S(x,y)=d_\mathrm{ trunk}^S(\phi_\mathcal{T}(x),\break \phi_\mathcal
{T}(y))$ for every $x,y\in\mathcal{T}^o$. This confirms that $\phi
_\mathcal{T}$ is an isometry, as desired.
\end{pf}

Before we complete the proof of Theorem~\ref{main12}, we mention
another property of the trunk
that will be needed. This is that
the trunk can be used to reconstruct the \emph{dual trunk}, that is, the
(subsequential)
scaling limit of the dual graph of the UST (see \cite{Schramm}, Remarks 10.13 and
10.14).
More precisely, for any two points $x,y\in\mathbb{S}^2\setminus\mathrm{
trunk}$, there
exists a unique path in $\mathbb{S}^2\setminus\mathrm{ trunk}$ between
them. Denote this path by $\gamma_\mathrm{ trunk^\dagger}(x,y)$, and set
$\mathrm{ trunk}^\dagger:=\bigcup_{x,y\in\mathbb{S}^2\setminus\mathrm{ trunk}}
\gamma_\mathrm{ trunk^\dagger}(x,y)\setminus\{x,y\}$. This is the dual
trunk, which
is distributed identically to $\mathrm{ trunk}$.

\begin{pf*}{Proof of Theorem~\ref{main12}} It readily follows from
Theorem~\ref{main1} and the unboundedness of $(\sU, d_{\sU})$ that the
diameter of $(\mathcal{T},d_\mathcal{T})$ is infinite, and so it has at
least one end at infinity. Thus, to complete the proof of part (a)(ii),
it will suffice to show that there can be no more than one end at
infinity. To this end, note that, for any $r>0$,
\[
\tilde{\mathbf{P}} \bigl((\mathcal{T},d_\mathcal{T})\mbox{ has $\geq$ 2
ends at infinity} \bigr)\leq\lim_{R\rightarrow\infty}\tilde{\mathbf {P}}
\bigl(C_\mathcal{T}(r,R)\geq2 \bigr),
\]
where $C_\mathcal{T}(r,R)$ is the event that there exist $x,y\notin
\mathcal{T}^{(R)}$ such that $\gamma_\mathcal{T}(x,y)\cap\mathcal
{T}^{(r)}\neq\varnothing$. By applying the coupling of Lemma~\ref{couplem},
it is possible to bound the inner probability by $\limsup_{\delta\rightarrow0}\mathbf{P}(C_\mathcal{U}(2\delta^{-\kappa}r,\delta
^{-\kappa}R/2))$, where $C_\mathcal{U}(r,R)$ is defined similarly to
$C_\mathcal{T}(r,R)$, with $\mathcal{T}$ replaced by $\sU$. (Since we
have already presented similar coupling arguments in the proofs of
Lemmas \ref{measlem}, \ref{dsbound} and \ref{push}, we omit the
details.) Now, for $\lambda>0$,
\begin{eqnarray*}
&&\mathbf{P} \bigl(C_\mathcal{U}\bigl(2\delta^{-\kappa}r,
\delta ^{-\kappa}R/2\bigr) \bigr)
\\
&&\qquad\leq \mathbf{P} \bigl(C^E_\mathcal{U}\bigl(\lambda
\delta^{-1}r^{1/\kappa},\delta ^{-1}R^{1/\kappa}/
\lambda\bigr) \bigr) +\mathbf{P} \bigl(B_\sU\bigl(0,2
\delta^{-\kappa}r\bigr)\nsubseteq B_E\bigl(0,\lambda
\delta^{-1}r^{1/\kappa}\bigr) \bigr)
\\
&& \qquad\quad{}+ \mathbf{P} \bigl(B_E\bigl(0,\lambda^{-1}
\delta^{-1}R^{1/\kappa}\bigr)\nsubseteq B_\sU\bigl(0,
\delta^{-\kappa}R/2\bigr) \bigr),
\end{eqnarray*}
where $C^E_\mathcal{U}(r,R)$ is the event that there exist $x,y\notin
B_E(0,R)$ such that $\gamma(x,y)\cap B_E(0,r)\neq\varnothing$. Hence,
from Theorem~\ref{vollemust}(a) and \cite{ABNW}, we obtain, for $R\geq
\lambda\geq2$ and $R\geq\lambda^{2\kappa} r$,
\[
\mathbf{P} \bigl(C_\mathcal{U}\bigl(2\delta^{-\kappa}r,
\delta^{-\kappa
}R/2\bigr) \bigr)\leq c_1 \biggl(
\frac{\lambda^2r^{1/\kappa}}{R^{1/\kappa
}} \biggr)^{c_2}+c_3\lambda^{-1/6}.
\]
Taking $\lambda=R^{1/3\kappa}$, this converges to 0 as $R\rightarrow
\infty$. It follows that the $\tilde{\mathbf{P}}$-probability of
$(\mathcal{T},d_\mathcal{T})$ having $\geq$ 2 ends at infinity is zero,
as desired.

For part (b)(i), we begin by noting that
$A\subseteq\phi_\mathcal{T}^{-1}(\phi_\mathcal{T}(A))$ for any set $A$.
Consequently, from Lemma~\ref{measlem}, we obtain $\mu_{\mathcal{T}}(\{
x\})\leq\mathcal{L}(\{\phi_\mathcal{T}(x)\})=0$ for every $x\in\mathcal
{T}$. Moreover, by Lemma~\ref{trunklem}, $\mu_{\mathcal{T}}({\mathcal
{T}^o})\leq\mathcal{L}(\phi_\mathcal{T}(\mathcal{T}^o))=\mathcal
{L}(\mathrm{ trunk})=0$,
where the final equality is a consequence of the fact that $\mathrm{
trunk}$ has Hausdorff dimension strictly less than two [as recalled in
Lemma~\ref{schrammlem}(d)].

To establish part (b)(ii), it will suffice to check that, given $R>0$,
there exist constants $c_1,c_2,c_3,c_4,c_5,c_6\in(0,\infty)$ such that,
for every $r\in(0,1)$,
\begin{eqnarray*}
\tilde{\mathbf{P}} \Bigl(\inf_{x\in B_\mathcal{T}(\rho_\mathcal{T},R)}\mu _\mathcal{T}
\bigl(B_\mathcal{T}(x,r) \bigr) \le c_1r^{d_f}\bigl(
\log r^{-1}\bigr)^{-80} \Bigr)&\leq& c_2r^{c_3},
\\
\tilde{\mathbf{P}} \Bigl(\sup_{x\in B_\mathcal{T}(\rho_\mathcal{T},R)}\mu _\mathcal{T}
\bigl(B_\mathcal{T}(x,r) \bigr) \geq c_4r^{d_f}\bigl(
\log r^{-1}\bigr)^{80} \Bigr)&\leq& c_5r^{c_6}.
\end{eqnarray*}
Indeed, one can then apply a simple Borel--Cantelli argument along the
subsequence $r_n=2^{-n}$, $n\in\mathbb{N}$, to deduce the result. We
observe that the above inequalities can be deduced from the definition
of $\tilde{\mathbf{P}}$ and Corollary~\ref{P:measprop} by applying the
coupling of Lemma~\ref{couplem}. Furthermore, note that part (a)(i) is an elementary
consequence of (b)(ii) (see \cite{Edgar}, Proposition~1.5.15, e.g.).

Part (b)(iii) can also be obtained using a Borel--Cantelli argument in
conjunction with the following: there exist constants $c_1,c_2\in
(0,\infty)$ such that
%
\begin{eqnarray}
\label{bfirst} \tilde{\mathbf{P}} \bigl(\mu_{\mathcal{T}} \bigl(B_\mathcal{T}(
\rho _\mathcal{T},r) \bigr)\geq\lambda r^{d_f} \bigr)&\leq&
c_1e^{-c_2\lambda^{1/3}},
\\
\label{bsecond} \tilde{\mathbf{P}} \bigl(\mu_{\mathcal{T}} \bigl(B_\mathcal{T}(
\rho _\mathcal{T},r) \bigr)\leq\lambda^{-1} r^{d_f}
\bigr)&\leq& c_1e^{-c_2\lambda^{1/9}},
\end{eqnarray}
for all $r>0$, $\lambda\geq1$. Again applying the definition of $\tilde
{\mathbf{P}}$ and the coupling of
Lemma~\ref{couplem}, it is possible to deduce the bound at (\ref
{bfirst}) from
Theorem~\ref{vollemust}(b). The proof for the bound at (\ref{bsecond})
is similar.

The first statement of part (c)(i) depends on Lemmas \ref{homeo} and
\ref{trunklem}. In particular, these two results imply that $\phi
_\mathcal{T}$ is a homeomorphism from $(\mathcal{T}^o,d_\mathcal{T})$
to $(\mathrm{ trunk},d_\mathrm{ trunk}^S)$. To replace the topology generated
by $d_\mathrm{ trunk}^S$ with the Euclidean one, we will show that the
identity map from $(\mathrm{ trunk},d_\mathrm{ trunk}^S)$ to $(\mathrm{
trunk},d_E)$ is also a homeomorphism. Clearly, it is a continuous
bijection, and so we need to show its inverse is continuous. To do
this, suppose that $x_n, x\in\mathrm{ trunk}$ are such that
$d_E(x_n,x)\rightarrow0$. Now, in light of $\phi_\mathcal{T}:(\mathcal
{T}^o,d_\mathcal{T})\rightarrow(\mathrm{ trunk},d_\mathrm{ trunk}^S)$ being a
homeomorphism, the map $\phi_\mathcal{T}:\mathcal{T}\rightarrow\mathbb
{R}^2$ can be viewed as the extension of the identity map $\mathrm
{trunk}\rightarrow\mathbb{R}^2$ to a continuous map on the completion
of $(\mathrm{trunk},d_\mathrm{trunk}^S)$, and it therefore follows from
the discussion in \cite{Schramm}, Remark~10.15, that $|\phi_\mathcal
{T}^{-1}(x)|$ is equal to one if $x$ is not contained in $\mathrm{
trunk}^\dagger$. In particular, since $\mathrm{ trunk}\cap\mathrm{
trunk}^\dagger=\varnothing$, there exist unique $y_n,y\in\mathcal{T}$
such that $\phi_\mathcal{T}(y_n)=x_n$ and $\phi_\mathcal{T}(y)=x$.
Moreover, since $x_n,x\in B_E(0,r)$ for some $r<\infty$, there must
exist an $R<\infty$ such that $y_n,y\in B_\mathcal{T}(\rho_\mathcal
{T},R)$---this is an easy consequence of (\ref{toprove}). Hence, by
compactness, for any subsequence $n_i$, there exists a sub-subsequence
$y_{{n_i}_j}$ such that $d_\mathcal{T}(y_{{n_i}_j},y')\rightarrow0$
for some $y'\in\mathcal{T}$. By the continuity of $\phi_\mathcal{T}$,
it follows that $d_E(\phi_\mathcal{T}(y'),\phi_\mathcal{T}(y))=\lim_{j\rightarrow\infty}
d_E(\phi_\mathcal{T}(y_{{n_i}_j}),\phi_\mathcal{T}(y))=\lim_{j\rightarrow\infty}
d_E(x_{{n_i}_j},x)=0$, and so $y'=y$. Noting that $y_n,y$ are
necessarily in $\mathcal{T}^o$, Lemmas \ref{homeo} and \ref{trunklem}
thus yield $d_\mathrm{ trunk}^S(x_{{n_i}_j},x)=d_{\mathcal
{T}}^S(y_{{n_i}_j},y)\rightarrow0$. Since the initial subsequence
$(n_i)$ was arbitrary, this implies $d_\mathrm{
trunk}^S(x_{{n}},x)\rightarrow0$, as desired. The denseness of $\phi
_\mathcal{T}(\mathcal{T}^o)$ in $\mathbb{R}^2$ follows from
Lemmas \ref{schrammlem}(b) and \ref{trunklem}. Furthermore, applying
Lemma~\ref{schrammlem}(c) together with the homeomorphism between
$\mathcal{T}^o$ and $\mathrm{ trunk}$ yields that $\max_{x\in\mathcal
{T}}\mathrm{ deg}_\mathcal{T}(x)=3$, which is the first claim of part
(c)(ii). To check the remaining claim of part (c)(ii), we note that if
$x$ is contained in $\mathrm{ trunk}^\dagger$, then $|\phi_\mathcal
{T}^{-1}(x)|$ is equal to the degree of $x$ in $\mathrm{ trunk}^\dagger$
(again, see the discussion in \cite{Schramm}, Remark~10.15). Since
$\mathrm{ trunk}^\dagger$ also has maximum degree 3 [by
Lemma~\ref{schrammlem}(c) again], this establishes the desired result.
[Recall that $|\phi_\mathcal{T}^{-1}(x)|=1$ for $x\notin\mathrm{
trunk}^\dagger$.] Finally, point (c)(iii) will be a simple consequence
of Lemma~\ref{measlem}, at least if we can show that $\mathcal{L}(\{
x:|\phi_\mathcal{T}^{-1}(x)|\geq2\})=0$. However, from our previous
observations, we know that the set $\{x:|\phi_\mathcal{T}^{-1}(x)|\geq
2\}$ is contained in $\mathrm{ trunk}^\dagger$, which has Hausdorff
dimension strictly less than two
[by Lemma~\ref{schrammlem}(d)].
\end{pf*}

\section{Simple random walk scaling limit}\label{srwsec}

In this section, we prove Theorem~\ref{main2}. The general convergence
result for simple random walks on graph trees (see Theorem~\ref{srwthm}
below) that we apply extends \cite{Croy1,Croy2,Croy3} from the setting
of ordered graph trees (in particular, in those articles graph trees
and real trees were encoded by functions). Work is also needed to
extend to the noncompact setting of this article.

Let us start by introducing some notation. Let $(T_n)_{n\geq1}$ be a
sequence of finite graph trees.
Write $d_{T_n}$ for the shortest path graph distance on $T_n$, and $\mu
_{T_n}$ for the counting measure on the vertices of $T_n$. Suppose that
$\phi_n$ is a map from the vertices of $T_n$ into $M$---until
otherwise noted, we assume that $M$ is a separable normed vector space,
and write the metric induced by its norm as $d_M$. Fix a distinguished
vertex $\rho_{T_n}$ of $T_n$. We extend $({T}_n,d_{{T}_n},\mu
_{{T}_n},\phi_{{T}_n},\rho_{{T}_n})$ to an element of $\mathbb{T}_c$ by
adding line segments of unit length along edges of the graph tree, and
(isometrically) interpolating $\phi_{T_n}$ between these along the
relevant geodesics. The process $(X^{T_n}_t)_{t\geq0}$ is the discrete
time simple random walk on $T_n$, and $P^{T_n}_x$ its law started from
$x$. We extend $(\phi_{{T}_n}(X^{{T}_n}_{t}))_{t\geq0}$ to an element
of $C(\mathbb{R}_+,M)$ by interpolation along geodesics.

The limit space we consider in Theorem~\ref{srwthm} is the natural
generalisation of that of \cite{Croy3}. In particular, let $\mathbb
{T}_c^*$ be the collection of those elements $\uT$ of $\mathbb{T}_c$
such that $\mu_{\mathcal{T}}$ is nonatomic, supported on the leaves of
$\mathcal{T}$ [recall that the leaves of a real tree $\mathcal{T}$ are
those points $x\in\mathcal{T}$ such that $\mathcal{T}\setminus\{x\}$
is connected, i.e., which have $\operatorname{ deg}_\mathcal{T}(x)=1$], and
also there exists a constant $c$ such that
%
\begin{equation}
\label{measlower} \liminf_{r\rightarrow0}{\inf_{x\in\mathcal{T}}r^{-c}
\mu_\mathcal {T} \bigl(B_\mathcal{T}(x,r) \bigr)}>0.
\end{equation}

For a locally compact real tree $(\mathcal{T},d_\mathcal{T})$ equipped
with a locally finite Borel measure $\mu_\mathcal{T}$ of full support,
it is shown in \cite{AEW} how to construct an associated ``Brownian
motion'' (cf. \cite{Kigami}, which deals with the case when $(\mathcal
{T},d_\mathcal{T})$ is complete). For readers' convenience, let us
briefly summarise this construction. In particular, first define the
length measure $\lambda_\mathcal{T}$ on $\mathcal{T}$ to be the
restriction of one-dimensional Hausdorff measure to $\mathcal
{T}^o:=\mathcal{T}\setminus\{x\in\mathcal{T}: \operatorname{ deg}_\mathcal
{T}(x)=1\}$.
Moreover, let $\mathcal{A}$ be the collection of locally absolutely
continuous functions on $\mathcal{T}$, where we say a function
$f:\mathcal{T}\rightarrow\mathbb{R}$ is locally absolutely continuous
if and only if for every $\varepsilon>0$ and subset $A\subseteq\mathcal
{T}$ with $\lambda_\mathcal{T}(A)<\infty$, there exists a $\delta>0$
such that: if $\gamma_\mathcal{T}(x_1,y_1),\ldots,\gamma_\mathcal
{T}(x_k,y_k)\subseteq A$, $k\in\bN$
are disjoint arcs with $\sum_{i=1}^kd_\mathcal{T}(x_i,y_i)<\delta$, then
$\sum_{i=1}^k|f(x_i)-f(y_i)|<\varepsilon$.
Given a function $f\in\mathcal{A}$, there exists a unique (up to
$\lambda_\mathcal{T}$-null sets) function $g$ that is locally in
$L^1(\mathcal{T},\lambda_\mathcal{T})$ such that
\[
f(y)-f(x) = - \int_{\gamma_\mathcal{T}(b_\mathcal{T}(\rho_\mathcal
{T},x,y),x)}{g(z)}\lambda_\mathcal{T}(dz)+
\int_{\gamma_\mathcal{T}(b_\mathcal{T}(\rho_\mathcal
{T},x,y),y)}{g(z)}\lambda_\mathcal{T}(dz),
\]
for all $x,y\in\mathcal{T}$, where $b_\mathcal{T}(\rho_\mathcal
{T},x,y)$ is the unique branch-point of $\rho_\mathcal{T}$, $x$ and $y$
in $\mathcal{T}$ \cite{AEW}, Proposition~1.1. (Note that the above
difference can be interpreted as an oriented integral from $x$ to $y$,
and is independent of the choice of $\rho_\mathcal{T}$.) The function
$g$ in the previous sentence is called the gradient of $f$, and is
denoted $\nabla f$. Next, define
\[
\mathcal{E}_\mathcal{T}(f,g):=\frac{1}{2}\int\nabla f (x)\nabla
g(x)\lambda_\mathcal{T}(dx),
\]
for all $f,g\in\mathcal{F}_\mathcal{T}$, where $\mathcal{F}:=\{f\in
\mathcal{A}: \nabla f\in L^2(\mathcal{T},\lambda_\mathcal{T})\}\cap
L^2(\mathcal{T},\mu_\mathcal{T}) \cap C_\infty(\mathcal{T})$. Here,
$C_\infty(\mathcal{T})$ is the space of continuous functions on
$\mathcal{T}$ that vanish at infinity.
By \cite{AEW}, Propositions 2.4 and 4.1, and the proof of Theorem~1,
$(\mathcal{E}_\mathcal{T},\mathcal{F}_\mathcal{T})$ is a local, regular
Dirichlet form on
$L^2(\mathcal{T},\mu_\mathcal{T})$.
The Brownian motion on $(\mathcal{T},d_\mathcal{T},\mu_\mathcal{T})$
is the continuous, $\mu_\mathcal{T}$-symmetric, strong Markov process
$((X^\mathcal{T}_t)_{t\geq0},
(P^\mathcal{T}_x)_{x\in\mathcal{T}})$ associated with this Dirichlet
form (see \cite{FOT}).
Clearly, this construction applies to elements of $\mathbb{T}_c^*$, and
the additional restriction
(\ref{measlower}) allows one to deduce that this Brownian motion has
local times\break 
$(L^\mathcal{T}_t(x))_{x\in\mathcal{T},t\geq0}$ which
are jointly continuous in $t$ and $x$ (see \cite{Croy3}, Lemma~2.2).

\begin{thmm}\label{srwthm}
Let $(a_n)_{n\geq1}$, $(b_n)_{n\geq1}$, $(c_n)_{n\geq1}$ be null sequences
with $b_n=o(a_n)$ such that
%
\begin{equation}
\label{tnconv} ({T}_n,a_nd_{{T}_n},b_n
\mu_{{T}_n},c_n\phi_{{T}_n},\rho _{{T}_n} )
\rightarrow\uT
\end{equation}
in $(\mathbb{T}_c,\Delta_c)$, where $\uT$ is an element of $\mathbb
{T}^*_c$. Then
%
\begin{equation}
\label{srwconv} \bigl(c_n\phi_{T_n} \bigl(X^{T_n}_{t/a_nb_n}
\bigr) \bigr)_{t\geq
0}\rightarrow \bigl(\phi_{\mathcal{T}}
\bigl(X^{\mathcal{T}}_{t} \bigr) \bigr)_{t\geq0}
\end{equation}
in distribution in $C(\mathbb{R}_+,M)$, where we assume $X^{T_n}_0=\rho
_{T_n}$ for each $n$, and also $X^\mathcal{T}_0=\rho_\mathcal{T}$.
\end{thmm}

Since the proof of this result is close to the arguments of \cite
{Croy1,Croy2,Croy3}, we will not
give all of the details.
For clarity, though, we will break it into three lemmas. The basic idea
is to approximate the processes
of interest by processes on trees which have finite total length, for
which convergence is more
straightforward to prove. So, let $\uT$ be an element of $\mathbb
{T}^*_c$, and $(x_i)_{i\geq1}$ be a dense sequence of vertices in
$\mathcal{T}$---these will be fixed throughout the current discussion.
(To avoid trivialities, we assume that $\mathcal{T}$ consists of more
than one point.) We suppose that $(x_i)_{i\geq1}$ are distinct, and
none is equal to the root $\rho_\mathcal{T}$, which we will sometimes
also denote by $x_0$. For each $k\geq1$, define $\mathcal{T}(k):=\bigcup_{i=1}^k\gamma_\mathcal{T}(\rho_\mathcal{T},x_i)$, and let $\pi
_{\mathcal{T},\mathcal{T}(k)}$ be the natural projection from $\mathcal
{T}$ to $\mathcal{T}(k)$, that is, for $x\in\mathcal{T}$, $\pi_{\mathcal
{T},\mathcal{T}(k)}(x)$ is the closest point in $\mathcal{T}(k)$ to
$x$. Taking $\mu^{(k)}:=\mu_\mathcal{T}\circ\pi_{\mathcal{T},\mathcal
{T}(k)}^{-1}$, we define $X^{\mathcal{T}(k),{\mu^{(k)}}}$ to be
Brownian motion on $(\mathcal{T}(k),d_\mathcal{T}|_{\mathcal{T}(k)},\mu
^{(k)})$. By \cite{Croy3}, Proposition~2.1, if we assume that
$X^{\mathcal{T}(k),\mu^{(k)}}$ and $X^\mathcal{T}$ are both started
from $\rho_\mathcal{T}$, then
$(X^{\mathcal{T}(k),\mu^{(k)}}_t)_{t\geq0}\rightarrow(X^{\mathcal
{T}}_t)_{t\geq0}$ in distribution in $C(\mathbb{R}_+,\mathcal{T})$.
(This step is one of the places in the proof that the existence of
jointly continuous local times for the Brownian motion $X^\mathcal{T}$
is used.) Hence, the continuous mapping theorem implies that
%
\begin{equation}
\label{almost} \bigl(\phi_\mathcal{T} \bigl(X^{\mathcal{T}(k),\mu^{(k)}}_t
\bigr) \bigr)_{t\geq0}\rightarrow \bigl(\phi_{\mathcal{T}}
\bigl(X^{\mathcal
{T}}_t \bigr) \bigr)_{t\geq0}
\end{equation}
in distribution in $C(\mathbb{R}_+,M)$. Moreover, if we define
$\mathcal{B}(k):=\{b_\mathcal{T}(\rho_\mathcal{T},x_i,x_j): i,j\in\{
0,\ldots,k\}\}$ and $\phi^{(k)}_{\mathcal{T}}:\mathcal{T}(k)\rightarrow
M$ by setting $\phi^{(k)}_\mathcal{T}=\phi_\mathcal{T}$ on $\mathcal
{B}(k)$ and interpolating along geodesics between these vertices, then
one can deduce from the denseness of $(x_i)_{i\geq1}$ and continuity
of $\phi_\mathcal{T}$ that\break  $\lim_{k\rightarrow\infty}\sup_{x\in\mathcal
{T}(k)}d_M(\phi_\mathcal{T}^{(k)}(x),\phi_\mathcal{T}(x))=0$
(cf. \cite{Croy2}, Theorem~8.2, and the following discussion).
Consequently, (\ref{almost}) yields the following lemma.

\begin{lem}\label{t1} As $k\rightarrow\infty$,
\[
\bigl(\phi_\mathcal{T}^{(k)} \bigl(X^{\mathcal{T}(k),\mu^{(k)}}_t
\bigr) \bigr)_{t\geq0} \rightarrow \bigl(\phi_{\mathcal{T}}
\bigl(X^{\mathcal{T}}_t \bigr) \bigr)_{t\geq0}
\]
in distribution in $C(\mathbb{R}_+,M)$.
\end{lem}

Before describing the connection with discrete objects, let us note
that $X^{\mathcal{T}(k),\mu^{(k)}}$ can also be represented as a time
change of another Brownian motion on $\mathcal{T}(k)$. In particular,
let $\lambda^{(k)}$ be the one-dimensional Hausdorff measure on
$(\mathcal{T}(k),d_\mathcal{T}|_{\mathcal{T}(k)})$, and $X^{\mathcal
{T}(k),\lambda^{(k)}}$ be the associated Brownian motion. Since $\lambda
^{(k)}$ satisfies (\ref{measlower}), this process admits jointly
continuous local times $(L^{(k)}_t(x))_{x\in\mathcal{T}(k),t\geq0}$,
from which we define an additive functional $\hat{A}^{(k)}_t:=\int_{\mathcal{T}(k)}L_t^{(k)}(x)\mu^{(k)}(dx)$,
and its inverse $\hat{\tau}^{(k)}(t):=\inf\{s: \hat{A}^{(k)}_s>t\}$.
(We use hatted notation for consistency with \cite{Croy3}.) From
\cite{Croy3}, Lemma~2.4, we then obtain that if $X^{\mathcal
{T}(k),\lambda^{(k)}}$ is started from $\rho_\mathcal{T}$, then
%
\begin{equation}
\label{xlam} \bigl(X^{\mathcal{T}(k),\lambda^{(k)}}_{\hat{\tau}^{(k)}(t)} \bigr)_{t\geq0}
\end{equation}
is distributed identically to $X^{\mathcal{T}(k),\mu^{(k)}}$ started
from $\rho_\mathcal{T}$.

For the next part of the proof of Theorem~\ref{srwthm}, we fix a
sequence of metric spaces $Z_n$, isometric embeddings $\psi_n:\mathcal
{T}\rightarrow Z_n$, $\psi'_n:({T}_n,a_nd_{T_n})\rightarrow Z_n$ and
correspondences $\mathcal{C}_n$ between $\mathcal{T}$ and ${T}_n$
containing $(\rho_\mathcal{T},\rho_{{T}_n})$ and such that (\ref{yy})
holds with $\uT_n$ replaced by $({T}_n,a_nd_{{T}_n},b_n\mu
_{{T}_n},c_n\phi_{{T}_n},\rho_{{T}_n})$ for some sequence $\varepsilon
_n\rightarrow0$. [This is possible if we suppose that (\ref{tnconv})
holds.] Moreover, let $x_i^n\in T_n$ be such that $(x_i,x_i^n)\in
\mathcal{C}_n$, and define the subtree $T_n(k)\subseteq T_n$ and
projection $\pi_{n,k}:T_n\rightarrow T_n(k)$ similarly to the
continuous case. Using elementary arguments, one can check that
%
\begin{eqnarray}
\label{deltank} &&\lim_{k\rightarrow\infty}\limsup_{n\rightarrow\infty}a_n
\max_{x\in
T_n}d_{T_n}\bigl(x,\pi_{n,k}(x)\bigr)
\nonumber
\\[-8pt]
\\[-8pt]
\nonumber
&&\qquad\leq\lim_{k\rightarrow\infty}\limsup_{n\rightarrow\infty} \Bigl(\sup
_{x\in\mathcal{T}}d_\mathcal{T}\bigl(x,\pi _{\mathcal{T},\mathcal{T}(k)}(x)
\bigr)+5\varepsilon_n \Bigr)=0
\end{eqnarray}
(cf. \cite{Croy1}, Lemma~2.7), which says that the subtrees $T_n(k)$
are uniformly good approximations of the full trees $T_n$. As for
processes, we define $X^{n,k}:=\pi_{n,k}(X^{T_n})$ and $J^{n,k}$ to be
the corresponding jump chain, that is, if $A^{n,k}_0:=0$ and
$A^{n,k}_t:=\min\{s\geq A_{t-1}^{n,k}: X^{T_n}_s\in T_n(k)\setminus\{
X^{T_n}_{A_{t-1}^{n,k}}\}\}$, then $J^{n,k}_t=X^{T_n}_{A^{n,k}_t}$.
Conversely, if $\tau^{n,k}(t):=\max\{s: A^{n,k}_s\leq t\}$, then we
can write
%
\begin{equation}
\label{xnk} X^{n,k}_t=J^{n,k}_{\tau^{n,k}(t)}.
\end{equation}
Define the local times of $J^{n,k}$ by setting $L^{n,k}_t(x):=\frac
{2}{\operatorname{ deg}_{n,k}(x)}\sum_{s=0}^{t-1}\mathbf{1}_x(J_s^{n,k})$ for
$x\in T_n(k)$ and $t\geq0$, where $\operatorname{ deg}_{n,k}(x)$ is the usual
graph degree of $x$ in $T_n(k)$. We use these to define an associated
additive functional $\hat{A}^{n,k}$ by setting $\hat{A}^{n,k}_0:=0$ and
$\hat{A}^{n,k}_t:=\int_{T_n(k)}L^{n,k}_{t}(x)\mu_{n,k}(dx)$, where $\mu
_{n,k}:=\mu_{T_n}\circ\pi_{n,k}^{-1}$. Finally, from the inverse $\hat
{\tau}^{n,k}(t):=\max\{s: \hat{A}^{n,k}_s\leq t\}$, we define an
alternative time-change of $J^{n,k}$ by setting
%
\begin{equation}
\label{hatxnk} \hat{X}^{n,k}_t=J^{n,k}_{\hat{\tau}^{n,k}(t)}.
\end{equation}
In the next lemma, we describe how methods of \cite{Croy3}, Section~3,
can be applied to deduce a scaling limit for these processes. The map
$\phi_{T_n}^{(k)}:T_n(k)\rightarrow M$ is
defined analogously to the continuous case.

\begin{lem}\label{t0} Suppose that (\ref{tnconv}) holds, and fix $k\geq1$. Then
\[
\bigl(c_n\phi_{T_n}^{(k)} \bigl(
\hat{X}^{n,k}_{t/a_nb_n} \bigr) \bigr)_{t\geq0}\rightarrow
\bigl(\phi_\mathcal{T}^{(k)} \bigl(X^{\mathcal
{T}(k),\mu^{(k)}}_t
\bigr) \bigr)_{t\geq0},
\]
in distribution in $C(\mathbb{R}_+,M)$.
\end{lem}

\begin{pf} In this proof, we will use an embedding of trees into $\ell^1$,
the Banach space of infinite sequences of real numbers equipped with
the metric $d_{\ell^1}$ induced by the norm $\sum_{i\geq1}|x(i)|$ for
$x\in\ell^1$ (the procedure was originally described in \cite
{Aldous3}). In particular, given a sequence $(\mathcal{T}(k))_{k\geq
1}$ as above it is possible to construct a distance-preserving map
$\tilde{\psi}:(\mathcal{T},d_\mathcal{T})\rightarrow(\ell^1, d_{\ell
^1})$ that satisfies $\tilde{\psi}(\rho)=0$ and
%
\begin{equation}
\label{psichar} \pi_k\bigl(\tilde{\psi}(\sigma)\bigr)=\tilde{\psi}
\bigl(\pi_{\mathcal{T},\mathcal
{T}(k)}(\sigma)\bigr)
\end{equation}
for every $\sigma\in\mathcal{T}$ and $k\geq1$, where $\pi_k$ is the
projection map on $\ell^1$, that is, $\pi_k(x(1),x(2),\ldots)=(x(1),\ldots
,x(k),0,0,\ldots)$. Roughly speaking, we first map $\mathcal{T}(1)$ to a
line segment of length $d_\mathcal{T}(\rho_\mathcal{T},x_1)$ in the
first coordinate direction of~$\ell^1$. Then, given the map $\tilde{\psi
}$ on $\mathcal{T}(k)$, map the additional line segment in $\mathcal
{T}(k+1)$ to a line segment in the $(k+1)$st coordinate direction of
$\ell^1$ [i.e., orthogonally to $\tilde{\psi}(\mathcal{T}(k)$)],
attached at the image in $\ell^1$ of the appropriate branch-point. Such
a map is determined uniquely by insisting that $\tilde{\psi}(\mathcal
{T})\subseteq\{(x(1),x(2),\ldots)\in\ell^1:x(i)\geq0, i=1,2,\ldots\}$.
We can of course embed the discrete trees similarly, and we denote the
corresponding embeddings by $\tilde{\psi}_n$. It is not difficult to
check from our construction that, for every $i\geq0$,
%
\begin{equation}
\label{c1} a_n\tilde{\psi}_n\bigl(x_i^n
\bigr)\rightarrow\tilde{\psi}(x_i).
\end{equation}
As a consequence of this and the fact that the maps $\tilde{\psi}$ and
$\tilde{\psi}_n$ are isometries, we find that
$({T}_n,a_nd_{{T}_n},b_n\mu_{{T}_n},a_n\tilde{\psi}_{n},\rho_{{T}_n})
\rightarrow(\mathcal{T},d_\mathcal{T},\mu_{\mathcal{T}},\tilde{\psi
},\rho_\mathcal{T})$
in the version of $(\mathbb{T}_c,\Delta_c)$ where maps are into
$(M,d_M)=(\ell^1,d_\ell^1)$. Moreover, taking projections $\pi_k$
yields $({T}_n,a_nd_{{T}_n},b_n\mu_{{T}_n},\pi_k\circ a_n\tilde{\psi
}_{n},\rho_{{T}_n})
\rightarrow(\mathcal{T},d_\mathcal{T},\mu_{\mathcal{T}},\pi_k\circ\tilde
{\psi},\rho_\mathcal{T})$
in the version of $(\mathbb{T}_c,\Delta_c)$ where maps are into $\mathbb
{R}^k$. Hence, by applying Lemma~\ref{measconv}, and noting the
characterisation of $\tilde{\psi}$ at (\ref{psichar}), it follows that
%
\begin{equation}
\label{c2} b_n\mu_{n,k}\circ(a_n\tilde{
\psi}_n)^{-1}\rightarrow\mu^{(k)}\circ \tilde{
\psi}^{-1}
\end{equation}
weakly as measures on $\mathbb{R}^k$, and the same conclusion also
holds in terms of weak convergence of measures on $\ell_1$. The two
conditions (\ref{c1}) and (\ref{c2}) enable us
to obtain from \cite{Croy3}, Proposition~3.1, (see also the proof of
\cite{Croy3}, Lemma~4.2) that
%
\begin{equation}
\label{2conv} \bigl(a_n\tilde{\psi}_n
\bigl(J^{n,k}_{t/a_n^2} \bigr), a_nb_n
\hat{A}^{n,k}_{t/a_n^2} \bigr)_{t\geq0}\rightarrow \bigl(
\tilde {\psi} \bigl(X^{\mathcal{T}(k),\lambda^{(k)}}_t \bigr),\hat
{A}^{(k)}_{t} \bigr)_{t\geq0}
\end{equation}
in distribution in $C(\mathbb{R}_+,\ell^1\times\mathbb{R}_+)$. Since
the functions $t\mapsto\hat{A}^{(k)}_{t}$ are almost-surely continuous
and strictly increasing (by \cite{Croy3}, Lemma~2.5), one can take an
inverse in the second coordinate and compose with the first to obtain
%
\begin{equation}
\label{c3} \bigl(a_n\tilde{\psi}_n \bigl(
\hat{X}^{n,k}_{t/a_nb_n} \bigr) \bigr)_{t\geq0}\rightarrow
\bigl(\tilde{\psi} \bigl(X^{\mathcal{T}(k),\mu
^{(k)}}_t \bigr)
\bigr)_{t\geq0}
\end{equation}
in distribution in $C(\mathbb{R}_+,\ell^1)$, for which it is helpful to
recall the expressions at (\ref{xlam}) and (\ref{hatxnk}). Now, from
our choice of $x_i^n$, one can check that, for every $i,j\geq0$,
$c_n\phi_{T_n}(b_{T_n}(\rho_{T_n},x_i^n,x^n_j))\rightarrow\phi
_{\mathcal{T}}(b_{\mathcal{T}}(\rho_{\mathcal{T}},x_i,x_j))$, where the
function $b_{T_n}$ returns the branch-point of three vertices of $T_n$.
This allows one to transfer the convergence of (\ref{c3}) into $M$, and
so obtain the result.\end{pf}

In light of Lemmas \ref{t1} and \ref{t0}, the proof of Theorem~\ref
{srwthm} is completed by the following lemma (see \cite{Bill}, Theorem~3.2, e.g.).

\begin{lem}Suppose that (\ref{tnconv}) holds. For every $\varepsilon>0$,
\[
\lim_{k\rightarrow\infty}\limsup_{n\rightarrow\infty}P^{T_n}_{\rho
_{T_n}}
\Bigl( \sup_{s\in[0,t]}d_M \bigl(c_n
\phi_{T_n} \bigl({X}^{T_n}_{s/a_nb_n} \bigr),
c_n\phi_{T_n}^{(k)} \bigl(\hat{X}^{n,k}_{s/a_nb_n}
\bigr) \bigr)> \varepsilon \Bigr)=0.
\]
\end{lem}

\begin{pf} From the definition of $\varepsilon_n$, one can obtain
\begin{eqnarray*}
&&\lim_{\delta\rightarrow0} \lim_{k\rightarrow\infty} \limsup
_{n\rightarrow\infty} \mathop{\max_{x\in T_n,y\in T_n(k):}}_{a_nd_{T_n}(x,y) \leq\delta}
c_nd_M \bigl(\phi_{T_n}(x),
\phi_{T_n}^{(k)}(y) \bigr)
\\
&&\qquad\leq\lim_{\delta\rightarrow0}\lim_{k\rightarrow\infty}\limsup
_{n\rightarrow\infty} \Bigl(\mathop{\sup_{x,y\in\mathcal{T}:}}_{d_{\mathcal{T}}(x,y)\leq
\delta+\delta_k+18\varepsilon_n}d_M
\bigl(\phi_{\mathcal{T}}(x),\phi _{\mathcal{T}}(y) \bigr)+2
\varepsilon_n \Bigr).
\end{eqnarray*}

Here, $\delta_k$ is the maximum $d_\mathcal{T}$-distance between two
adjacent vertices of $\mathcal{B}(k)$, where by saying $x,y\in\mathcal
{B}(k)$ are adjacent, we mean that $\gamma_{\mathcal{T}}(x,y)$ contains
no element of $\mathcal{B}(k)$ other than $x$ and $y$. By the
continuity of $\phi_\mathcal{T}$ and denseness of $(x_i)_{i\geq1}$,
the upper bound above is equal to 0. Hence, the lemma will follow from
%
\begin{equation}
\label{target} \lim_{k\rightarrow\infty}\limsup_{n\rightarrow\infty}P^{T_n}_{\rho_{T_n}}
\Bigl( a_n\sup_{s\in[0,t]}d_{T_n}
\bigl({X}^{n,k}_{s/a_nb_n}, \hat{X}^{n,k}_{s/a_nb_n}
\bigr)> \varepsilon \Bigr)=0,
\end{equation}
where we have applied (\ref{deltank}) to replace $X^{T_n}$ by $X^{n,k}$
in this requirement. Now, by making the change from $\alpha_n^{-1}$ to
$a_n$ and from $n^{-1}$ to $b_n$, one can follow the argument of \cite{Croy3},
Lemma~4.3, exactly to deduce that
\[
\lim_{k\rightarrow\infty}\limsup_{n\rightarrow\infty}P^{T_n}_{\rho
_{T_n}}
\Bigl( a_nb_n\sup_{s\leq t}\bigl\llvert
{A}^{n,k}_{s/a_n^2}-\hat {A}^{n,k}_{s/a_n^2}\bigr
\rrvert > \varepsilon \Bigr)=0.
\]
Note that we needed the fact that $b_n=o(a_n)$,
and also (\ref{deltank}), (\ref{2conv}) and $b_n\mu
_{T_n}(T_n)\rightarrow\mu_\mathcal{T}(\mathcal{T})$ [which follows from
our assumption at (\ref{tnconv})]. Recalling the characterisations of
$X^{n,k}$ and $\hat{X}^{n,k}$ given in (\ref{xnk}) and (\ref{hatxnk}),
respectively,
we can complete the proof of (\ref{target}) by combining the limit above
with the convergence statements of (\ref{c3}) and Lemma~\ref{t1} (cf.
\cite{Croy3}, Proposition~4.1).\end{pf}

The following measurability result will be useful when we come to
look at random walks on random trees.
Its proof is similar to that of \cite{Croy1}, Lemma~8.1(b).

\begin{lem}\label{measlemm} The map $\uT\mapsto P^\mathcal{T}_{\rho_{\mathcal
{T}}}\circ\phi_{\mathcal{T}}^{-1}$ defines a measurable function from
$\mathbb{T}_c^*$ (equipped with the subspace $\sigma$-algebra) to the
space of probability measures on $C(\mathbb{R}_+,M)$.
\end{lem}

\begin{pf} Suppose that
$\uT_n\rightarrow\uT$ in $\mathbb{T}_c^*$. A straightforward
adaptation of \cite{Miermont}, Proposition~10, then yields that if
$(x_i^n)_{i\geq1}$ is a sequence of $\mu_{\mathcal{T}_n}$-random
vertices of $\mathcal{T}_n$ and $(x_i)_{i\geq1}$ is a sequence of $\mu
_\mathcal{T}$-random vertices of $\mathcal{T}$, then for each fixed
$k\geq1$,
\[
\bigl(\mathcal{T}_n,d_{\mathcal{T}_n},\mu_{\mathcal{T}_n},
\phi_{\mathcal
{T}_n},\rho_{\mathcal{T}_n},\bigl(x_i^n
\bigr)_{i=1}^k \bigr)\rightarrow \bigl(
\mathcal{T},d_\mathcal{T},\mu_{\mathcal{T}},\phi_\mathcal{T},\rho
_\mathcal{T},(x_i)_{i=1}^k \bigr),
\]
in distribution in a version of $(\mathbb{T}_c,\Delta_c)$ where metric
spaces are marked by $k$ points [so that the supremum in (\ref
{deltacdef}) is taken over correspondences that include not only the
root pairs $(\rho_\mathcal{T},\rho_\mathcal{T}')$, but also the marked
pairs $(x_i,x_i')$, $i=1,\ldots,k$, say]. Since the latter space can be
checked to be separable in the same way as was discussed for $(\mathbb
{T}_c,\Delta_c)$ in the proof of Proposition~\ref{sepprop}, one can
apply a Skorohod representation argument to deduce that there exist
realisations of the relevant random variables such that the above
convergence occurs almost surely. As a consequence, the proof of Lemma~\ref{t0} can be applied to deduce that for fixed $k\geq1$, as
$n\rightarrow\infty$,
\[
P_{\rho_{\mathcal{T}_n}}^{\mathcal{T}_n(k),\mu_{\mathcal
{T}_n}^{(k)}}\circ\phi_{\mathcal{T}_n}^{-1}
\rightarrow P_{\rho
_{\mathcal{T}}}^{\mathcal{T}(k),\mu_\mathcal{T}^{(k)}}\circ\phi _{\mathcal{T}}^{-1}
\]
in distribution as probability measures on $C(\mathbb{R}_+,M)$, where
$P_{\rho_{\mathcal{T}}}^{\mathcal{T}(k),\mu_\mathcal{T}^{(k)}}$ is the
law of $X^{\mathcal{T}(k),\mu^{(k)}}$ started from $\rho_\mathcal{T}$,
and the objects indexed by $n$ are defined analogously to the limiting
ones. In particular, this establishes that the map from $\uT$ to the
law of $P_{\rho_{\mathcal{T}}}^{\mathcal{T}(k),\mu_\mathcal
{T}^{(k)}}\circ\phi_{\mathcal{T}}^{-1}$ is continuous, and therefore
measurable, on $\mathbb{T}_c^*$. Moreover, since $\mu_\mathcal{T}$ is
nonatomic and has full support, then we may assume that $(x_i)_{i\geq
1}$ is almost-surely dense in $\mathcal{T}$, and that all the vertices
are distinct (and not equal to the root $\rho_\mathcal{T}$). Hence, by
Lemma~\ref{t0}, it holds that
\[
P_{\rho_{\mathcal{T}}}^{\mathcal{T}(k),\mu_\mathcal{T}^{(k)}}\circ\phi _{\mathcal{T}}^{-1}
\rightarrow P_{\rho_{\mathcal{T}}}^{\mathcal
{T}}\circ\phi_{\mathcal{T}}^{-1}
\]
almost-surely as probability measures on $C(\mathbb{R}_+,M)$.
Thus, the map from $\uT$ to the law of $P_{\rho_{\mathcal{T}}}^{\mathcal
{T}}\circ\phi_{\mathcal{T}}^{-1}$ is
a limit of measurable functions, and so is also measurable on $\mathbb{T}_c^*$.
Since $P_{\rho_{\mathcal{T}}}^{\mathcal{T}}\circ\phi_{\mathcal
{T}}^{-1}$ is a function of only $\uT$ [and not the particular sequence
$(x_i)_{i\geq1}$], the result follows from a standard argument.
\end{pf}

Suppose that $\uT$ is a random element of $\mathbb{T}_c$, built on a
probability space with probability measure $\mathbf{P}$, and $\mathbf
{P}$-a.s. takes values in $\mathbb{T}_c^*$. The previous lemma tells us
that the annealed law of the process $\phi_\mathcal{T}(X^\mathcal{T})$,
where the Brownian motion $X^\mathcal{T}$ is started from the root,
that is, $\int_{\mathbb{T}_c}P_{\rho_\mathcal{T}}^\mathcal{T}\circ\phi
_\mathcal{T}^{-1}(\cdot) d{\mathbf{P}}$ [cf. (\ref{annp})], is a
well-defined probability measure on $C(\mathbb{R}_+,M)$. By a Skorohod
representation argument, we also obtain the following as an immediate
corollary of Theorem~\ref{srwthm}.

\begin{cor}\label{randomsrwcor}
Let $({T}_n,d_{{T}_n},\mu_{{T}_n},\phi_{{T}_n},\rho_{{T}_n})$, $n\geq
1$ be a random sequence,
and $(a_n)_{n\geq1}$, $(b_n)_{n\geq1}$, $(c_n)_{n\geq1}$ be null sequences
with $b_n=o(a_n)$, such that (\ref{tnconv}) holds in distribution, and
the limit $\uT$ almost-surely
takes values in $\mathbb{T}_c^*$. Then the annealed laws of the processes
\[
\bigl(c_n\phi_{T_n} \bigl(X^{T_n}_{t/a_nb_n}
\bigr) \bigr)_{t\geq0}
\]
[cf. (\ref{annlaw})] converge to the annealed law of
$\phi_\mathcal{T}(X^\mathcal{T})$, where we assume that $X^{T_n}_0=\rho
_{T_n}$ for each $n$, and also $X^\mathcal{T}_0=\rho_\mathcal{T}$.
\end{cor}

As with Lemmas \ref{measconv} and \ref{tlem}, Theorem~\ref{srwthm} and
Corollary~\ref{randomsrwcor} can be extended to the noncompact case
with an additional assumption. To begin with the deterministic case,
suppose $(T_n)_{n\geq1}$ is a deterministic sequence of locally finite
graph trees for which (\ref{tnconv}) holds in $(\mathbb{T},\Delta)$,
where $\uT$ is such that $\mu_\mathcal{T}$ is nonatomic, supported on
the leaves of $\mathcal{T}$, and satisfies (\ref{measlower}) when the
infimum is taken over $B_\mathcal{T}(\rho_\mathcal{T},R)$ for any $R$;
we denote the subset of $\mathbb{T}$ whose elements satisfy these
properties by~$\mathbb{T}^*$. (Note that for an element of $\mathbb
{T}^*$, it is possible to define Brownian motion $X^\mathcal{T}$ on
$(\mathcal{T},d_{\mathcal{T}},\mu_\mathcal{T})$ by the procedure of
\cite{AEW}, as described above the statement of Theorem~\ref{srwthm}.)
Moreover, assume that, for $t>0$,
\[
\lim_{R\rightarrow\infty}\limsup_{n\rightarrow\infty}{P}_{\rho
_{T_n}}^{T_n}
\bigl(\tau \bigl(X^{T_n},B_{T_n}\bigl(\rho
_{T_{n}},a_n^{-1}R\bigr) \bigr)\leq
t/a_nb_n \bigr)=0,
\]
where, here and in the following, $\tau(X,A):=\inf\{t\geq0: X_t\notin A\}$ is the exit time of a process $X$ from a set $A$. It is then
the case that (\ref{srwconv}) holds. Similarly to Remarks~\ref{measrem}
and \ref{otherrem}, we do not include the full details of this
argument, but instead restrict our presentation to describing the
probabilistic version of the argument needed to handle the simple
random walk on $\mathcal{U}$. In particular, given Theorems \ref{main1}
and \ref{main12}, the key additional ingredient for this is the
following, where we recall $\mathbb{P}$ is the annealed law of the
simple random walk on $\mathcal{U}$, as introduced at (\ref{annlaw}).

\begin{propn}\label{exitprop} For $t>0$,
\[
\lim_{R\rightarrow\infty}\limsup_{\delta\rightarrow0}\mathbb{P} \bigl(
\tau \bigl(X^\mathcal{U},B_\mathcal{U}\bigl(0,\delta^{-\kappa}
R\bigr) \bigr)\leq t\delta^{-\kappa d_w} \bigr)=0.
\]
\end{propn}

\begin{pf} Given the volume growth and resistance estimates of \cite{BM11},
Proposition~4.2, we can apply an identical argument to that
used to prove the corresponding limit in \cite{CroyKum}, Theorem~1.1.
\end{pf}

We are now ready to prove the main simple random walk convergence
result for the UST. In the proof, we let $(\mathbf{P}_{\delta
_n})_{n\geq1}$ be a convergent sequence with limit $\tilde{\mathbf
{P}}$, and suppose that $\uT$ is a random variable with law $\tilde
{\mathbf{P}}$. Unless otherwise noted, we take $(M,d_M)$ to be $\mathbb
{R}^2$ equipped with the Euclidean distance.

\begin{pf*}{Proof of Theorem~\ref{main2}\normalfont{(a), (b)}} As in the proof of Lemma~\ref{measlem},
the separability of $(\mathbb{T},\Delta)$ allows us to
find realisations of $\uU_{\delta_n}$, $n\geq1$, and $\uT$ built on a
common probability space with probability measure $\mathbf{P}^*$ such
that $\uU_{\delta_n}\rightarrow\uT$ holds in $(\mathbb{T},\Delta)$,
$\mathbf{P}^*$-a.s. This yields the existence of a subsequence
$(n_i)_{i\geq1}$ and divergent sequence $(r_j)_{j\geq1}$ such that,
for every $r_j$, $\mathbf{P}^*$-a.s., we have convergence in $(\mathbb
{T}_c,\Delta_c)$ of the radius $r_j$ restrictions along the subsequence
$(n_i)_{i\geq1}$, that is, $\uU_{\delta_{n_i}}^{(r_j)}$ converges to
$\uT^{(r_j)}$. Moreover, since the Lebesgue measure of those $r\geq0$
for which $\mu_\mathcal{T}(\partial\mathcal{T}^{(r)})>0$ is zero,
$\tilde{\mathbf{P}}$-a.s., we may further assume that $\mu_\mathcal
{T}(\partial\mathcal{T}^{(r_j)})=0$ for every $r_j$, $\tilde{\mathbf
{P}}$-a.s. Noting that the map $\uT\mapsto\uT^{(r)}$ is continuous at
those elements of $\mathbb{T}$ satisfying $\mu_{\mathcal{T}}(\partial
\mathcal{T}^{(r)})=0$ [cf. the proof of Proposition~\ref{deltametric},
and equation (\ref{extra}) in particular], this final condition ensures
that, for every $r_j$, $\uT^{(r_j)}$ is $\uT$-measurable. Now, from
Theorem~\ref{main12}(b)(ii), it is the case that $\uT^{(r)}$ takes
values in $\mathbb{T}^*_c$, $\tilde{\mathbf{P}}$-a.s. [Observe that if
$r'\leq r/2$ and $x\in\mathcal{T}^{(r)}$, then there exists an $x'\in
\mathcal{T}^{(r)}$ such that $B_{\mathcal{T}^{(r)}}(x,r')\supseteq
B_\mathcal{T}(x',r'/2)$, and so there is no problem in deducing the
lower volume estimate of (\ref{measlower}) for the ball $\mathcal
{T}^{(r)}$ from the lower volume estimate from $\mathcal{T}$. The
further two properties---that $\mu^{(r)}_\mathcal{T}$ is nonatomic
and supported on the leaves of $\mathcal{T}^{(r)}$---are immediate
from the full $\mathcal{T}$ statements.] Hence, by Corollary~\ref
{randomsrwcor} (with $a_{n_i}=\delta_{n_i}^{\kappa}$, $b_{n_i}=\delta
_{n_i}^2$, $c_{n_i}=\delta_{n_i}$), it follows that, for every $r_j$,
the annealed law of
\[
\bigl(\delta_{n_i}X^{{B}_\mathcal{U}(0,\delta_{n_i}^{-\kappa
}r_j)}_{\delta_{n_i}^{-\kappa d_w}t}
\bigr)_{t\geq0},
\]
where $X^{{B}_\mathcal{U}(0,\delta_{n_i}^{-\kappa}r_j)}$ is the simple
random walk on ${B}_\mathcal{U}(0,\delta_{n_i}^{-\kappa}r_j)$,
converges as $i\rightarrow\infty$ to the annealed law of $(\phi
^{(r_j)}_\mathcal{T}(X^{\mathcal{T}^{(r_j)}}_t))_{t\geq0}$, where
$X^{\mathcal{T}^{(r_j)}}$ is Brownian motion on the measured real tree
$(\mathcal{T}^{(r_j)},d^{(r_j)}_\mathcal{T},\mu^{(r_j)}_\mathcal{T})$,
and we assume processes are started from the roots of the relevant
trees. Given the measurability of $\uT^{(r_j)}$ described above, this
limit law can be expressed as
%
\begin{equation}
\label{annt} \int P^{\mathcal{T}^{(r_j)}}_{\rho_\mathcal{T}}\circ\phi_\mathcal
{T}^{-1} (\cdot ) \,d\tilde{\mathbf{P}}.
\end{equation}
Furthermore, since the limit does not depend on the subsequence, we
obtain annealed distributional convergence of the rescaled discrete
processes along the full sequence $(\delta_n)_{n\geq1}$.

Next, suppose $X^\mathcal{U}$ and $X^{{B}_\mathcal{U}(0,\delta
_{n_i}^{-\kappa}r_j)}$ are coupled so that their sample paths agree up
to $\tau(X^\mathcal{U}, {B}_\mathcal{U}(0,\delta_{n_i}^{-\kappa}r_j))$
[e.g., by taking $X^{{B}_\mathcal{U}(0,\delta_{n_i}^{-\kappa}r_j)}$ to
be $X^\mathcal{U}$ observed on ${{B}_\mathcal{U}(0,\delta_{n_i}^{-\kappa
}r_j)}$]. It then holds that, for $t<\infty$ and $\varepsilon>0$,
\[
\mathbb{P} \Bigl(\sup_{s\in[0,t]} \delta_n\bigl\llvert
X^\mathcal{U}_{\delta_{n}^{-\kappa d_w}s}-X^{{B}_\mathcal
{U}(0,\delta_n^{-\kappa}r_j)}_{\delta_{n}^{-\kappa d_w}s}\bigr
\rrvert >\varepsilon \Bigr) \leq\mathbb{P} \bigl( \tau \bigl(X^\mathcal{U},{B}_\mathcal{U}
\bigl(0,\delta_n^{-\kappa} r_j\bigr) \bigr)\leq t
\delta^{-\kappa d_w} \bigr).
\]
Hence, we deduce from Proposition~\ref{exitprop} that the left-hand
side converges
to 0 as $n\rightarrow\infty$ and then $j\rightarrow\infty$.

Finally, we need to take care of the situation when $r_j\rightarrow
\infty$ for the continuous trees. We have already established that
$P^{\mathcal{T}^{(r_j)}}_{\rho_\mathcal{T}}\circ\phi^{-1}_{\mathcal
{T}}$ is $\uT$-measurable for every $r_j$. To show that this is also
the case for $P^{\mathcal{T}}_{\rho_\mathcal{T}}\circ\phi^{-1}_{\mathcal
{T}}$, it will suffice to check that $P^{\mathcal{T}^{(r)}}_{\rho
_\mathcal{T}}\rightarrow P^{\mathcal{T}}_{\rho_\mathcal{T}}$ as
$r\rightarrow\infty$, $\tilde{\mathbf{P}}$-a.s. This will follow if we
can check that
%
\begin{equation}
\label{toch} P^{\mathcal{T}}_{\rho_\mathcal{T}} \Bigl(\lim_{r\rightarrow\infty}
\tau \bigl(X^\mathcal{T},B_\mathcal{T}(\rho_\mathcal{T},r)
\bigr)=\infty \Bigr)=1,
\end{equation}
for $\tilde{\mathbf{P}}$-a.e. realisation of $\mathcal{T}$. To this
end, first note that
%
\begin{eqnarray}
\label{trip} &&P^{\mathcal{T}}_{\rho_\mathcal{T}} \Bigl(\lim_{r\rightarrow\infty}
\tau \bigl(X^\mathcal{T},B_\mathcal{T}(\rho_\mathcal{T},r)
\bigr)<\infty \Bigr)
\nonumber
\\[-8pt]
\\[-8pt]
\nonumber
&&\qquad= \lim_{t\rightarrow\infty}\lim_{r\rightarrow\infty}
\lim_{j\rightarrow
\infty}P^{\mathcal{T}^{(r_j)}}_{\rho_\mathcal{T}} \bigl(\tau
\bigl(X^{\mathcal{T}^{(r_j)}},B_\mathcal{T}(\rho_\mathcal{T},r) \bigr)\leq
t \bigr),
\end{eqnarray}
where we note that the laws of $X^\mathcal{T}$ and $X^{\mathcal
{T}^{(r_j)}}$ agree up to the exit time of $B_\mathcal{T}(\rho_\mathcal
{T},r)$ whenever $r_j>r$. Now, suppose that $\uU_{\delta_n}$, $n\geq
1$, and $\uT$ are coupled as in the first part of the proof. It is not
difficult to check from the definition of $\Delta_c$ that the
convergence in $(\mathbb{T}_c,\Delta_c)$ of the radius $r_j$
restrictions still holds $\mathbf{P}^*$-a.s. if $\delta_n\phi_\mathcal
{U}$ is replaced by $\delta_n^\kappa\tilde{\phi}_\mathcal{U}$ and $\phi
_\mathcal{T}$ is replaced by $\tilde{\phi}_\mathcal{T}$, where $\tilde
{\phi}_\mathcal{U}(x):=(d_\mathcal{U}(0,x),0)$ and $\tilde{\phi
}_\mathcal{T}(x):=(d_\mathcal{T}(0,x),0)$, respectively. Thus, an
application of Theorem~\ref{srwthm} (with $M=\mathbb{R}$) and a
subsequence argument yields that, for every $r_j$, $\mathbf{P}^*$-a.s.,
$(\delta_n^\kappa d_\mathcal{U}(0,X^{{B}_\mathcal{U}(0,\delta_n^{-\kappa
} r_j)}_{\delta_{n}^{-\kappa d_w}t}))_{t\geq0}$ converges to
$(d_\mathcal{T}(\rho_\mathcal{T},X^{\mathcal{T}^{(r_j)}}_t))_{t\geq0}$
in distribution in $C(\mathbb{R}_+,\mathbb{R})$ as $n\rightarrow\infty
$. Consequently,
\begin{eqnarray*}
&&\int P^{\mathcal{T}^{(r_j)}}_{\rho_\mathcal{T}} \bigl(\tau \bigl(X^{\mathcal{T}^{(r_j)}},B_\mathcal{T}(
\rho_\mathcal{T},r) \bigr)\leq t \bigr) \,d\tilde{\mathbf{P}}
\\
&&\qquad=\int P^{\mathcal{T}^{(r_j)}}_{\rho_\mathcal{T}} \Bigl( \sup_{s\in[0,t]}d_\mathcal{T}
\bigl(\rho_\mathcal{T},X^{\mathcal
{T}^{(r_j)}}_s \bigr)\geq r
\Bigr) \,d\tilde{\mathbf{P}}
\\
&&\qquad\leq\int\liminf_{n\rightarrow\infty}P^{{B}_\mathcal{U}(0,\delta
_{n_i}^{-\kappa}r_j)}_{0}
\Bigl(\sup_{s\in[0,t]}d_\mathcal{U} \bigl(0,
X^{{B}_\mathcal{U}(0,\delta_{n}^{-\kappa}r_j)}_{\delta_{n}^{-\kappa
d_w}s} \bigr)\geq\delta_n^{-\kappa}r/2
\Bigr)\,d{\mathbf{P}}
\\
&&\qquad\leq \liminf_{n\rightarrow\infty}\mathbb{P} \bigl(\tau
\bigl(X^{{B}_\mathcal{U}(0,\delta_n^{-\kappa}r_j)},{B}_\mathcal{U}\bigl(0,\delta _n^{-\kappa}r
\bigr) \bigr)\leq t\delta_n^{-\kappa d_w} \bigr),
\end{eqnarray*}
where we note that the necessary measurability of the law of $d_\mathcal
{T}(\rho_\mathcal{T},X^{\mathcal{T}^{(r_j)}})$ can be checked similarly
to the measurability of $P^{\mathcal{T}^{(r_j)}}_{\rho_\mathcal{T}}\circ
\phi^{-1}_{\mathcal{T}}$. Taking limits as $j\rightarrow\infty$,
$r\rightarrow\infty$ and then $t\rightarrow\infty$, we obtain from
another application of Proposition~\ref{exitprop} that the upper bound
above converges to zero. Thus, the dominated convergence theorem yields
that the $\tilde{\mathbf{P}}$-expectation of the left-hand side of (\ref
{trip}) is equal to zero, and so we have established (\ref{toch}), as
desired. In summary, we have now shown that $P^{\mathcal{T}}_{\rho
_\mathcal{T}}\circ\phi^{-1}_{\mathcal{T}}$ is $\uT$-measurable, and so
$\tilde{\mathbb{P}}$---the annealed law of $\phi_\mathcal{T}(X^\mathcal
{T})$, as introduced at (\ref{annp})---is well-defined. This
establishes part (a) of the theorem. Moreover, since $P^{\mathcal
{T}^{(r)}}_{\rho_\mathcal{T}}\rightarrow P^{\mathcal{T}}_{\rho_\mathcal
{T}}$, $\tilde{\mathbf{P}}$-a.s., the continuous mapping theorem and
the dominated convergence theorem yield that the annealed law at (\ref
{annt}) converges as $j\rightarrow\infty$ to $\tilde{\mathbb{P}}$.
Together with the conclusions of the previous two paragraphs, this
completes the proof of part (b) (see \cite{Bill}, Theorem~3.2, e.g.).
\end{pf*}

\begin{rem}
In \cite{Croy2}, Theorem~8.1, it was shown that the convergence
of rescaled graph tree ``tours'' (i.e., functions encoding trees and
embeddings into~$\mathbb{R}^d$) implies convergence of spatial trees.
It only requires a simple extension of the proof of that result to add
the measure, and thereby deduce that convergence of tours also implies
convergence in the topology of $(\mathbb{T}_c,\Delta_c)$ (with $\mathbb
{R}^2$ replaced by~$\mathbb{R}^d$). Consequently, (the $\mathbb{R}^d$
version of) Theorem~\ref{srwthm} provides an extension of the random
walk convergence result of \cite{Croy2}, Theorem~8.1, and can also be
used to deduce a scaling limit for the simple random walks on critical
branching random walks satisfying the various assumptions on the
offspring and step distribution detailed in \cite{Croy2}, Section~10.
\end{rem}

\begin{rem}
In \cite{Croy3}, Theorem~1.1 (see also \cite{Croy1}, Theorem~1.1), a
scaling limit was established for simple random walks on ordered graph
trees. The proofs of these results used the convenience of encoding
ordered graph trees by functions to simplify various technical details.
Theorem~\ref{srwthm} provides the additional framework needed to handle
unordered graph trees. Indeed, suppose that we have a sequence of
finite rooted graph trees such that
%
\begin{equation}
\label{unor} ({T}_n,a_nd_{{T}_n},b_n
\mu_{{T}_n},\rho_{{T}_n} )\rightarrow (\mathcal{T},d_\mathcal{T},
\mu_{\mathcal{T}},\rho_\mathcal{T} ),
\end{equation}
for some null sequences $(a_n)_{n\geq1}$ and $(b_n)_{n\geq1}$ with
$b_n=o(a_n)$ in the Gromov--Hausdorff--Prohorov topology of \cite{ADH}
for rooted, compact metric spaces\break  equipped with Borel measures.
Moreover, assume that the limit satisfies the additional properties on
the measure that hold for elements of $\mathbb{T}_c^*$. It is then the
case that, by applying the procedure described in the proof of Lemma~\ref{t0}, one can find isometric embeddings $\tilde{\psi}:\mathcal
{T}\rightarrow\ell^1$ and $\tilde{\psi}_n:{T}_n\rightarrow\ell^1$,
$n\geq1$, such that
\[
\bigl(a_n\tilde{\psi}_n \bigl(X^{T_n}_{t/a_nb_n}
\bigr) \bigr)_{t\geq
0}\rightarrow \bigl(\tilde{\psi} \bigl(X^{\mathcal{T}}_{t}
\bigr) \bigr)_{t\geq0},
\]
in distribution in $C(\mathbb{R}_+,\ell^1)$, where $X^{T_n}$ is the
simple random walk on $T_n$ started from $\rho_{T_n}$, and $X^{\mathcal
{T}}$ is the Brownian motion on $(\mathcal{T},d_\mathcal{T},\mu
_{\mathcal{T}})$ started from $\rho_\mathcal{T}$. We note that a
similar, but slightly stronger result was recently established
independently in \cite{ALW}---the most important difference being that
the argument of the latter work required a weaker lower volume growth condition.

The above result can also be extended to the random case by embedding
trees into $\ell^1$ in a canonical random way---specifically, by
choosing the vertices $(x_i)_{i\geq1}\in\mathcal{T}$ to be a $\mu
_\mathcal{T}$-random sequence. The only additional complications come
from some measurability issues, but these can be resolved using similar
ideas to those applied in the proof of Lemma~\ref{measlemm}. To
summarise the conclusion, suppose that $(T_n)_{n\geq1}$ is a sequence
of finite rooted (unordered) graph trees for which (\ref{unor}) holds
in distribution for some random measured compact real tree $(\mathcal
{T},d_\mathcal{T},\mu_\mathcal{T},\rho_\mathcal{T})$. Moreover, assume
that $\mu_\mathcal{T}$ is a nonatomic probability measure, supported
on the leaves of $\mathcal{T}$ and satisfies (\ref{measlower}),
almost-surely. It can then be checked that the annealed laws of
$(a_n\tilde{\psi}_n(X^{T_n}_{t/a_nb_n}))_{t\geq0}$, where $\tilde{\psi
}_n$ is the canonical random isometric embedding of $T_n$ into $\ell
^1$, converge to the annealed law of Brownian motion on $\mathcal{T}$
randomly isometrically embedded into $\ell^1$. For example, taking
$a_n=n^{-1/2}$, $b_n=n^{-1}$, this result applies to the model of
uniformly random unordered trees with $n$ vertices, in which each
vertex has at most $m\in\{2,3,\ldots,\infty\}$ children, as studied in
\cite{HM}. Moreover, applying the argument of \cite{CHK}, Section~5.2,
(where critical Galton--Watson trees were considered), a corollary of
this is that the mixing times of the simple random walks on the random
discrete trees, when rescaled by $n^{3/2}$, converge in distribution to
(a constant multiple of) the mixing time of the limiting diffusion.
\end{rem}

\section{Heat kernel bounds for the scaling limit}\label{hklimsec}

As in Section~\ref{propsec} and the proof of Theorem~\ref{main2}(a), (b),
let $(\mathbf{P}_{\delta_n})_{n\geq1}$ be a convergent sequence with limit
$\tilde{\mathbf{P}}$, and suppose that
$(\mathcal{T},d_\mathcal{T},\mu_\mathcal{T},\phi_\mathcal{T},\rho
_\mathcal{T})$ is a random variable with law $\tilde{\mathbf{P}}$.
It follows from~\cite{AEW}, Remark~3.1 and \cite{FOT}, Theorem~1.5.2,
that the Dirichlet form $(\mathcal{E}_\mathcal{T},\mathcal{F}_\mathcal{T})$
given in Section~\ref{srwsec} is the same as that of \cite{Kigami}, Section~5.
In particular, this is the form
associated with the natural ``resistance form'' on $(\mathcal
{T},d_\mathcal{T})$, and so we can apply \cite{Kigami2}, Theorem~10.4,
to deduce the existence of a jointly continuous transition density
$(p_t^\mathcal{T}(x,y))_{x,y\in\mathcal{T},t>0}$ for the process
$X^\mathcal{T}$.

Let $R_\mathcal{T}$ be the resistance associated with
$(\mathcal{E}_\mathcal{T},\mathcal{F}_\mathcal{T})$, defined
by setting, for two disjoint subsets $A,B\subseteq\mathcal{T}$,
\[
R_\mathcal{T}(A,B)^{-1}:=\inf\bigl\{\mathcal{E}_\mathcal{T}(f,f):
f\in \mathcal{F}_\mathcal{T}, f|_A=0,f|_B=1
\bigr\}.
\]
Since $\mathcal{T}$ is complete and, by Theorem~\ref{main12}(a)(ii),
has a single end at infinity, we deduce from \cite{AEW}, Theorem~4, that
$X^\mathcal{T}$ is recurrent. As a consequence, combining \cite{FOT}, Theorem~1.6.3 and \cite{AEW}, Proposition~3.5, yields that the resistance
between two points corresponds to (a multiple of) the distance between
them, that is, $R_\mathcal{T}(\{x\},\{y\})=2d_\mathcal{T}(x,y)$ for all
$x,y\in\mathcal{T}$, $x\neq y$ (see \cite{AEW}, e.g.).
Hence, we can use \cite{Croy-1}, Theorem~3, to obtain estimates for
$p_t^\mathcal{T}(x,y)$
from the volume estimates of Theorem~\ref{main12}(b)(ii).

\begin{pf*}{Proof of Theorem~\ref{main2}\normalfont{(c)}} We have already discussed
the claims about the existence and joint continuity of the heat kernel,
and the recurrence of $X^\mathcal{T}$. So, we will simply present here
a few key points that are needed to apply the heat kernel estimates of
\cite{Croy-1}, Theorem~3. As already noted, the resistance metric
coincides with (a multiple of) the tree metric $d_\mathcal{T}$ in our
setting, and so we can replace $R_\mathcal{T}$ in the conclusion of
\cite{Croy-1} by $d_\mathcal{T}$. Moreover, the fact that $(\mathcal
{T},d_\mathcal{T})$ is a real tree automatically means the ``chaining
condition'' of \cite{Croy-1} is satisfied, that is, there exists a
constant $c_1$ such that for all $x,y\in\mathcal{T}$ and all $n\in
\mathbb{N}$, there exist $x_0=x,x_1,\ldots,x_N=y$ such that $d_\mathcal
{T}(x_{i-1}, x_i) \leq c_1{d_\mathcal{T}(x,y)}/{n}$, $\forall i=1,\ldots,
n$. (Clearly, we can take $c_1=1$ and equality in the latter
statement.) Finally, note that in \cite{Croy-1} volume estimates were
assumed to hold uniformly over the entire space, but Theorem~\ref
{main12}(b)(ii) only gives uniformity over balls of finite radius.
However, it is straightforward to check that the arguments of \cite
{Croy-1} are enough to give the stated heat kernel estimates.
\end{pf*}

For the remaining heat kernel estimates, we derive the following tail
bound for the resistance from the root to the radius of a ball will be useful.

\begin{lem}There exist constants $c_1,c_2,\theta\in(0,\infty)$ such that for
all $r>0$, $\lambda\geq1$,
%
\begin{equation}
\label{bthird} \tilde{\mathbf{P}} \bigl(R_{\mathcal{T}} \bigl(
\rho_\mathcal{T},B_\mathcal {T}(\rho_\mathcal{T},r)^c
\bigr)\leq\lambda^{-1} r \bigr)\leq c_1e^{-c_2\lambda^\theta}.
\end{equation}
\end{lem}

\begin{pf} As we have done several times previously, we will apply a
coupling argument, and start by supposing that we have a realisation of
random variables such that
$\uU_{\delta_n}\rightarrow\uT$ holds almost-surely along the sequence
$(\delta_n)_{n\geq1}$.
Let $N_{\mathcal{T}}(r,r/\lambda)$ be the minimum number of $d_\mathcal
{T}$-balls of
radius $r/\lambda$ required to cover $B_{\mathcal{T}}(\rho_\mathcal{T},r)$.
From the definition of $(\mathbb{T}, \Delta)$, it is elementary to
check that if
$N_{\mathcal{T}}(r,r/\lambda)\geq N_0$, then so is
$N_{\mathcal{U}}(4\delta_n^{-\kappa}r,\delta_n^{-\kappa}r/8\lambda)$
for large $n$.
It follows that
%
\begin{equation}
\label{hhh} \tilde{\mathbf{P}} \bigl(N_{\mathcal{T}}(r,r/\lambda)\geq c
\lambda ^{5} \bigr) \leq\limsup_{n\rightarrow\infty} {\mathbf{P}}
\bigl(N_{\mathcal{U}}\bigl(4\delta_n^{-\kappa}r,
\delta_n^{-\kappa
}r/8\lambda\bigr)\geq c\lambda^{5}
\bigr),
\end{equation}
which by Remark~\ref{R:covers} is, for a suitable choice of $c$,
bounded above by
$c_1e^{-c_2\lambda^{1/80}}$. Now, the proof of \cite{K04}, Lemma~4.1,
gives that
$R_{\mathcal{T}}(\rho_\mathcal{T},B_\mathcal{T}(\rho_\mathcal
{T},r)^c)\geq r/8N_{\mathcal{T}}(r, r/4)$,
and so, for any $\lambda>4^{5}c$,
\[
\tilde{\mathbf{P}} \bigl(R_{\mathcal{T}} \bigl(\rho_\mathcal{T},B_\mathcal
{T}(\rho_\mathcal{T},r)^c \bigr)\leq r/\lambda \bigr) \leq
\tilde{\mathbf{P}} \bigl(N_{\mathcal{T}}\bigl(r,c^{1/5}r/
\lambda^{1/5}\bigr)\geq 8\lambda \bigr).
\]
Combining this bound with (\ref{hhh}) yields the result.
\end{pf}

Given (\ref{bfirst}), (\ref{bsecond}) and (\ref{bthird}), the next two
results can be proved in exactly the same way as the corresponding
parts of \cite{Croy0}, Theorem~1.6 and \cite{Croy0}, Proposition~1.7,
modulo a different volume growth exponent.

\begin{thmm}\textup{(a)} $\tilde{\mathbf{P}}$-a.s., there exists a random $t_0(\sT) \in
(0,\infty)$ and deterministic constants $c_1,c_2,\theta_1,\theta_2\in
(0,\infty)$ such that
\[
c_1t^{-8/13}\bigl(\log\log t^{-1}
\bigr)^{-\theta_1} \leq p_t^\mathcal{T}(
\rho_\mathcal{T},\rho_\mathcal{T})\leq c_2t^{-8/13}
\bigl(\log\log t^{-1}\bigr)^{\theta_2},
\]
for all $t\in(0,t_0)$.

\textup{(b)} There exist constants $c_1,c_2\in(0,\infty)$ such that
\[
c_1t^{-8/13} \leq\tilde{\mathbf{E}} p_t^\mathcal{T}(
\rho_\mathcal{T},\rho_\mathcal {T})\leq c_2t^{-8/13},
\]
for all $t\in(0,1)$.
\end{thmm}

\section*{Acknowledgement}
The authors thank Sunil Chhita, who wrote the code used to produce
Figure~\ref{srwfig}.






\printaddresses

\begin{thebibliography}{49}


\bibitem{ADH}
\begin{barticle}[mr]
\bauthor{\bsnm{Abraham},~\bfnm{Romain}\binits{R.}},
\bauthor{\bsnm{Delmas},~\bfnm{Jean-Fran{\c{c}}ois}\binits{J.-F.}} \AND
\bauthor{\bsnm{Hoscheit},~\bfnm{Patrick}\binits{P.}}
(\byear{2013}).
\btitle{A note on the {G}romov--{H}ausdorff--{P}rokhorov distance between (locally) compact metric measure spaces}.
\bjournal{Electron. J. Probab.}
\bvolume{18}
\bpages{1--21}.
\bid{doi={10.1214/EJP.v18-2116}, issn={1083-6489}, mr={3035742}}
\bptnote{check pages}%
\end{barticle}
%

\bptok{imsref}%
\endbibitem

\bibitem{ABNW}
\begin{barticle}[mr]
\bauthor{\bsnm{Aizenman},~\bfnm{Michael}\binits{M.}},
\bauthor{\bsnm{Burchard},~\bfnm{Almut}\binits{A.}},
\bauthor{\bsnm{Newman},~\bfnm{Charles~M.}\binits{C.~M.}} \AND
\bauthor{\bsnm{Wilson},~\bfnm{David~B.}\binits{D.~B.}}
(\byear{1999}).
\btitle{Scaling limits for minimal and random spanning trees in two dimensions}.
\bjournal{Random Structures Algorithms}
\bvolume{15}
\bpages{319--367}.
\bid{doi={10.1002/(SICI)1098-2418(199910/12)15:3/4<319::AID-RSA8>3.3.CO;2-7}, issn={1042-9832}, mr={1716768}}
\end{barticle}
%

\bptok{imsref}%
\endbibitem

\bibitem{AKM}
\begin{barticle}[mr]
\bauthor{\bsnm{Alberts},~\bfnm{Tom}\binits{T.}},
\bauthor{\bsnm{Kozdron},~\bfnm{Michael~J.}\binits{M.~J.}} \AND
\bauthor{\bsnm{Masson},~\bfnm{Robert}\binits{R.}}
(\byear{2013}).
\btitle{Some partial results on the convergence of loop-erased random walk to {$\mathrm{SLE}(2)$} in the natural parametrization}.
\bjournal{J. Stat. Phys.}
\bvolume{153}
\bpages{119--141}.
\bid{doi={10.1007/s10955-013-0816-7}, issn={0022-4715}, mr={3100817}}
\end{barticle}
%

\bptok{imsref}%
\endbibitem

\bibitem{Aldous3}
\begin{barticle}[mr]
\bauthor{\bsnm{Aldous},~\bfnm{David}\binits{D.}}
(\byear{1993}).
\btitle{The continuum random tree. {III}}.
\bjournal{Ann. Probab.}
\bvolume{21}
\bpages{248--289}.
\bid{issn={0091-1798}, mr={1207226}}
\end{barticle}
%

\bptok{imsref}%
\endbibitem

\bibitem{Aldous0}
\begin{barticle}[mr]
\bauthor{\bsnm{Aldous},~\bfnm{David~J.}\binits{D.~J.}}
(\byear{1990}).
\btitle{The random walk construction of uniform spanning trees and uniform labelled trees}.
\bjournal{SIAM J. Discrete Math.}
\bvolume{3}
\bpages{450--465}.
\bid{doi={10.1137/0403039}, issn={0895-4801}, mr={1069105}}
\bptnote{check volume}%
\end{barticle}
%

\bptok{imsref}%
\endbibitem

\bibitem{AEW}
\begin{barticle}[mr]
\bauthor{\bsnm{Athreya},~\bfnm{Siva}\binits{S.}},
\bauthor{\bsnm{Eckhoff},~\bfnm{Michael}\binits{M.}} \AND
\bauthor{\bsnm{Winter},~\bfnm{Anita}\binits{A.}}
(\byear{2013}).
\btitle{Brownian motion on {$\mathbb{R}$}-trees}.
\bjournal{Trans. Amer. Math. Soc.}
\bvolume{365}
\bpages{3115--3150}.
\bid{doi={10.1090/S0002-9947-2012-05752-7}, issn={0002-9947}, mr={3034461}}
\end{barticle}
%

\bptok{imsref}%
\endbibitem

\bibitem{ALW}
\begin{bmisc}[auto:parserefs-M02]
\bauthor{\bsnm{Athreya},~\bfnm{S.}\binits{S.}},
\bauthor{\bsnm{L\"ohr},~\bfnm{W.}\binits{W.}} \AND
\bauthor{\bsnm{Winter},~\bfnm{A.}\binits{A.}}
(\byear{2014}).
\bhowpublished{Invariance principle for variable speed random walks on trees.
Preprint. Available at \arxivurl{arXiv:1404.6290}}.
\end{bmisc}
%

\bptok{imsref}%
\endbibitem

\bibitem{barLRW}
\begin{barticle}[mr]
\bauthor{\bsnm{Barlow},~\bfnm{M.~T.}\binits{M.~T.}}
(\byear{2016}).
\btitle{Loop erased walks and uniform spanning trees}.
\bjournal{MSJ Memoirs}
\bvolume{34}
\bpages{1--32}.
\bid{mr={3525847}}
\end{barticle}
%

\bptok{imsref}%
\endbibitem

\bibitem{BJKS}
\begin{barticle}[mr]
\bauthor{\bsnm{Barlow},~\bfnm{Martin~T.}\binits{M.~T.}},
\bauthor{\bsnm{J{\'a}rai},~\bfnm{Antal~A.}\binits{A.~A.}},
\bauthor{\bsnm{Kumagai},~\bfnm{Takashi}\binits{T.}} \AND
\bauthor{\bsnm{Slade},~\bfnm{Gordon}\binits{G.}}
(\byear{2008}).
\btitle{Random walk on the incipient infinite cluster for oriented percolation in high dimensions}.
\bjournal{Comm. Math. Phys.}
\bvolume{278}
\bpages{385--431}.
\bid{doi={10.1007/s00220-007-0410-4}, issn={0010-3616}, mr={2372764}}
\end{barticle}
%

\bptok{imsref}%
\endbibitem

\bibitem{BM10}
\begin{barticle}[mr]
\bauthor{\bsnm{Barlow},~\bfnm{Martin~T.}\binits{M.~T.}} \AND
\bauthor{\bsnm{Masson},~\bfnm{Robert}\binits{R.}}
(\byear{2010}).
\btitle{Exponential tail bounds for loop-erased random walk in two dimensions}.
\bjournal{Ann. Probab.}
\bvolume{38}
\bpages{2379--2417}.
\bid{doi={10.1214/10-AOP539}, issn={0091-1798}, mr={2683633}}
\end{barticle}
%

\bptok{imsref}%
\endbibitem

\bibitem{BM11}
\begin{barticle}[mr]
\bauthor{\bsnm{Barlow},~\bfnm{Martin~T.}\binits{M.~T.}} \AND
\bauthor{\bsnm{Masson},~\bfnm{Robert}\binits{R.}}
(\byear{2011}).
\btitle{Spectral dimension and random walks on the two dimensional uniform spanning tree}.
\bjournal{Comm. Math. Phys.}
\bvolume{305}
\bpages{23--57}.
\bid{doi={10.1007/s00220-011-1251-8}, issn={0010-3616}, mr={2802298}}
\end{barticle}
%

\bptok{imsref}%
\endbibitem

\bibitem{BLPS}
\begin{barticle}[mr]
\bauthor{\bsnm{Benjamini},~\bfnm{Itai}\binits{I.}},
\bauthor{\bsnm{Lyons},~\bfnm{Russell}\binits{R.}},
\bauthor{\bsnm{Peres},~\bfnm{Yuval}\binits{Y.}} \AND
\bauthor{\bsnm{Schramm},~\bfnm{Oded}\binits{O.}}
(\byear{2001}).
\btitle{Uniform spanning forests}.
\bjournal{Ann. Probab.}
\bvolume{29}
\bpages{1--65}.
\bid{doi={10.1214/aop/1008956321}, issn={0091-1798}, mr={1825141}}
\end{barticle}
%

\bptok{imsref}%
\endbibitem

\bibitem{Bill}
\begin{bbook}[mr]
\bauthor{\bsnm{Billingsley},~\bfnm{Patrick}\binits{P.}}
(\byear{1999}).
\btitle{Convergence of Probability Measures},
\bedition{2nd} ed.
\bpublisher{Wiley},
\blocation{New York}.
\bid{doi={10.1002/9780470316962}, mr={1700749}}
\end{bbook}
%

\bptok{imsref}%
\endbibitem

\bibitem{Broder}
\begin{bincollection}[auto:parserefs-M02]
\bauthor{\bsnm{Broder},~\bfnm{A.}\binits{A.}}
(\byear{1989}).
\btitle{Generating random spanning trees}.
In \bbooktitle{Proceedings of the 30th Annual Symposium on Foundations of Computer Science}
\bpages{442--447}.
\bpublisher{IEEE Computer Society},
\blocation{Washington, DC}.
\end{bincollection}
%

\bptok{imsref}%
\endbibitem

\bibitem{BBI}
\begin{bbook}[mr]
\bauthor{\bsnm{Burago},~\bfnm{Dmitri}\binits{D.}},
\bauthor{\bsnm{Burago},~\bfnm{Yuri}\binits{Y.}} \AND
\bauthor{\bsnm{Ivanov},~\bfnm{Sergei}\binits{S.}}
(\byear{2001}).
\btitle{A Course in Metric Geometry}.
\bseries{Graduate Studies in Mathematics}
\bvolume{33}.
\bpublisher{Amer. Math. Soc.},
\blocation{Providence, RI}.
\bid{doi={10.1090/gsm/033}, mr={1835418}}
\end{bbook}
%

\bptok{imsref}%
\endbibitem

\bibitem{Croy1}
\begin{barticle}[mr]
\bauthor{\bsnm{Croydon},~\bfnm{David}\binits{D.}}
(\byear{2008}).
\btitle{Convergence of simple random walks on random discrete trees to {B}rownian motion on the continuum random tree}.
\bjournal{Ann. Inst. Henri Poincar\'e Probab. Stat.}
\bvolume{44}
\bpages{987--1019}.
\bid{doi={10.1214/07-AIHP153}, issn={0246-0203}, mr={2469332}}
\end{barticle}
%

\bptok{imsref}%
\endbibitem

\bibitem{CroyKum}
\begin{barticle}[mr]
\bauthor{\bsnm{Croydon},~\bfnm{David}\binits{D.}} \AND
\bauthor{\bsnm{Kumagai},~\bfnm{Takashi}\binits{T.}}
(\byear{2008}).
\btitle{Random walks on {G}alton--{W}atson trees with infinite variance offspring distribution conditioned to survive}.
\bjournal{Electron. J. Probab.}
\bvolume{13}
\bpages{1419--1441}.
\bid{doi={10.1214/EJP.v13-536}, issn={1083-6489}, mr={2438812}}
\end{barticle}
%

\bptok{imsref}%
\endbibitem

\bibitem{Croy-1}
\begin{barticle}[mr]
\bauthor{\bsnm{Croydon},~\bfnm{D.~A.}\binits{D.~A.}}
(\byear{2007}).
\btitle{Heat kernel fluctuations for a resistance form with non-uniform volume growth}.
\bjournal{Proc. Lond. Math. Soc. (3)}
\bvolume{94}
\bpages{672--694}.
\bid{doi={10.1112/plms/pdl025}, issn={0024-6115}, mr={2325316}}
\end{barticle}
%

\bptok{imsref}%
\endbibitem

\bibitem{Croy0}
\begin{barticle}[mr]
\bauthor{\bsnm{Croydon},~\bfnm{David~A.}\binits{D.~A.}}
(\byear{2008}).
\btitle{Volume growth and heat kernel estimates for the continuum random tree}.
\bjournal{Probab. Theory Related Fields}
\bvolume{140}
\bpages{207--238}.
\bid{doi={10.1007/s00440-007-0063-4}, issn={0178-8051}, mr={2357676}}
\end{barticle}
%

\bptok{imsref}%
\endbibitem

\bibitem{Croy2}
\begin{barticle}[mr]
\bauthor{\bsnm{Croydon},~\bfnm{David~A.}\binits{D.~A.}}
(\byear{2009}).
\btitle{Hausdorff measure of arcs and {B}rownian motion on {B}rownian spatial trees}.
\bjournal{Ann. Probab.}
\bvolume{37}
\bpages{946--978}.
\bid{doi={10.1214/08-AOP425}, issn={0091-1798}, mr={2537546}}
\end{barticle}
%

\bptok{imsref}%
\endbibitem

\bibitem{Croy3}
\begin{barticle}[mr]
\bauthor{\bsnm{Croydon},~\bfnm{D.~A.}\binits{D.~A.}}
(\byear{2010}).
\btitle{Scaling limits for simple random walks on random ordered graph trees}.
\bjournal{Adv. in Appl. Probab.}
\bvolume{42}
\bpages{528--558}.
\bid{doi={10.1239/aap/1275055241}, issn={0001-8678}, mr={2675115}}
\end{barticle}
%

\bptok{imsref}%
\endbibitem

\bibitem{CHK}
\begin{barticle}[mr]
\bauthor{\bsnm{Croydon},~\bfnm{D.~A.}\binits{D.~A.}},
\bauthor{\bsnm{Hambly},~\bfnm{B.~M.}\binits{B.~M.}} \AND
\bauthor{\bsnm{Kumagai},~\bfnm{T.}\binits{T.}}
(\byear{2012}).
\btitle{Convergence of mixing times for sequences of random walks on finite graphs}.
\bjournal{Electron. J. Probab.}
\bvolume{17}
\bpages{1--32}.
\bid{issn={1083-6489}, mr={2869250}}
\end{barticle}
%

\bptok{imsref}%
\endbibitem

\bibitem{DL}
\begin{barticle}[mr]
\bauthor{\bsnm{Duquesne},~\bfnm{Thomas}\binits{T.}} \AND
\bauthor{\bsnm{Le Gall},~\bfnm{Jean-Fran{\c{c}}ois}\binits{J.-F.}}
(\byear{2005}).
\btitle{Probabilistic and fractal aspects of L\'evy trees}.
\bjournal{Probab. Theory Related Fields}
\bvolume{131}
\bpages{553--603}.
\bid{doi={10.1007/s00440-004-0385-4}, issn={0178-8051}, mr={2147221}}
\end{barticle}
%

\bptok{imsref}%
\endbibitem

\bibitem{Edgar}
\begin{bbook}[mr]
\bauthor{\bsnm{Edgar},~\bfnm{Gerald~A.}\binits{G.~A.}}
(\byear{1998}).
\btitle{Integral, Probability, and Fractal Measures}.
\bpublisher{Springer},
\blocation{New York}.
\bid{doi={10.1007/978-1-4757-2958-0}, mr={1484412}}
\end{bbook}
%

\bptok{imsref}%
\endbibitem

\bibitem{Evans}
\begin{bbook}[mr]
\bauthor{\bsnm{Evans},~\bfnm{Steven~N.}\binits{S.~N.}}
(\byear{2008}).
\btitle{Probability and Real Trees}.
\bseries{Lecture Notes in Math.}
\bvolume{1920}.
\bpublisher{Springer},
\blocation{Berlin}.
\bid{doi={10.1007/978-3-540-74798-7}, mr={2351587}}
\end{bbook}
%

\bptok{imsref}%
\endbibitem

\bibitem{EPW}
\begin{barticle}[mr]
\bauthor{\bsnm{Evans},~\bfnm{Steven~N.}\binits{S.~N.}},
\bauthor{\bsnm{Pitman},~\bfnm{Jim}\binits{J.}} \AND
\bauthor{\bsnm{Winter},~\bfnm{Anita}\binits{A.}}
(\byear{2006}).
\btitle{Rayleigh processes, real trees, and root growth with re-grafting}.
\bjournal{Probab. Theory Related Fields}
\bvolume{134}
\bpages{81--126}.
\bid{doi={10.1007/s00440-004-0411-6}, issn={0178-8051}, mr={2221786}}
\end{barticle}
%

\bptok{imsref}%
\endbibitem

\bibitem{FOT}
\begin{bbook}[mr]
\bauthor{\bsnm{Fukushima},~\bfnm{Masatoshi}\binits{M.}},
\bauthor{\bsnm{{O}shima},~\bfnm{Y{o}ichi}\binits{Y.}} \AND
\bauthor{\bsnm{Takeda},~\bfnm{Masayoshi}\binits{M.}}
(\byear{1994}).
\btitle{Dirichlet Forms and Symmetric {M}arkov Processes}.
\bseries{De Gruyter Studies in Mathematics}
\bvolume{19}.
\bpublisher{de Gruyter},
\blocation{Berlin}.
\bid{doi={10.1515/9783110889741}, mr={1303354}}
\end{bbook}
%

\bptok{imsref}%
\endbibitem

\bibitem{HM}
\begin{barticle}[mr]
\bauthor{\bsnm{Haas},~\bfnm{B{\'e}n{\'e}dicte}\binits{B.}} \AND
\bauthor{\bsnm{Miermont},~\bfnm{Gr{\'e}gory}\binits{G.}}
(\byear{2012}).
\btitle{Scaling limits of {M}arkov branching trees with applications to {G}alton--{W}atson and random unordered trees}.
\bjournal{Ann. Probab.}
\bvolume{40}
\bpages{2589--2666}.
\bid{doi={10.1214/11-AOP686}, issn={0091-1798}, mr={3050512}}
\end{barticle}
%

\bptok{imsref}%
\endbibitem

\bibitem{Haggstrom}
\begin{barticle}[mr]
\bauthor{\bsnm{H{\"a}ggstr{\"o}m},~\bfnm{Olle}\binits{O.}}
(\byear{1995}).
\btitle{Random-cluster measures and uniform spanning trees}.
\bjournal{Stochastic Process. Appl.}
\bvolume{59}
\bpages{267--275}.
\bid{doi={10.1016/0304-4149(95)00042-6}, issn={0304-4149}, mr={1357655}}
\end{barticle}
%

\bptok{imsref}%
\endbibitem

\bibitem{Kall}
\begin{bbook}[mr]
\bauthor{\bsnm{Kallenberg},~\bfnm{Olav}\binits{O.}}
(\byear{2002}).
\btitle{Foundations of Modern Probability},
\bedition{2nd} ed.
\bpublisher{Springer},
\blocation{New York}.
\bid{doi={10.1007/978-1-4757-4015-8}, mr={1876169}}
\end{bbook}
%

\bptok{imsref}%
\endbibitem

\bibitem{Kenyon}
\begin{barticle}[mr]
\bauthor{\bsnm{Kenyon},~\bfnm{Richard}\binits{R.}}
(\byear{2000}).
\btitle{The asymptotic determinant of the discrete {L}aplacian}.
\bjournal{Acta Math.}
\bvolume{185}
\bpages{239--286}.
\bid{doi={10.1007/BF02392811}, issn={0001-5962}, mr={1819995}}
\end{barticle}
%

\bptok{imsref}%
\endbibitem

\bibitem{Kigami}
\begin{barticle}[mr]
\bauthor{\bsnm{Kigami},~\bfnm{Jun}\binits{J.}}
(\byear{1995}).
\btitle{Harmonic calculus on limits of networks and its application to dendrites}.
\bjournal{J. Funct. Anal.}
\bvolume{128}
\bpages{48--86}.
\bid{doi={10.1006/jfan.1995.1023}, issn={0022-1236}, mr={1317710}}
\end{barticle}
%

\bptok{imsref}%
\endbibitem

\bibitem{Kigami2}
\begin{barticle}[mr]
\bauthor{\bsnm{Kigami},~\bfnm{Jun}\binits{J.}}
(\byear{2012}).
\btitle{Resistance forms, quasisymmetric maps and heat kernel estimates}.
\bjournal{Mem. Amer. Math. Soc.}
\bvolume{216}
\bpages{vi+132}.
\bid{doi={10.1090/S0065-9266-2011-00632-5}, issn={0065-9266}, mr={2919892}}
\end{barticle}
%

\bptok{imsref}%
\endbibitem

\bibitem{Kirchoff}
\begin{barticle}[auto:parserefs-M02]
\bauthor{\bsnm{Kirchhoff},~\bfnm{G.}\binits{G.}}
(\byear{1847}).
\btitle{Ueber die Aufl\"osung der Gleichungen, auf welche man bei der Untersuchung der linearen Vertheilung galvanischer Str\"ome gef\"uhrt wird}.
\bjournal{Ann. Phys.}
\bvolume{148}
\bpages{497--508}.
\end{barticle}
%

\bptok{imsref}%
\endbibitem

\bibitem{K04}
\begin{barticle}[mr]
\bauthor{\bsnm{Kumagai},~\bfnm{Takashi}\binits{T.}}
(\byear{2004}).
\btitle{Heat kernel estimates and parabolic {H}arnack inequalities on graphs and resistance forms}.
\bjournal{Publ. Res. Inst. Math. Sci.}
\bvolume{40}
\bpages{793--818}.
\bid{issn={0034-5318}, mr={2074701}}
\end{barticle}
%

\bptok{imsref}%
\endbibitem

\bibitem{K14}
\begin{bbook}[mr]
\bauthor{\bsnm{Kumagai},~\bfnm{Takashi}\binits{T.}}
(\byear{2014}).
\btitle{Random Walks on Disordered Media and Their Scaling Limits}.
\bseries{Lecture Notes in Math.}
\bvolume{2101}.
\bpublisher{Springer},
\blocation{Cham}.
\bid{doi={10.1007/978-3-319-03152-1}, mr={3156983}}
\end{bbook}
%

\bptok{imsref}%
\endbibitem

\bibitem{KM}
\begin{barticle}[mr]
\bauthor{\bsnm{Kumagai},~\bfnm{Takashi}\binits{T.}} \AND
\bauthor{\bsnm{Misumi},~\bfnm{Jun}\binits{J.}}
(\byear{2008}).
\btitle{Heat kernel estimates for strongly recurrent random walk on random media}.
\bjournal{J. Theoret. Probab.}
\bvolume{21}
\bpages{910--935}.
\bid{doi={10.1007/s10959-008-0183-5}, issn={0894-9840}, mr={2443641}}
\end{barticle}
%

\bptok{imsref}%
\endbibitem

\bibitem{Law99}
\begin{bincollection}[mr]
\bauthor{\bsnm{Lawler},~\bfnm{Gregory~F.}\binits{G.~F.}}
(\byear{1999}).
\btitle{Loop-erased random walk}.
In \bbooktitle{Perplexing Problems in Probability}
(\beditor{\binits{M.}\bfnm{M.}~\bsnm{Bramson}}
\AND
\beditor{\binits{R.}\bfnm{R.} \bsnm{Durrett}}, eds.).
\bseries{Progress in Probability}
\bvolume{44}
\bpages{197--217}.
\bpublisher{Birkh\"auser},
\blocation{Boston, MA}.
\bid{mr={1703133}}
\bptnote{check pages}%
\end{bincollection}
%

\bptok{imsref}%
\endbibitem

\bibitem{Lawlerbook}
\begin{bbook}[mr]
\bauthor{\bsnm{Lawler},~\bfnm{Gregory~F.}\binits{G.~F.}}
(\byear{2013}).
\btitle{Intersections of Random Walks}.
\bpublisher{Birkh\"auser/Springer},
\blocation{New York}.
\bid{doi={10.1007/978-1-4614-5972-9}, mr={2985195}}
\bptnote{check year}%
\end{bbook}
%

\bptok{imsref}%
\endbibitem

\bibitem{LE}
\begin{barticle}[mr]
\bauthor{\bsnm{Lawler},~\bfnm{Gregory~F.}\binits{G.~F.}}
(\byear{2014}).
\btitle{The probability that planar loop-erased random walk uses a given edge}.
\bjournal{Electron. Commun. Probab.}
\bvolume{19}
\bpages{1--13}.
\bid{doi={10.1214/ECP.v19-2908}, issn={1083-589X}, mr={3246970}}
\bptnote{check pages}%
\end{barticle}
%

\bptok{imsref}%
\endbibitem

\bibitem{LSW}
\begin{barticle}[mr]
\bauthor{\bsnm{Lawler},~\bfnm{Gregory~F.}\binits{G.~F.}},
\bauthor{\bsnm{Schramm},~\bfnm{Oded}\binits{O.}} \AND
\bauthor{\bsnm{Werner},~\bfnm{Wendelin}\binits{W.}}
(\byear{2004}).
\btitle{Conformal invariance of planar loop-erased random walks and uniform spanning trees}.
\bjournal{Ann. Probab.}
\bvolume{32}
\bpages{939--995}.
\bid{doi={10.1214/aop/1079021469}, issn={0091-1798}, mr={2044671}}
\end{barticle}
%

\bptok{imsref}%
\endbibitem

\bibitem{LS}
\begin{barticle}[mr]
\bauthor{\bsnm{Lawler},~\bfnm{Gregory~F.}\binits{G.~F.}} \AND
\bauthor{\bsnm{Sheffield},~\bfnm{Scott}\binits{S.}}
(\byear{2011}).
\btitle{A natural parametrization for the {S}chramm--{L}oewner evolution}.
\bjournal{Ann. Probab.}
\bvolume{39}
\bpages{1896--1937}.
\bid{doi={10.1214/10-AOP560}, issn={0091-1798}, mr={2884877}}
\end{barticle}
%

\bptok{imsref}%
\endbibitem

\bibitem{LZ}
\begin{barticle}[mr]
\bauthor{\bsnm{Lawler},~\bfnm{Gregory~F.}\binits{G.~F.}} \AND
\bauthor{\bsnm{Zhou},~\bfnm{Wang}\binits{W.}}
(\byear{2013}).
\btitle{{\textit{SLE}} curves and natural parametrization}.
\bjournal{Ann. Probab.}
\bvolume{41}
\bpages{1556--1584}.
\bid{doi={10.1214/12-AOP742}, issn={0091-1798}, mr={3098684}}
\end{barticle}
%

\bptok{imsref}%
\endbibitem

\bibitem{rrt}
\begin{barticle}[mr]
\bauthor{\bsnm{Le Gall},~\bfnm{Jean-Fran{\c{c}}ois}\binits{J.-F.}}
(\byear{2006}).
\btitle{Random real trees}.
\bjournal{Ann. Fac. Sci. Toulouse Math. (6)}
\bvolume{15}
\bpages{35--62}.
\bid{issn={0240-2963}, mr={2225746}}
\end{barticle}
%

\bptok{imsref}%
\endbibitem

\bibitem{Mas09}
\begin{barticle}[mr]
\bauthor{\bsnm{Masson},~\bfnm{Robert}\binits{R.}}
(\byear{2009}).
\btitle{The growth exponent for planar loop-erased random walk}.
\bjournal{Electron. J. Probab.}
\bvolume{14}
\bpages{1012--1073}.
\bid{doi={10.1214/EJP.v14-651}, issn={1083-6489}, mr={2506124}}
\end{barticle}
%

\bptok{imsref}%
\endbibitem

\bibitem{Miermont}
\begin{barticle}[mr]
\bauthor{\bsnm{Miermont},~\bfnm{Gr{\'e}gory}\binits{G.}}
(\byear{2009}).
\btitle{Tessellations of random maps of arbitrary genus}.
\bjournal{Ann. Sci. \'Ec. Norm. Sup\'er. (4)}
\bvolume{42}
\bpages{725--781}.
\bid{issn={0012-9593}, mr={2571957}}
\end{barticle}
%

\bptok{imsref}%
\endbibitem

\bibitem{Pemantle}
\begin{barticle}[mr]
\bauthor{\bsnm{Pemantle},~\bfnm{Robin}\binits{R.}}
(\byear{1991}).
\btitle{Choosing a spanning tree for the integer lattice uniformly}.
\bjournal{Ann. Probab.}
\bvolume{19}
\bpages{1559--1574}.
\bid{issn={0091-1798}, mr={1127715}}
\end{barticle}
%

\bptok{imsref}%
\endbibitem

\bibitem{Schramm}
\begin{barticle}[mr]
\bauthor{\bsnm{Schramm},~\bfnm{Oded}\binits{O.}}
(\byear{2000}).
\btitle{Scaling limits of loop-erased random walks and uniform spanning trees}.
\bjournal{Israel J. Math.}
\bvolume{118}
\bpages{221--288}.
\bid{doi={10.1007/BF02803524}, issn={0021-2172}, mr={1776084}}
\end{barticle}
%

\bptok{imsref}%
\endbibitem

\bibitem{Wilson}
\begin{binproceedings}[mr]
\bauthor{\bsnm{Wilson},~\bfnm{David~Bruce}\binits{D.~B.}}
(\byear{1996}).
\btitle{Generating random spanning trees more quickly than the cover time}.
In \bbooktitle{Proceedings of the {T}wenty-Eighth {A}nnual ACM {S}ymposium on the {T}heory of {C}omputing ({P}hiladelphia, PA, 1996)}
\bpages{296--303}.
\bpublisher{ACM},
\blocation{New York}.
\bid{doi={10.1145/237814.237880}, mr={1427525}}
\end{binproceedings}
%

\bptok{imsref}%
\endbibitem
\end{thebibliography}
\end{document}